\documentclass[a4paper]{report}
\usepackage[latin1]{inputenc}
\usepackage{amsfonts, amsmath, amssymb, amsthm, fancyhdr}
\usepackage[left=3cm,right=3cm,top=3cm,bottom=2cm,includeheadfoot]{geometry}
\usepackage[arrow, matrix, curve]{xy}
\usepackage{makeidx}

\newtheorem{theo}{Theorem}[section]
\newtheorem*{theonn}{Theorem}
\newtheorem{lemm}[theo]{Lemma}
\newtheorem{prop}[theo]{Proposition}
\newtheorem{coro}[theo]{Corollary}
\theoremstyle{definition} 
\newtheorem{defi}[theo]{Definition}
\newtheorem{cons}[theo]{Construction}
\newtheorem{exam}[theo]{Example}

\newtheorem{rema}[theo]{Remark}

\renewcommand{\theenumi}{(\roman{enumi})} 
\renewcommand{\labelenumi}{\theenumi}

\DeclareMathOperator{\End}{End}
\DeclareMathOperator{\GL}{GL}

\DeclareMathOperator{\ad}{ad}
\DeclareMathOperator{\Ad}{Ad}
\DeclareMathOperator{\im}{im}
\DeclareMathOperator{\Hom}{Hom}
\DeclareMathOperator{\Id}{Id}
\DeclareMathOperator{\Ham}{Ham}
\DeclareMathOperator{\Symp}{Symp}
\DeclareMathOperator{\supp}{supp}
\DeclareMathOperator{\Ver}{Ver}
\DeclareMathOperator{\Hor}{Hor}
\DeclareMathOperator{\Hol}{Hol}

\makeindex

\begin{document}
	
	\title{The Orbit Method for Compact Connected Lie Groups}
	\author{Matthias Peter}
	\date{Paderborn, October 2007\\ (Today is \today)}
			
	\begin{center} \Large
		Paderborn University\\
		Department of Mathematics
	\end{center}
	
	\vspace{0.5cm}
	
	\begin{center}
		
	\end{center}
	
	\vspace{0.5cm}
	
	\begin{center}
		{\Large Diploma Thesis}
		
		\vspace{0.5cm}
		
		{\huge \textbf{The Orbit Method for Compact Connected Lie Groups}}
		
		\vspace{0.5cm}
		
		{\Large Matthias Peter}

		October 2007
	\end{center}
	
	\vspace{0.5cm}
	
	\begin{center}
		
	\end{center}
	
	\vspace{4cm}
	
	\begin{center} \large
		\begin{tabular}{rl} 
			Advisor: & Prof. Dr. Joachim Hilgert \\
								& Department of Mathematics	\\
			Co-Advisor:  & Prof. Dr. Sönke Hansen \\
									 & Department of Mathematics \\
		\end{tabular}		
	\end{center}
	
	\begin{center}
		(Revised version of \today)
	\end{center}
	
	\thispagestyle{empty}
	
	\newpage
		\thispagestyle{empty}
		{I would very much like to thank the people at the Department of Mathematics for five great years in Paderborn. Most of all, I would like to thank both Prof. Joachim Hilgert and Prof. Sönke Hansen who had a major influence on my mathematical development and who always had time for me and my questions. I am especially indebted to Prof. Hilgert for his guidance. He encouraged me to spend a semester abroad and proposed the topic of this thesis to me. I am also grateful to the "Studienstiftung des deutschen Volkes" for two brilliant summer academies and the opportunity to meet a lot of great people. Last but not least, I would like to thank my family and my friends for their unrelenting support and encouragement.}
		
		\vfill
		With the following declaration I affirm to have written this thesis on my own and using only the resources listed in the bibliography.
		
		Hiermit versichere ich die vorliegende Arbeit selbstständig und nur unter Verwendung der angegebenen Quellen und Hilfsmittel verfasst zu haben.
		
		Paderborn, den 31.10.2007
		\vspace{2 cm}

	\newpage
	\pagenumbering{roman}
  
	\tableofcontents
	\vspace{0.5cm}
	\textbf{Bibliography \hfill \pageref{biblio}}
	
	\vspace{0.5cm}
	\textbf{Index \hfill 55} 
	
	
		
		
	\newpage
	\pagenumbering{arabic}

	\begin{center}
	{\Large \textbf{Introduction}}
	\end{center}
	
	The goal of this diploma thesis is to give a detailed description of Kirillov's Orbit Method for the case of compact connected Lie groups. The theory of Kirillov aims at finding all irreducible unitary representations of a given Lie group $G$. By a representation we mean a homomorphism from $G$ to the group of linear transformations of a vector space. The classical result is the following (cf. \cite[Theorem 10.1]{vogur}).
	\begin{theonn}[Kirillov]
			Suppose $G$ is a connected, simply connected nilpotent Lie group with Lie algebra $\mathfrak g$. Then the irreducible unitary representations of $G$ are in a natural one-to-one correspondence with the set of orbits of $G$ in $\mathfrak g^*$.
		\end{theonn}
	For more complicated Lie groups this statement is no longer true. One has to make amendments such as restricting to a subset of the set of orbits or trying to attach multiple representations to a single orbit (cf. \cite[§10]{vogur}). Though people have worked on this topic for decades there is still no generalization which covers all Lie groups or even the case of reductive groups.
	
	Nevertheless there are techniques which work for certain subclasses such as, for example, compact connected Lie groups. Here one can see what made the Orbit Method so interesting in the first place: its geometric cleanness which inspires hope that there is still a deep reason to be found why it works. At the end of this thesis we will have seen that the Orbit Method for compact connected Lie groups can be directly related to the classical Borel-Weil approach. This also means that it yields no new results, but it is a nice example for this beautiful theory and will, perhaps, help to understand representation theory for other classes of groups.
	
	The first chapter is intended to recall some facts about Lie groups. The most important result is that for a functional $\lambda$ on the Lie algebra $\mathfrak t$ of a maximal torus in a compact connected Lie group the stabilizer is of the form $C(T_1)$ for a smaller torus $T_1$. This allows us to define a complex structure on the coadjoint orbit $\mathcal O_\lambda = G/C(T_1)$. In the second chapter we present the Borel-Weil Theorem which tells us that all irreducible unitary representations of a compact connected Lie group can be realized in the space of holomorphic sections of the bundle $G\times_\tau\mathbb C$, where $\tau$ is a character of a maximal torus. Finally, in the third chapter we show that coadjoint orbits are symplectic manifolds for which prequantization yields a complex line bundle. Then we investigate for which classes of coadjoint orbits one obtains representations in the space of polarized sections of the prequantum bundle. This will lead us back to the result of Borel-Weil.
	
	\chapter{Structure Theory of Compact Lie Groups}
	
	\section{Fundamentals}
	
	\subsection{Smooth Manifolds}
		
		This section is not intended as an introduction to differentiable manifolds. We merely introduce notations and some tools. For a comprehensive introduction into this topic consult, for example, \cite{nara}, \cite{war} or \cite{hga}.
		
		\begin{defi} 
			Let $M$ be a second countable Hausdorff topological space. We say that $M$ is a \textbf{smooth manifold}\index{manifold!smooth} of dimension $n$ if there is given a family of pairs $\{U_j,\varphi_j\}_{j\in J}$, where $U_j$ is an open set in $M$ and $\varphi_j:U_j\rightarrow\varphi(U_j)$ is a homeomorphism of $U_j$ onto an open set in $\mathbb R^n$ such that
			\renewcommand{\theenumi}{(\roman{enumi})}
			\renewcommand{\labelenumi}{\theenumi}
			\begin{enumerate}
				\item $\bigcup_{j\in J}U_j=M$; and
				\item whenever $j,\ k\in J$ such that $U_j\cap U_k\neq\emptyset$, then
					\begin{equation} \label{eq:atlas}
						\varphi_k\circ\varphi_j^{-1}\colon \varphi_j(U_j\cap U_k)\to\varphi_k(U_j\cap U_k)
					\end{equation}
				is a smooth map.
			\end{enumerate}
			By \textbf{smooth} maps, we always mean infinitely differentiable maps.
		\end{defi}
		
		The family $\mathcal A=\{U_j,\varphi_j\}_{j\in J}$ is called an \textbf{atlas} of $M$ and its elements are called \textbf{charts}. If $\mathcal B=\{V_k,\psi_k\}_{k\in K}$ is another atlas of $M$, then $\mathcal A$ and $\mathcal B$ are said to be \textbf{compatible}, if $\mathcal A\cup\mathcal B$ satisfies (\ref{eq:atlas}). Then $\mathcal A$ and $\mathcal B$ are contained in the larger atlas $\mathcal A\cup\mathcal B$. An atlas which is maximal with respect to inclusion is called a \textbf{differentiable structure}\index{differentiable structure}. One can show (cf. \cite[p.5]{war}) that each atlas of $M$ uniquely determines a differentiable structure on $M$.
		
		\begin{defi}
			Let $M$ be a smooth manifold with atlas $\{U_j,\varphi_j\}_{j\in J}$. Let $f\colon M\to\mathbb C$ be a continuous function. We say that $f$ is \textbf{smooth}\index{function!smooth} if for each chart $(U_j,\varphi_j)$ the map
			$$f\circ\varphi_j^{-1}\colon\varphi_j(U_j)\to\mathbb R$$
			is smooth. Similarly, if $N$ is another smooth manifold with atlas $\{V_k,\psi_k\}_{k\in K}$, then a continuous map $f\colon M\to N$ is called smooth\index{map!smooth} if for all $p\in M$ there are $j,\ k$ such that $p\in U_j$, $f(p)\in V_k$ and the map
			$$\psi_k\circ f\circ\varphi_j^{-1}\colon\varphi_j(U_j)\to\psi_k(V_k)$$
			is smooth. The vector space of smooth real- or complex-valued functions on M is denoted by $C^\infty(M,\mathbb R)$ and $C^\infty(M,\mathbb C)$ respectively. The \textbf{support} $\supp f$\index{$\supp$, support} of $C^\infty(M,\mathbb C)$ is defined as the closure of $M\backslash f^{-1}(\{0\})$. The group of diffeomorphisms on $M$ is denoted by $D(M)$\index{$D(M)$, group of diffeomorphisms}.
		\end{defi}
		
		\begin{defi} \label{defi:coverings} 
			A family $\{U_j\}_{j\in J}$ of subsets of $M$ is called a \textbf{covering}\index{covering} of $M$ if the union of all $U_j$ is $M$. It is called an \textbf{open covering}\index{covering!open} if each $U_j$ is open. A \textbf{refinement}\index{refinement} $\{V_k\}_{k\in K}$ of $\{U_j\}_{j\in J}$ is a covering such that for each $k$ there is an $j$ such that $V_k\subseteq U_j$. A family $\{U_j\}_{j\in J}$ of subsets of $M$ is called \textbf{locally finite}\index{covering!locally finite} if each point $p\in M$ has a neighborhood $V_p$ such that $V_p\cap U_j\neq\emptyset$ for only finitely many $j\in J$. A topological space is \textbf{paracompact}\index{paracompact}, if every open covering has an open locally finite refinement. A covering $\{U_j\}_{j\in J}$ is called \textbf{contractible}\index{covering!contractible} if any non-empty intersection of finitely many $U_j$ is contractible.
		\end{defi}
		
		\begin{defi}
			A \textbf{partition of unity}\index{partition of unity} on $M$ is a family $\{f_j\}_{j\in J}\subseteq C^\infty(M,\mathbb R)$ such that
			\renewcommand{\theenumi}{(\roman{enumi})}
			\renewcommand{\labelenumi}{\theenumi}
			\begin{enumerate}
				\item $f_j\ge 0\quad\forall j\in J$,
				\item $\{\supp f_j\}_{j\in J}$ is locally finite; and
				\item $\sum_{j\in J}f_j(p)=1\quad\forall p\in M$.
			\end{enumerate}
			It is called \textbf{subordinate} to the covering $\{V_k\}_{k\in K}$ if for each $j$ there exists a $k$ such that $\supp f_j\subseteq V_k$.
		\end{defi}
		
		\begin{theo} \label{theo:countpara}  
			Let $M$ be a smooth manifold. According to our definition $M$ is a second countable Hausdorff space. Then $M$ is paracompact, every open covering of $M$ has a subordinate partition of unity and any open covering has a contractible open refinement.
		\end{theo}
		
		\begin{proof}
			\cite[Lemma 1.9]{war}, \cite[Theorem 2.2.14]{nara} and \cite[Theorem 1.2]{walhomo}.
		\end{proof}
		
		For a verification of well-definedness and an extensive introduction to the following concepts, please confer \cite{hga}, \cite{nara} or \cite{war}.

		Let $M$ be a smooth manifold of dimension $n$ and $p\in M$. On the set of pairs $(f,U)$, where $f$ is a smooth real-valued function which is defined on the open neighborhood $U$ of $p$, we introduce the following equivalence relation.
		$$(f,U)\simeq(g,V)\quad\Leftrightarrow\quad f_{|U\cap V}=g_{|U\cap V}$$
		We denote the equivalence class of $(f,U)$ by $f_p$ and we denote the set of equivalence classes by $C_p^\infty$. We call $C_p^\infty$ the set of  \textbf{smooth germs}\index{germ} in $p$ and observe that it carries a natural algebra structure inherited by addition and multiplication of smooth functions. The \textbf{tangent space}\index{tangent space} $T_pM$\index{$T_pM$, tangent space} of $M$ at $p$ is the $n$-dimensional real vector space of derivations of $C_p^\infty$. By a derivation\index{derivation} we mean a real-linear map $v\colon C_p^\infty\to\mathbb R$ such that
		$$v(f_p\cdot g_p)=g(p)\cdot v(f_p)+f(p)\cdot v(g_p)\quad\forall f_p,\ g_p\in C^\infty_p.$$
		If $N$ is another smooth manifold and  $\tau\colon M\to N$ is a smooth map, then the \textbf{derivative}\index{derivative} of $\tau$ at $p\in M$ is the linear map
		$$\tau^\prime(p)\colon T_pM\to T_{\tau(p)}M,\quad v\mapsto \tau^\prime(p)v,$$
		where the latter is defined as
		$$(\tau^\prime(p)v)f_p:=v(f\circ\tau)_{\tau(p)}\quad\forall f_{\tau(p)}\in C_{f(p)}^\infty.$$
		The sets $TM:=\bigcup_{p\in M}T_pM$ and $TM_{\mathbb C}:=\bigcup_{p\in M}(T_pM)_\mathbb C$ inherit a natural differentiable structure from $M$ and are called the \textbf{tangent bundle}\index{tangent bundle} and the \textbf{complexified tangent bundle}. More generally one can define the \textbf{tensor bundle} of type $(r,s)$ denoted by $\bigotimes^s_r TM:=\bigcup_{p\in M}\bigotimes^s_r T_pM$ which is again a smooth manifold.
		
		A \textbf{vector field}\index{vector field} $\mathfrak X$ is a smooth map $\mathfrak X\colon M\to TM$ such that $\pi\circ\mathfrak X=\Id_M$. Here $\pi\colon TM\to M$ is the natural projection which associates to $v\in T_pM$ its base point $p$. The $C^\infty(M,\mathbb R)$-module of smooth vector fields is denoted by $\mathcal X(M)$\index{$\mathcal X(M)$, module of tensor fields}. We will interpret $\mathcal X(M)$ as the derivations of the algebra $C^\infty(M,\mathbb R)$. By $\mathcal X_\mathbb C(M)$ we mean the complexification of $\mathcal X(M)$ which we interpret as derivations of $C^\infty(M,\mathbb C)$. The \textbf{bracket}\index{bracket of vector fields} $[\mathfrak X,\mathfrak Y]$ of two vector fields $\mathfrak X,\ \mathfrak Y$ is defined as the commutator $\mathfrak X\circ\mathfrak Y-\mathfrak Y\circ\mathfrak X$. In an analogous manner, a \textbf{tensor field}\index{tensor field} $\mathfrak T$ is a smooth map $\mathfrak T\colon M\to\bigotimes^s_r TM$ with $\pi\circ\mathfrak T=\Id_M$. We denote by $\Omega^k(M)$\index{$\Omega^k(M)$, module of $k$-forms}, $k\in\mathbb N_0$, the $C^\infty(M,\mathbb R)$-module of smooth differential $k$-forms. A \textbf{differential $k$-form} or \textbf{$k$-form}\index{form} is an alternating $C^\infty(M,\mathbb R)$-$k$-linear map from $\mathcal X(M)$ to $C^\infty(M,\mathbb R)$. We write $d\colon\Omega^k(M)\to\Omega^{k+1}(M)$ for the \textbf{exterior derivative}\index{exterior derivative}\index{$d$, exterior derivative} and recall its definition for the cases $k=0$ and $k=1$:
		
		\begin{defi}
			Let $\mathfrak X,\ \mathfrak Y\in\mathcal X(M)$, $f\in\Omega^0(M)=C^\infty(M)$ and $\alpha\in\Omega^1(M)$. Then
			\begin{align*}
				(df)(\mathfrak X) & =\mathfrak X(f)\\
				(d\alpha)(\mathfrak X,\mathfrak Y) & = \mathfrak X(\alpha(\mathfrak Y))-\mathfrak Y(\alpha(\mathfrak X))-\alpha([\mathfrak X,\mathfrak Y]).
			\end{align*}
		A form $\omega\in\Omega^k(M)$ is called \textbf{closed}\index{form! closed} if $d\omega=0$. It is called \textbf{exact}\index{form! exact} if there is a $(k-1)$-form $\alpha$ such that $d\alpha=\omega$. Exact forms are always closed because $d^2=0$. The \textbf{interior product}\index{interior product} or \textbf{contraction}\index{contraction}\index{$i_\mathfrak X\omega$, contraction} of $\omega$ with respect to $\mathfrak X$ is defined as the $(k-1)$-form
			$$(i_\mathfrak X\omega)(\mathfrak X_1,\dots,\mathfrak X_{k-1})=\omega(\mathfrak X,\mathfrak X_1,\dots,\mathfrak X_{k-1})\quad \forall \mathfrak X_1,\dots,\mathfrak X_{k-1}\in\mathcal X(M).$$
			We define the \textbf{Lie derivative}\index{Lie derivative}\index{$\mathcal L_{\mathfrak X}$, Lie derivative} of $\mathfrak Y$ and the Lie derivative of $\omega$ with respect to $\mathfrak X$ by the formulae
			\begin{align*}
				\mathcal L_{\mathfrak X}(\mathfrak Y)(p) & = \lim_{t\rightarrow 0}\frac{1}{t}((\Phi_{\mathfrak X,t})^*(\mathfrak Y)-\mathfrak Y)(p)\quad \forall p\in M, \\
				\mathcal L_{\mathfrak X}(\omega)(p) & = \lim_{t\rightarrow 0}\frac{1}{t}( (\Phi_{\mathfrak X,t})^*\omega-\omega )(p)\quad\forall p\in M,
			\end{align*}
			where $\Phi_{\mathfrak X,t}$ is the \textbf{flow} of $\mathfrak X$.
		\end{defi}
		
		\begin{theo}[Cartan Identity]\index{Cartan Identity}
			The three operations are related by the following formula.
			$$\mathcal L_\mathfrak X\omega=i_\mathfrak X(d\omega)+d(i_\mathfrak X\omega)\quad\forall\mathfrak X\in\mathcal X(M)\ \forall\omega\in\Omega^k(M),\ k\in\mathbb N.$$
		\end{theo}
		
		\begin{proof}
			\cite[§2]{war}.
		\end{proof}
		
		We briefly recall the definition of the de Rham cohomology. For a broad introduction confer \cite[5.28]{war}.
		
		\begin{defi}
			Let $M$ be a smooth manifold. The \textbf{de Rham cohomology}\index{de Rham cohomology}\index{$H^k_{dR}(M,\mathbb R)$, de Rham cohomology} is defined as
			$$H^k_{dR}(M,\mathbb R)=\ker d_{k+1}/\im d_k, k\in\mathbb N_0,$$
			where $d_k\colon\Omega^k(M)\to\Omega^{k+1}(M)$ denotes the exterior derivative.
		\end{defi}
		
		\begin{theo}[Poincaré Lemma]\index{Poincaré Lemma}
			Let $M$ be a smooth manifold, $U\subseteq M$ an open contractible set, $k\in\mathbb N$, and $\omega\in\Omega^k(U)$. Then there exists $\alpha\in\Omega^{k-1}(U)$ such that $d\alpha=\omega$.
		\end{theo}
		
		\begin{proof}
			\cite[2.13.1]{nara}.
		\end{proof}
		
	\subsection{Lie Groups and their Lie Algebras}
		
		\begin{defi}
			Let $\mathbb K$ denote the real or complex numbers. A $\mathbb K$ vector space $\mathfrak g$ together with a $\mathbb K$-bilinear mapping $[.,.]\colon\mathfrak g\times\mathfrak g\to\mathfrak g$ is called a \textbf{Lie algebra}\index{Lie algebra} over $\mathbb K$, if
			\renewcommand{\theenumi}{(\roman{enumi})}
			\renewcommand{\labelenumi}{\theenumi}
			\begin{enumerate}
				\item $[X,Y]=-[Y,X]$\quad\textit{(anti-commutativity)}
				\item $[X, [Y,Z]] + [Y, [Z,X]] + [Z, [X, Y]] = 0$\quad\textit{(Jacobi identity)}
			\end{enumerate}
			The operation $[.,.]$ is called the \textbf{Lie bracket}\index{Lie bracket}\index{$[.,.]$, Lie bracket} of $\mathfrak g$.
		\end{defi}
		
		\begin{defi} \label{def:liegroup}
			A \textbf{Lie group}\index{Lie group} $G$ is a smooth manifold which also carries a group structure such that multiplication and inversion
			$$m\colon G\times G\to G,\quad (g,h)\mapsto gh,\quad inv\colon G\to G,\quad g\mapsto g^{-1}$$
			are smooth maps. The neutral element of $G$ is denoted by $e$ or $e_G$.
		\end{defi}
		
		On a Lie group $G$ one can define for any $g\in G$ a \textbf{left translation}\index{$\lambda_g$, left translation by $g$} $\lambda_g\colon G\to G, \quad h\mapsto gh$. In the same way, one defines the \textbf{right translation} $\rho_g$\index{$\rho_g$, right translation by $g$}. By Definition \ref{def:liegroup} the map $\lambda_g$ is a smooth self-map of $G$ for each $g\in G$. We call a vector field $\mathcal X(G)$ \textbf{left-invariant}\index{vector field!left-invariant}, if
		\begin{equation} \label{eq:leftinvariant}
			\mathfrak X(\lambda_g(h))=\lambda_g^\prime(h)(\mathfrak X(h))\quad\forall g,\ h\in G.
		\end{equation}
		The set of left-invariant vector fields is denoted by $\mathcal X_l(G)$. We see from (\ref{eq:leftinvariant}) that a left-invariant vector field is uniquely determined by its value at $e$. On the other hand, for every $X\in T_eG$ we can define a left-invariant vector field $\widetilde{X}$\index{$\widetilde{X}$, left-invariant vector field} by putting
		$$\widetilde{X}(g)=\lambda_g^\prime(e)X\quad\forall g\in G.$$
		The resulting bijection between $\mathcal X_l(G)$ and $T_eG$ gives rise to the following definition.
		
		\begin{defi} \label{defi:liealg}
			Let $G$ be a Lie group. The \textbf{associated Lie algebra}\index{Lie algebra!associated to a Lie group} $\mathfrak g$\index{$\mathfrak g$, Lie algebra of $G$} is the vector space $T_eG$ together with the bilinear form
			$$[,.,]\colon\mathfrak g\times\mathfrak g\to\mathfrak g, \quad [X,Y]:=[\widetilde{X},\widetilde{Y}](e),$$
			where we have used the bracket of vector fields.
		\end{defi}
		
		\begin{prop}
			$\mathcal X_l(G)$ is closed under the bracket operation and $\mathfrak g$ is a Lie algebra.
		\end{prop}
		
		\begin{proof}
			\cite[Proposition 3.7]{war}.
		\end{proof}
		
		\begin{defi}
			Let $G$ and $H$ be two Lie groups with Lie algebras $\mathfrak g$ and $\mathfrak h$. A \textbf{Lie group homomorphism}\index{Lie group!homomorphism} from $G$ to $H$ is a smooth map between $G$ and $H$ which is also a group homomorphism. By a \textbf{Lie algebra homomorphism}\index{Lie algebra!homomorphism} from $\mathfrak g$ to $\mathfrak h$ we mean a linear map of vector spaces which preserves the Lie bracket.
		\end{defi}
		
		\begin{rema} \label{rema:contsmooth}
			A continuous homomorphism between Lie groups is automatically smooth (cf. \cite[Theorem 3.39]{war}). Therefore we could have defined Lie group homomorphisms as continuous maps but assumed them to be infinitely differentiable.
		\end{rema}
		
		\begin{prop} \label{prop:indlahom}
			Let $G$ and $H$ be two Lie groups and let $\varphi\colon G\to H$ be a homomorphism of Lie groups. Then by taking the derivative of $\varphi$ at $e$ we obtain a homomorphism of Lie algebras $L\varphi:=\varphi^\prime(e_G)\colon\mathfrak g\to\mathfrak h$\index{$L\varphi$, $\varphi^\prime(e)$}. That means the following holds.
			$$[\varphi^\prime(e_G)X,\varphi^\prime(e_G)Y]_H=\varphi^\prime(e_G)([X,Y]_G)\quad\forall X,\ Y\in\mathfrak g$$
		\end{prop}
		
		\begin{proof}
			\cite[Theorem 3.14]{war}.
		\end{proof}
		
		\begin{defi}
			Let $G$ and $H\subseteq G$ be Lie groups. We say that $H$ is a \textbf{Lie subgroup}\index{Lie group!subgroup} of $G$ if $H$ is a submanifold of $G$ such that the inclusion $i$ is a Lie group homomorphism.
		\end{defi}
		
		\begin{theo} \label{theo:existsubgroups}
			Let $G$ be a Lie group with Lie algebra $\mathfrak g$ and let $\mathfrak h_0$ be a subalgebra. Then up to isomorphism there exists a unique connected Lie subgroup $H$ of $G$ such that the associated Lie algebra homomorphism satisfies $Li(\mathfrak h)=\mathfrak h_0$.
		\end{theo}
		
		\begin{proof}
			\cite[Theorem 3.19]{war}.
		\end{proof}
		
		One can interpret the Lie algebra $\mathfrak g$ of a Lie group $G$ as its local version in the following sense. Let $X\in\mathfrak g$ and let $\widetilde{X}$ be the associated left-invariant vector field. Then by moving a solution around via left-invariance one finds that each integral curve is defined on $\mathbb R$. One then defines the \textbf{exponential map}\index{exponential map}\index{$\exp$, exponential map}
		$$\exp\colon\mathfrak g\to G,\quad X\mapsto\gamma_{X}(1),$$
		where $\gamma_X\colon\mathbb R\rightarrow G$ is the unique curve which solves $\widetilde X$ and satisfies $\gamma_X(0)=e$.
		
		\begin{prop} \label{prop:expprop}
			The exponential map satisfies
			\renewcommand{\theenumi}{(\roman{enumi})}
			\renewcommand{\labelenumi}{\theenumi}
			\begin{enumerate}
				\item $\exp(0)=e$ and $\exp$ is smooth with derivative $\exp^\prime(0)=\Id_\mathfrak g$;
				\item there is a neighborhood $U$ of $0$ such that $\exp_{|U}$ is a diffeomorphism onto its image; and
				\item if $\varphi\colon G\to H$ is a homomorphism of Lie groups then the following diagram is commutative.
				$$\begin{xy}
     				\xymatrix{
         			G \ar[r]^-\varphi     &   H \\
         			\mathfrak g \ar[r]_-{L\varphi} \ar[u]^-{\exp_G} &   \mathfrak{h} \ar[u]_-{\exp_H}
         			}
   			\end{xy}$$
			\end{enumerate}
		\end{prop}
		
		\begin{proof}
			\cite[Theorem 3.31, Theorem 3.32]{war}.
		\end{proof}
		
		One obvious question at this point is whether a homomorphism of Lie algebras can be lifted to a homomorphism of Lie groups. In general this is not possible, but we have the following two propositions.
		
		\begin{prop} \label{prop:liftexistforsc}
			Let $G$ and $H$ be Lie groups with Lie algebras $\mathfrak g$ and $\mathfrak h$ respectively, and with $G$ simply connected. Let $\psi:\mathfrak g\to\mathfrak h$ be a homomorphism. Then there exists a unique homomorphism $\varphi:G\to H$ such that $L\varphi=\psi$.
		\end{prop}
		
		\begin{proof}
			\cite[Theorem 3.27]{war}.
		\end{proof}
		
		\begin{prop} \label{prop:uniqueforconnected}
			Let $G$ be a connected Lie group and let $\varphi$ and $\psi$ be Lie group homomorphisms from $G$ to a Lie group $H$ such that $L\varphi=L\psi$. Then $\varphi=\psi$.
		\end{prop}
		
		\begin{proof}
			\cite[Theorem 3.16]{war}.
		\end{proof}
		
		\begin{prop} \label{prop:la}
			Let $G$ be a Lie group and $H$ a closed subgroup. Then the Lie algebra of $H$ is given by
			$$\mathfrak h=\{X\in\mathfrak g\ |\ \exp(tX)\in H\quad\forall t\in\mathbb R\}.$$
		\end{prop}
		
		\begin{proof}
			\cite[Proposition 3.33]{war}.
		\end{proof}
		
		\subsection{Basic Representation Theory}
		
		\begin{prop} \label{prop:glv} 
			Let $V$ be a finite-dimensional real or complex vector space. Then the group $\GL(V)$\index{$\GL(V)$, invertible linear transformations} of real- or complex-linear invertible transformations is a Lie group with Lie algebra $\End(V)$. Here the Lie bracket on $\End(V)$ is the commutator.
		\end{prop}
		
		\begin{proof}
			\cite[Example 3.10(b)]{war}.
		\end{proof}
		
		\begin{defi}
			Let $G$ be a Lie group and $V$ a real or complex topological vector space. Then a \textbf{representation}\index{representation!Lie group} $(\pi,V)$ of $G$ on $V$ is a strongly continuous group homomorphism $\pi\colon G\to\mathcal B(V)$. Here $\mathcal B(V)$ denotes the algebra of continuous linear operators on $V$. In particular, we have
			$$\pi(e)=\Id_V,\quad \pi(gh)=\pi(g)\circ\pi(h)\quad\forall g,\ h\in G.$$
			If $V$ is finite-dimensional, then this definition reduces to $\pi\colon G\to\GL(V)$, where $\GL(V)$ is a Lie group according to Proposition \ref{prop:glv} and $\pi$ is a Lie group homomorphism according to Remark \ref{rema:contsmooth}. The representation $(\pi,V)$ is said to be \textbf{unitary}\index{representation!unitary} if $V$ is a complex Hilbert space and $\pi(g)$ is unitary for each $g\in G$.
		\end{defi}
		
		As a consequence of Proposition \ref{prop:expprop} every finite-dimensional representation $(\pi,V)$ induces a Lie algebra homomorphism $L\pi\colon\mathfrak g\to\End(V)$ such that the following diagram is commutative.
		$$\begin{xy}
     				\xymatrix{
         			G \ar[r]^-\pi     &   \GL(V) \\
         			\mathfrak g \ar[r]_-{L\pi} \ar[u]^-{\exp} &   \End(V)\ar[u]_-{\exp}  
     				}
   	\end{xy}$$
		We call $(L\pi,V)$ the associated \textbf{Lie algebra representation}\index{representation!Lie algebra}.
		
		\begin{exam} \label{exam:char}
			A Lie group homomorphism from $G$ to the multiplicative group $\mathbb C^*\cong\GL(\mathbb C)$ is called a \textbf{character}\index{character of a Lie group} of $G$. Let $\pi:G\rightarrow S^1$ be a Lie group homomorphism. We view the sphere $S^1$\index{$S^1$, the 1-sphere} as the set of unitary operators on a one-dimensional complex Hilbert space. If we interpret $i\mathbb R$ as the tangent space at 1 and hence the Lie algebra of $S^1$, we can define the exponential map as the usual exponential function. The induced Lie algebra representation $L\pi\colon\mathfrak g\to i\mathbb R$ is then a functional on $\mathfrak g$.
		\end{exam}
		
		Given representations can be used to construct new ones. The easiest example is taking two representations $(\pi_1,V_1)$ and $(\pi_2,V_2)$ of $G$ and consider the direct sum $V=V_1\oplus V_2$. Then we define a representation $\pi:=\pi_1\oplus\pi_2$ of $G$ on $V$ by putting for each $g\in G$
		$$\pi(g)(v_1,v_2):=(\pi_1(g)v_1,\pi_2(g)v_2)\quad\forall(v_1,v_2)\in V.$$
		This technique works in an obvious way also for finitely many representations.
		
		\begin{defi}
			Let $(\pi,V)$ be a representation of a Lie group $G$. A closed subspace $S\subseteq V$ is called \textbf{invariant} if
			$$\pi(g)S\subseteq S\quad\forall g\in G.$$
			We say that $(\pi,V)$ is \textbf{irreducible}\index{representation!irreducible} if the only invariant subspaces are $\{0\}$ and $V$. Otherwise it is called \textbf{reducible}.
		\end{defi}
		
		If $(\pi,V)$ is a reducible representation of $G$ and $S$ is a non-trivial invariant subspace, then one can define a representation on $S$ by restriction. More precisely, one has the representation $(\pi_{|S},S)$, where $\pi_{|S}(g):=\pi(g)_{|S}$ for all $g\in G$. In the case that $(\pi,V)$ is unitary one can even decompose the representation in the following way. Let $S^\bot$ denote the orthogonal complement of $S$ in $V$ with respect to the Hermitian inner product $(.,.)$. Then $S^\bot$ is invariant, because for all $g\in G$
		\begin{align*}
			(v,\pi(g)w) & = (\pi(g)^*v,w) \\
									& = (\pi(g^{-1})v,w) = 0\quad\forall v\in S,\ w\in S^\bot.
		\end{align*}
		For the two invariant subspaces $S$ and $S^\bot$ we obtain representations $\pi_{|S}$ and $\pi_{|S^\bot}$. Furthermore we note that $\pi=\pi_{|S}\oplus\pi_{|S^\bot}$ which means that we have found a decomposition of $\pi$ into representations in smaller spaces. By employing induction this idea gives the following proposition.
		
		\begin{prop} \label{prop:finunirep}
			Every finite-dimensional unitary representation of a Lie group is the (finite) sum of irreducible representations.
		\end{prop}
		
		In the compact case, one also has the converse.
		
		\begin{prop} \label{prop:irredisfinite}
			Every irreducible unitary representation of a compact Lie group is finite-dimensional.
		\end{prop}
		
		\begin{proof}
			\cite[Corollary 9.5]{knapp}.
		\end{proof}
		
		Using the following result, one obtains a similar statement for abelian groups.
		
		\begin{theo}[Schur's Lemma]\index{Schur's Lemma}
			Let $(\pi,V)$ be a finite-dimensional unitary representation of a Lie group $G$. Set
			$$M:=\{A\in\GL(V)\ |\ A\pi(g)=\pi(g)A\quad\forall g\in G\}$$
			Then one has
			$$\pi\text{ irreducible}\quad\Leftrightarrow\quad M=\mathbb C\Id_V.$$
		\end{theo}
		
		\begin{proof}
			\cite[Corollary 4.9]{knapp}.
		\end{proof}
		
		\begin{coro} \label{coro:irreunitabel}
			Let $(\pi,V)$ be an irreducible unitary representation of an abelian Lie group $G$. Then $V$ is a one-dimensional complex vector space.
		\end{coro}
		
		In the following, we illustrate the existence of a left-invariant measure and consider its properties for compact groups.
		
		\begin{theo}[Haar]
			Let $G$ be a locally compact topological group. Then up to a positive scalar there exists a unique regular Borel measure $\mu$ such that
			$$\int_G f(gh)d\mu(h)=\int_G f(h)d\mu(h)$$
			for all measurable real or complex-valued functions $f$ on $G$.
		\end{theo}
		
		\begin{proof}
			If $G$ is an $n$-dimensional Lie group, then integration is realized by integrating differential forms of degree $n$ (cf. \cite[§4]{war}). To prove the theorem we need to find a left-invariant non-vanishing $n$-form. This is always possible by picking a non-zero $n$-form on $\mathfrak g$ from the one-dimensional space $\Lambda^n\mathfrak g$ of alternating $n$-forms and considering the associated left-invariant differential $n$-form. A proof of the general case can be found in \cite[§58]{halmos}.
		\end{proof}
		
		\begin{prop} \label{prop:haar}
			If $G$ is a compact Lie group, then $\mu$ is also right translation-invariant. Since $G$ has finite measure, we can assume $\mu$ to be a probability measure. We call the unique normed translation-invariant measure the \textbf{Haar measure}\index{Haar measure} on $G$.
		\end{prop}
		
		\begin{proof}
			\cite[p.239]{knapp}.
		\end{proof}
		
		\begin{prop} \label{prop:iip}
			Let $(\pi,V)$ be a finite-dimensional representation of a compact Lie group $G$ on a real or complex vector space. Then one can find an inner product on $V$ such that $G$ acts by orthogonal, respectively unitary transformations.
		\end{prop}
		
		\begin{proof}
			We start with an arbitrary inner product or Hermitian inner product $(.,.)$ on $V$ and average over the group:
			$$(v,w)_{inv}:=\int_G(\pi(h)v,\pi(h)w)d\mu(h).$$
			This is well-defined since $h\mapsto(\pi(h)v,\pi(h))$ is smooth and hence measurable. In particular, it is integrable over the compact group $G$. The invariance property follows from the translation-invariance of $\mu$ with the following calculation.
			\begin{align*}
				\left(\pi(g)v,\pi(g)w\right)_{inv} & = \int_G(\pi(h)(\pi(g)v),\pi(h)(\pi(g)w))d\mu(h) \\
																& = \int_G(\pi(hg)v,\pi(hg)w)d\mu(h) \\
																& = \int_G(\pi(h)v,\pi(h)w)d\mu(h) \\
																& = \left(v,w\right)_{inv}.
			\end{align*}
		\end{proof}
		
		\begin{defi}
			We say that two representations $(\pi_1,V_1)$ and $(\pi_2,V_2)$ are \textbf{equivalent}\index{representation!equivalence} if there is a vector space isomorphisms $A\colon V_1\to V_2$ such that for every $g\in G$ the following diagram is commutative.
			$$\begin{xy}
     				\xymatrix{
         			V_1 \ar[r]^-A     &   V_2 \\
         			V_1 \ar[r]_-A \ar[u]^-{\pi_1(g)} &   V_2\ar[u]_-{\pi_2(g)}  
     				}
   		\end{xy}$$
   		In the case of unitary representations one additionally demands that $A$ is a Hilbert space isomorphism.
		\end{defi}
		
		The relation between two representations to be equivalent is clearly reflexive, symmetric and transitive. Therefore it is an equivalence relation on the set of representations and we can talk about equivalence classes. The same is true for the equivalence of unitary representations. Furthermore if $(\pi_1,V_1)$ and $(\pi_2,V_2)$ are equivalent, then $\pi_1$ is irreducible if and only if $\pi_2$ is irreducible, since invariant subspaces stay invariant under the equivalence isomorphism.
		
		Wallach shows in \cite[Lemma 2.3.8]{walhomo} that if $(\pi,V)$ is an irreducible finite-dimensional unitary representation of $G$, then any $G$-invariant Hermitian inner product on $V$ is a real multiple of the first Hermitian inner product. As a corollary we obtain that two unitary irreducible finite-dimensional representations are equivalent as representations if and only if they are equivalent as unitary representations. With these observations the following definition makes sense.
		
		\begin{defi}
			Let $G$ be a Lie group. Then we define the \textbf{unitary dual}\index{unitary dual} $\hat G_u$\index{$\hat G_u$, unitary dual} as the set of equivalence classes of irreducible unitary representations of $G$.
		\end{defi}
		
		One major problem in representation theory is to describe $\hat{G}_u$ for a given Lie group $G$. In the case of a compact Lie group, the knowledge of $\hat{G}_u$ is especially significant. First of all, we know from Proposition \ref{prop:irredisfinite} that every irreducible representation is finite-dimensional. Furthermore, Proposition \ref{prop:iip} says that any finite-dimensional representation is unitary. Combined with Proposition \ref{prop:finunirep} we find that every finite-dimensional representations of $G$ is equivalent to the direct sum of certain elements of $\hat{G}_u$. If $G$ is in addition abelian, then all elements of $\hat{G}_u$ are of dimension one (cf. Corollary \ref{coro:irreunitabel}) and can be identified with characters of $G$ (cp. Example \ref{exam:char}).
		
	\subsection{Adjoint and Coadjoint Representation} \label{ss:adcoad}
		
		Let $G$ be Lie group. Since multiplication and inversion are smooth maps, the conjugation $I_g$\index{$I_g$, conjugation by $g$} by $g\in G$
		$$I_g:G\rightarrow G,\quad x\mapsto gxg^{-1}$$
		is smooth and we can calculate its derivative. Since $I_g$ is a Lie group automorphism the induced map $LI_g$ is a Lie algebra automorphism of $\mathfrak g$. We define a mapping $\Ad:G\rightarrow\GL(\mathfrak g)$, $g\mapsto LI_g$\index{$\Ad$, adjoint representation}. If $g,\ h\in G$, then $I_g\circ I_h=I_{gh}$ and the chain rule yields
		\begin{align*}
			\Ad(gh) & = I_{gh}^\prime(e) \\
							& = (I_g\circ I_h)^\prime(e) \\
							& = I_g^\prime(e)\circ I_h^\prime(e) \\
							& = \Ad(g)\circ \Ad(h).
		\end{align*}
		Therefore $\Ad$ is a Lie group homomorphism to $\GL(\mathfrak g)$ and hence a representation of $G$ in $\mathfrak g$. It is called the \textbf{adjoint representation}\index{adjoint representation}. The associated Lie algebra representation is called $\ad:\mathfrak g\rightarrow\End(\mathfrak g)$\index{$\ad$, $L\Ad$}. For $\Ad$ and $\ad$ we have the following commutative diagram
		$$\begin{xy}
   		\xymatrix{
     		G \ar[r]^-\Ad     &   \GL(\mathfrak g) \\
     		\mathfrak g \ar[r]_-\ad \ar[u]^-{\exp} &   \End{(\mathfrak g)}\ar[u]_-{\exp}  
   		}
 		\end{xy}$$
		and one can prove (cf. \cite[Lemma 8.1.9]{hga}) that $\ad$ is actually given by
		$$\ad(X)Y=\frac{d}{dt}|_{t=0}\Ad(\exp tX)(Y)=[X,Y].$$
			
		By $\mathfrak g^*$\index{$\mathfrak g^*$, dual space of $\mathfrak g$} we denote the \textbf{dual space of $\mathfrak g$}, which is the vector space of real-valued \textbf{functionals}\index{functional} on $\mathfrak g$, i.e. linear mappings from $\mathfrak g$ to $\mathbb R$. The dual space is a vector space of the same dimension as $\mathfrak g$. The value of a functional $\lambda$ at $v$ is often written as $\lambda(v)=<\lambda,v>$. This notation is employed to define for $A\in\GL(\mathfrak g)$ the \textbf{dual operator $A^*\in\GL(\mathfrak g^*)$} by the equation
		$$<\lambda,AX>=<A^*\lambda,X>\quad\forall\lambda\in\mathfrak g^*\ \forall X\in\mathfrak g.$$
		The dual or contragrediant version $\Ad^*$\index{$\Ad^*$, coadjoint representation} of the adjoint representation is given by
		$$<\Ad^*(g)\lambda,X>=<\lambda,\Ad(g^{-1})X>.$$
		Note that $\Ad^*(g)$ is not dual to $\Ad(g)$ as the notation might lead to believe, but rather $\Ad^*(g)=(\Ad(g^{-1}))^*$, and in particular
		$$<\Ad^*(g)\lambda,\Ad(g)X>=<\lambda,X>\quad\forall\lambda\in\mathfrak g^*\ \forall X\in\mathfrak g.$$
		Using the properties of $\Ad$ we observe that $\Ad^*$ is a representation of $G$ in $\mathfrak g^*$. It is called the \textbf{coadjoint representation}\index{coadjoint!representation} or \textbf{coadjoint action}\index{coadjoint!action}. The associated Lie algebra representation\index{$\ad^*$, $L\Ad^*$} is then given by
		$$\ad^*(X)\lambda(Y)=-\lambda([X,Y]).$$
		
	\subsection{Homogeneous Spaces}
		
		\begin{defi} \label{defi:action}
			Let $G$ be a Lie group and $M$ a set. An \textbf{action}\index{action} of $G$ on $M$ is a map $\sigma\colon G\times M\to M$ such that
			$$\sigma(gh,p)=\sigma(g,(\sigma(h,p))),\quad \sigma(e,p)=p\quad\forall g,\ h\in G\ \forall p\in M.$$
			If $M$ is a smooth manifold and $\sigma$ is a smooth map, then we call it a \textbf{smooth group action}.	If we fix $g\in G$, then $\ p\mapsto\sigma(g,p)$ is a diffeomorphism of $M$ which we denote by $\sigma(g)$\index{$\sigma(g)$, diffeomorphism induced by $\sigma$}. In fact, the map $g\mapsto\sigma(g)$ is a group homomorphism from $G$ to $D(M)$. If the action is clear from the context, one often writes $g\cdot p$ instead of $\sigma(g)p$. We say that the action is \textbf{transitive}\index{action!transitive} if for all $p,\ q\in M$ there exists a $g\in G$ such that $g\cdot p=q$. In this case we call $M$ a \textbf{homogeneous}\index{homogeneous space} $G$-space. By the derivative of $\sigma$ we mean the map $L{\sigma}\colon\mathfrak g\to\mathcal X(M)$\index{$L\sigma$, derivative of $\sigma$} defined by the formula
			$$L{\sigma}(X)f(p):=\frac{d}{dt}_{|t=0}f(\sigma(\exp -tX)p)\quad\forall p\in M\ \forall f\in C^\infty(M,\mathbb R).$$
			From \cite[Theorem 3.45]{war} we obtain that $L{\sigma}$ is a Lie algebra homomorphism. If $N$ is another homogeneous $G$-spaces, then we say that a map $f\colon M\to N$ is \textbf{$G$-equivariant}\index{equivariance} if $f(g\cdot p)=g\cdot f(p)$ for all $p\in M$ and all $g\in G$.
		\end{defi}
		
		\begin{rema}
			Sometimes an action defined as above is called a \textbf{left action}\index{action!left action}. It is put in contrast to right actions which are defined in an analogous manner. A \textbf{right action}\index{action!right action} $\sigma\colon M\times G\to M$ is an assignment $(p,g)\mapsto p\cdot g$ which satisfies $(p\cdot g)\cdot h=p\cdot(gh)$ and $p\cdot e=p$. Right actions intuitively express that an product $(gh)$ acts on $p$ by first applying $g$ and then $h$. Otherwise, the concept is completely equivalent to left actions: one can always obtain a left action via a right action by putting $g\cdot p:=p\cdot g^{-1}$.
		\end{rema}
		
		\begin{exam}
			Let $G$ be a Lie group and $(\pi,V)$ a finite-dimensional representation. Then $\sigma\colon G\times V\to V$ defined by $\sigma(g,v)=\pi(g)v$ is a smooth action by linear transformations. $G$-equivariant bijective linear maps from $V$ to other representation spaces correspond to equivalence isomorphisms.
		\end{exam}
		
		For an action of $G$ on $M$ and a point $p\in M$ we define the \textbf{orbit}\index{orbit} of $p$ under $G$ as
		$$G\cdot p=\{g\cdot p\ |\ g\in G\}.$$
		We observe that $p\in G\cdot p$ and that if $q\in G\cdot p$, then $G\cdot q=G\cdot p$. In fact we note that lying in the same orbit is an equivalence relation on $M$. The action is transitive if and only if $G\cdot p=M$ for one and hence all $p\in M$. The restriction of the action to an orbit is always transitive. Therefore any orbit which is a submanifold of $M$ is automatically a homogeneous $G$-space.
		
		\begin{exam} \label{exam:co}
			Given a Lie group $G$, we have defined in Subsection \ref{ss:adcoad} the coadjoint action $\Ad^*\colon G\to\GL(\mathfrak g^*)$.  Let $\lambda\in\mathfrak g^*$ be a functional on $\mathfrak g$. Then the set
			$$\mathcal O_\lambda:=\Ad^*(G)\lambda=\{\Ad^*(g)\lambda\ |\ g\in G\}$$
			is called the \textbf{coadjoint orbit}\index{coadjoint!orbit}\index{orbit!coadjoint}\index{$\mathcal O_\lambda$, coadjoint orbit of $\lambda$} of $\lambda$.
		\end{exam}
		
		\begin{theo} \label{theo:quotishomo}
			Let $G$ be a Lie group and $H\subseteq G$ a closed Lie group. Let $M:=G/H=\{gH\ |\ g\in G\}$\index{$G/H$, space of left cosets} be the space of left cosets of $G$ with respect to $H$ and $\pi\colon G\to G/H, \quad \pi(g)=gH$ the quotient map. Then there is a unique differentiable structure on $G/H$ such that
			\renewcommand{\theenumi}{(\roman{enumi})}
			\renewcommand{\labelenumi}{\theenumi}
			\begin{enumerate}
				\item $\pi$ is smooth.
				\item There is a neighborhood $U$ of $eH$ in $M$ and a smooth map $\tau\colon U\to G$ such that $\pi\circ\tau=\Id_U$.
				\item $M$ is a homogeneous $G$-space w.r.t. the smooth action $\sigma(g_1)g_2H:=g_1g_2H$.
			\end{enumerate}
		\end{theo}
		
		\begin{proof}
			\cite[Theorem 3.58]{war}.
		\end{proof}
		
		The following theorem gives a description of all homogeneous $G$-spaces. Let $G$ be a Lie group and $M$ a homogeneous $G$-space. For $p\in M$ we define the \textbf{stabilizer}\index{stabilizer} or \textbf{isotropy group}\index{isotropy group} as\index{$G_p$, stabilizer}
		$$G_p:=\{g\in G\ |\ \sigma(g)p=p\}.$$
		Then the following holds.
		
		\begin{prop} \label{prop:stab}
			$G_p$ is a closed Lie subgroup of $G$ whose Lie algebra is given by $\mathfrak g_p=\{X\in\mathfrak g\ |\ L{\sigma}(X)(p)=0\}$\index{$\mathfrak g_p$, Lie algebra of $G_p$}.
		\end{prop}
		
		\begin{proof}
			We observe that $G_p$ is a closed subgroup since the action is continuous. Then \cite[Theorem 3.42]{war} asserts that $G_p$ is a Lie subgroup of $G$. Furthermore the Lie algebra $\mathfrak g_p$ of $G_p$ consists of elements $X\in\mathfrak g$ such that for all $t\in\mathbb R$ we have $\exp(tX)\in G_p$ (cp. Proposition \ref{prop:la}). Therefore $X\in\mathfrak g_p$ satisfies $\sigma(\exp(tX))p=p$. By taking the derivative we conclude $L{\sigma}(X)(p)=0$. Conversely, if $L{\sigma}(X)(p)=0$, then the flow of $L{\sigma}(X)$ is the constant function $t\mapsto p$. Since on the other hand, the flow is given by $t\mapsto\sigma(\exp-tX)p$, we find $\sigma(\exp X)p=p$.
		\end{proof}
		
		\begin{theo} \label{theo:homoisquot}
			Let $\sigma\colon G\times M\to M$ be a smooth transitive action and $G_p$ the isotropy group of a point in $M$. Then the mapping
			$$G/G_p\to M,\quad gH\mapsto \sigma(g)p$$
			is a diffeomorphism.
		\end{theo}
		
		\begin{proof}
			\cite[Theorem 3.62]{war}.
		\end{proof}
		
		We can therefore identify any homogeneous $G$-space with a quotient space.
		
	\subsection{Basic Structure Theory} \label{structtheo}
		
		\begin{defi}
			Let $G$ be a Lie group. By a \textbf{torus}\index{torus} $T$ we mean a compact connected abelian subgroup of $G$. A torus $T$ is called \textbf{maximal}\index{torus!maximal} if there is no torus properly containing $T$. 
		\end{defi}
		
		In a compact Lie group $G$ one can always obtain a torus by considering the subgroup $\overline{\exp(\mathbb RX)}$ for a non-zero $X\in\mathfrak g$. Note that the so obtained torus does not need to be one-dimensional. In fact, the closure might be all of $G$ even if $G$ has dimension greater than one. One example is the so-called \textbf{dense wind}\index{dense wind} (cf. \cite[p.33]{hine}). A precise argument for the existence of maximal tori is \cite[Prop. 4.30]{knapp} which tells us that maximal tori correspond to maximal abelian subalgebras of $\mathfrak g$ (cf. Theorem \ref{theo:existsubgroups}). Moreover, we find the following properties of maximal tori (cf. \cite[Theorem 4.35, Theorem 4.50]{knapp}).
		
		\begin{theo} \label{theo:toriconjugated}
			In a compact connected Lie group, any two maximal tori are conjugate to each other.
		\end{theo}
		
		\begin{theo} \label{theo:existmaxitorus}
			Let $G$ be a compact connected Lie group and $T_1$ a torus in $G$. Then $T_1$ is contained in a maximal torus. Moreover, if $g\in G$ satisfies $gt=tg$ for all $t\in T_1$, then there is a torus $T$ which contains both $T_1$ and $g$.
		\end{theo}
		
		Let $G$ be a compact connected Lie group and $T$ a maximal torus. Proposition \ref{prop:iip} tells us that there is an inner product on $\mathfrak g$ such that $\Ad$ is an orthogonal representation. With respect to this inner product we can decompose $\mathfrak g=\mathfrak t\oplus\mathfrak t^\bot$. In this way, we will view functionals $\lambda\in\mathfrak t^*$ as functionals on $\mathfrak g$ extending them by zero on $\mathfrak t^\bot$.
		
		Now if we complexify $\mathfrak g$ to $\mathfrak g_\mathbb C:=\mathfrak g\oplus i\mathfrak g$, then we can extend the linear transformations $\Ad(G)$ to complex-linear transformations on $\mathfrak g_\mathbb C$. If we extend the $G$-invariant inner product on $\mathfrak g$ to a Hermitian inner product on $\mathfrak g_\mathbb C$, we find that $\Ad$ extends to a unitary representation of $G$ on $\mathfrak g_\mathbb C$. We will also denote it by $\Ad$.
		
		A good summary of the following ideas can also be found in \cite[§1]{vogur}.
		
		\begin{cons} (cf. \cite[§IV.5]{knapp} or \cite[§3.2]{walhomo}) \label{cons:rsd}
			Consider the restriction of $\Ad\colon G\to\GL(\mathfrak g_\mathbb C)$ to $T$. Then $\Ad(T)$ is a family of commuting diagonalizable unitary operators on $\mathfrak g_\mathbb C$. Linear Algebra tells us that $\Ad(T)$ must have a simultaneous eigenspace decomposition. One can show that it has the following form:
			$$\mathfrak g_\mathbb C=\mathfrak t_\mathbb C \oplus \sum_{\alpha\in\Delta}\mathfrak g^\alpha.$$
			Here for any functional $\alpha\in\mathfrak t_\mathbb C^*$ the space $\mathfrak g^\alpha$\index{$\mathfrak g^\alpha$ root space} is defined by the formula
			$$\mathfrak g^\alpha=\{X\in\mathfrak g_\mathbb C\ |\ [H,X]=\alpha(H)X\quad\forall H\in\mathfrak t_\mathbb C\}.$$
			$\Delta$\index{$\Delta$, set of roots} is defined as the finite subset of $\mathfrak t_\mathbb C^*$ of non-zero functionals $\alpha$ such that $\mathfrak g^\alpha\neq\{0\}$. We call the elements of $\Delta$ the \textbf{roots}\index{root} of $G$ with respect to $T$ and the spaces $\mathfrak g^\alpha$ are called \textbf{root spaces}\index{root!root space}.
		\end{cons}
		
		The relation between the roots and the eigenspaces of $\Ad(T)$ is the following. For each $\alpha\in\Delta$ there is a character $\chi_\alpha\colon T\to S^1$ such that $L\chi_\alpha=\alpha_{|\mathfrak t}$ and for each $t\in T$
		$$\Ad(t)X=\chi_\alpha(t)X\quad\forall X\in\mathfrak g^\alpha.$$
		The space $\mathfrak t_\mathbb C$ is exactly the eigenspace for the trivial character $t\mapsto 1$ whose differential is the zero functional which we did not include in $\Delta$. In the following proposition we summarize some properties of the roots. The proofs can be found in \cite[§IV.5]{knapp} and \cite[3.5]{wal}.
		
		\begin{prop} \label{prop:rootprop}
			For all $\alpha,\ \beta\in\Delta$ holds: 
			\renewcommand{\theenumi}{(\roman{enumi})}
			\renewcommand{\labelenumi}{\theenumi}
			\begin{enumerate}
				\item $\alpha_{|i\mathfrak t}$ is real-valued.
				\item $\mathfrak g^\alpha$ has complex dimension 1.
				\item $\overline{\mathfrak g^\alpha}=\mathfrak g^{-\alpha}$ and hence $-\alpha\in\Delta$.
				\item $[\mathfrak t_\mathbb C,\mathfrak g^\alpha]\subseteq\mathfrak g^\alpha.$
				\item $[\mathfrak g^\alpha,\mathfrak g^\beta]\left\{
\begin{array}{ll}
	=\mathfrak g^{\alpha+\beta} & \text{ if }\alpha+\beta\in\Delta\\ \subseteq\mathfrak t_{\mathbb C} & \text{ if }\beta=-\alpha\\ =\{0\} &\text{ else.}\\
\end{array}
\right.$
				\item $$\mathfrak g=\mathfrak t+\sum_{\alpha\in\Delta^+}\mathfrak g\cap(\mathfrak g^\alpha+\overline{\mathfrak g^\alpha}).$$
			\end{enumerate}
		\end{prop}
		
		Let $(.,.)$ be an $\Ad(G)$-invariant inner product on $\mathfrak g_\mathbb C$. We can use it to obtain an identification of $\mathfrak t_\mathbb C$ with $\mathfrak t_\mathbb C^*$ and $\mathfrak g_\mathbb C$ with $\mathfrak g^*_\mathbb C$. More precisely, for each $\xi\in\mathfrak t_\mathbb C^*$ there is a unique $H_\xi\in\mathfrak t_\mathbb C$ such that $(H_\xi,H)=\xi(H)$ for all $H\in\mathfrak t_\mathbb C$. We employ \cite[Lemma 3.5.8]{walhomo} to find that for each $\alpha\in\Delta$ we can choose $X_\alpha\in\mathfrak g^{\alpha}$ and $X_{-\alpha}\in\mathfrak g^{-\alpha}$ such that $H_\alpha=[X_\alpha,X_{-\alpha}]$.
		
		We use $(.,.)$ to define an inner product on $i\mathfrak t^*$ by duality. That is, if $\xi,\ \eta\in i\mathfrak t^*$, then we define $(\xi,\eta):=(H_\xi,H_\eta)$. Then our considerations yield the following proposition.
		
		\begin{prop} \label{prop:existrootvector}
			For each $\alpha\in\Delta$ there exists $X_\alpha\in\mathfrak g^\alpha,\ X_{-\alpha}\in\mathfrak g^{-\alpha}$ such that
			$$\xi([X_\alpha,X_{-\alpha}])=(\xi,\alpha)\quad\forall\xi\in i\mathfrak t^*.$$
		\end{prop}
		
		We define $\mathfrak t_R$ by 
		$$\mathfrak t_R:=\sum_{\alpha\in\Delta}\mathbb R H_\alpha\subset i\mathfrak t.$$
		
		\begin{defi}
			A connected component $P$ of $\mathfrak t_R\backslash(\bigcup_{\alpha\in\Delta}\ker\alpha)$ is called a \textbf{Weyl chamber}\index{Weyl chamber} of $T$. From this definition it is clear that each root $\alpha$ is either positive or negative on $P$. Then
			$$\Delta^+_P=\{\alpha\in\Delta\ |\ \alpha(P)>0\}$$
			is called the set of \textbf{positive roots}\index{root!positive} with respect to $P$.
			
			Let $\pi$ be a representation of $T$. Apply construction \ref{cons:rsd} to $\pi(T)$ instead of $\Ad(T)$ to obtain an eigenspace decomposition and call the resulting functionals on $\mathfrak t_\mathbb C$ the \textbf{weights}\index{weight} of $\pi$. Here the the zero functional is not excluded as it was done in the definition of roots. A weight $\lambda$ is called \textbf{dominant}\index{weight!dominant}, if
			$$(H_\lambda,H_\alpha)\ge 0\quad\forall\alpha\in\Delta_P^+.$$
		\end{defi}
		
		\begin{defi}
			Let $G$ be a compact connected Lie group. Fix a maximal torus and consider the normalizer $N(T)=\{n\in G\ |\ nTn^{-1}=T\}$. Then the group $W(G,T):=N(T)/T$ is called the \textbf{Weyl group}\index{Weyl group}\index{$W(G,T)$, Weyl group} of $G$ with respect to $T$ (cf. \cite[§IV.6]{knapp}).
		\end{defi}
		
		From Theorem \ref{theo:toriconjugated} we know that any other maximal torus has the form $gTg^{-1}$. Its normalizer $N(gTg^{-1})$ has the form $gN(T)g^{-1}$ and therefore $W(G,T)$ and $W(G,gTg^{-1})$ are isomorphic via conjugation. Hence the Weyl group is independent of the choice of $T$ for our purposes.
		
		The Weyl group $W(G,T)$ acts on $T$ by conjugation: $nT\cdot t=ntn^{-1}$ for $t\in T$. The derivative of the conjugation by $n$ is the Lie algebra automorphism $\Ad(n)_{|\mathfrak t}$ and therefore we also obtain an action of $W(G,T)$ on $\mathfrak t$. For the Weyl group and its actions, we have the following theorem.
		
		\begin{theo}
			The Weyl group $W(G,T)$ satisfies: \label{theo:wgroup}
			\renewcommand{\theenumi}{(\roman{enumi})}
			\renewcommand{\labelenumi}{\theenumi}
			\begin{enumerate}
				\item $W(G,T)$ is a finite group.
				\item The following map between the set of Weyl orbits in $T$ and the set of conjugacy classes of $T$ in $G$ is a bijection. \label{theo:wgroup:conjclass}
				$$W\cdot t\mapsto\{gtg^{-1}\ |\ g\in G\},\quad t\in T$$
				\item The following map between the set of Weyl orbits in $\mathfrak t$ and the set of adjoint orbits is a bijection. \label{theo:wgroup:adorb}
				$$W\cdot H\mapsto G\cdot H, \quad H\in\mathfrak t$$
			\end{enumerate}
		\end{theo}
		
		\begin{proof}
			Statements $(i)$ and $(ii)$ are proved in \cite[Proposition 3.10.8]{walhomo} and \cite[Proposition 4.53]{knapp}. Statement $(iii)$ is a differentiated version of $(ii)$ and an idea for its proof can be found in \cite[§3.8]{dk}.
		\end{proof}
		
		From Theorem \ref{theo:existmaxitorus} we know that each element lies in a maximal torus. Theorem \ref{theo:toriconjugated} tells us that all maximal tori are conjugate. Consequently, if $T$ is a maximal torus in a compact connected Lie group, then each element in $G$ is conjugate to an element of $T$. Hence \ref{theo:wgroup}\ref{theo:wgroup:conjclass} says that the Weyl orbits in $T$ parametrize the conjugacy classes of $G$.
		
		\begin{coro} \label{coro:fots}
			Similar to our previous considerations, we use \cite[4.34]{knapp} which tells us that any two maximal abelian subalgebras of $\mathfrak g$ are conjugate via $\Ad(G)$. Then \ref{theo:wgroup}\ref{theo:wgroup:adorb} says that the Weyl orbits in $\mathfrak t$ parametrize the set of adjoint orbits of $G$ in $\mathfrak g$. Let $(.,.)$ be the $\Ad(G)$-invariant inner product on $\mathfrak g$ from Construction \ref{cons:rsd}. We use it to identify $\mathfrak t$ with $\mathfrak t^*$ and $\mathfrak g$ with $\mathfrak g^*$. Then the parametrization means that each coadjoint orbit intersects $\mathfrak t^*$ in at most $|W(G,T)|$ points. In particular, this tells us that each coadjoint orbit in $\mathfrak g^*$ is given by $\mathcal O_\lambda$ for some $\lambda\in\mathfrak t^*$.
		\end{coro}
	
	\subsection{Centralizers of Tori}
		
		\begin{defi}
			Let $G$ be a compact connected Lie group and $S\subseteq G$ a subset. Then the \textbf{centralizer}\index{centralizer} $C(S)$\index{$C(S)$, centralizer} of $S$ is the subgroup
			$$C(S):=\{g\in G\ |\ gsg^{-1}=s\quad\forall s\in S\}.$$
			If $T_1$ is a torus in $G$, then $G/C(T_1)$\index{$G/C(T_1)$, generalized flag manifold} is called a \textbf{generalized flag manifold}\index{generalized flag manifold}.
		\end{defi}
		
		We will show in Proposition \ref{prop:cent} that $C(T_1)$ is a Lie group and then Theorem \ref{theo:quotishomo} implies that generalized flag manifolds are homogeneous $G$-spaces.
		
		\begin{prop} \label{prop:maxtcent}
			Let $C(T_1)$ be a centralizer of a torus in a compact connected Lie group $G$. Then $C(T_1)$ is connected. If $T_1$ is a maximal torus, then $C(T_1)=T_1$.
		\end{prop}
		
		\begin{proof}
			We apply Theorem \ref{theo:existmaxitorus}. If $g\in G$ centralizes $T_1$ then there is a torus $T$ containing both $g$ and $T_1$. This implies $T\subseteq C(T_1)$ and therefore the centralizer is connected. If $T_1$ is a maximal torus, then any torus $T$ containing $T_1$ coincides with $T_1$ and hence $g\in T=T_1$ for each $g\in C(T_1)$.
		\end{proof}
		
		\begin{prop} \label{prop:cent}
			Let $T_1$ be a torus with Lie algebra $\mathfrak t_1$ in a compact connected Lie group $G$. Then the centralizer is of the form $C(T_1)=\{g\in G\ |\ \Ad(g)H=H\quad\forall H\in \mathfrak t_1\}$. Furthermore $C(T_1)$ is a Lie group with Lie algebra $\mathfrak c_1=\{X\in\mathfrak g\ |\ \ad(X)H=0\quad\forall H\in\mathfrak t_1\}$.
		\end{prop}
		
		\begin{proof}
			From the continuity of the group operations we see that $C(T_1)$ is a closed subgroup. Then we can use \cite[Theorem 3.42]{war} to see that $C(T_1)$ is indeed a Lie subgroup. Since $T_1$ is connected, we obtain $I_g(h)=ghg^{-1}=h$ if and only if $\Ad(g)H=I_g^\prime(e)H=H$ for all $H\in T_1$ (cf. \cite[§3.1]{dk}). The identification of the Lie algebra is then similar to the proof of Proposition \ref{prop:stab}.
		\end{proof}
		
		Let $G$ be a compact connected Lie group, $T_1$ a torus and $T$ a maximal torus containing $T_1$. Then $T$ is a maximal torus in $C(T_1)$. Therefore we can apply Construction \ref{cons:rsd} to obtain a root system $\Delta$ of $G$ relative to $T$ and a root system $\Delta_1$ of $C(T_1)$ relative to $T_1$. If $P$ is a Weyl chamber of $\Delta$ then it is said to be $T_1$-admissible, if there is a Weyl chamber $P_1$ of $\Delta_1$ such that
		\renewcommand{\theenumi}{(\roman{enumi})}
		\renewcommand{\labelenumi}{\theenumi}
		\begin{enumerate}
			\item $\Delta^+\cap\Delta_1=\Delta_{P_1}^+$; and
			\item if $\alpha\in\Delta^+\backslash\Delta_{P_1}^+$, $\beta\in\Delta_1$ and $\alpha+\beta\in\Delta$, then $\alpha+\beta\in\Delta^+\backslash\Delta_{P_1}^+$.
		\end{enumerate}
		
		\begin{prop} \label{prop:existwc}
			There always exist $T_1$-admissible Weyl chambers.
		\end{prop}
		
		\begin{proof}
			\cite[Lemma 6.2.9]{walhomo}.
		\end{proof}
		
		In light of Proposition \ref{prop:existwc} we give an alternative construction of centralizers of tori (cf. \cite[6.2.10]{walhomo}).
		
		\begin{cons} \label{cons:altdefcent}
			Let $G$ be a compact connected Lie group and $T$ a maximal torus in $G$. Let $\Delta$ be the root system of $G$ with respect to $T$ and let $\Delta^+$ be a system of positive roots. Choose a subset $\Delta_1^+$ of $\Delta^+$ such that
			\begin{equation} \label{eq:closedroots}
				\alpha,\ \beta\in\Delta_1^+, \alpha+\beta\in\Delta^+\quad\Rightarrow\quad \alpha+\beta\in\Delta_1^+.
			\end{equation}
			Define $\Delta_1:=-\Delta_1^+\cup\Delta_1^+$ and set $\mathfrak c_1:=\mathfrak g\cap\{\mathfrak t_\mathbb C+\sum_{\alpha\in\Delta_1}\mathfrak g^\alpha\}=\mathfrak t+\sum_{\alpha\in\Delta_1^+}\mathfrak g\cap(\mathfrak g^\alpha+\overline{\mathfrak g^\alpha})$ and $\mathfrak t_1:=\{H\in\mathfrak t\ |\ \alpha(H)=0\quad \forall\alpha\in\Delta_1\}$. Then $\mathfrak t_1$ is an abelian subalgebra of $\mathfrak g$ and corresponds to the torus $T_1=\{t\in T\ |\ \chi_\alpha(t)=1\quad\forall\alpha\in\Delta_1^+\}$.
		\end{cons}
		
		\begin{prop} \label{prop:altdefcent} 
			Let $\Delta_1$ and $\mathfrak c_1$ be as in Construction \ref{cons:altdefcent}. Then $C(T_1)$ is the uniquely determined connected Lie subgroup having Lie algebra $\mathfrak c_1$.
		\end{prop}
		
		\begin{proof}
			The uniqueness follows from Theorem \ref{theo:existsubgroups}. Let $\mathfrak c=\{X\in\mathfrak g\ |\ \ad(X)H=0\quad\forall H\in\mathfrak t_1\}$ be the Lie algebra of $C(T_1)$ (cp. Proposition \ref{prop:cent}). Let $H\in\mathfrak t_1$ and let $X=X_\mathfrak t+\sum_{\alpha\in\Delta^+}(X_\alpha+\overline{X_\alpha})$ be an arbitrary element of $\mathfrak g=\mathfrak t+\sum_{\alpha\in\Delta^+}\mathfrak g\cap(\mathfrak g^\alpha+\overline{\mathfrak g^\alpha})$ (cf. Proposition \ref{prop:rootprop}). Then we have
			\begin{align*}
				\ad(X)H & = -\ad(H)X \\
								& = -\ad(H)X_\mathfrak t - \sum_{\alpha\in\Delta^+} \ad(H)(X_\alpha+\overline{X_\alpha}) \\
								& = - \sum_{\alpha\in\Delta^+} \alpha(H)(X_\alpha-\overline{X_\alpha}).
			\end{align*}
			Hence we find $\mathfrak c_1\subseteq\mathfrak c$. Conversely, if $\ad(X)H=0$ for all $H\in\mathfrak t_1$, then for each $\alpha$ either $X_\alpha=0$ or $\forall H\in\mathfrak t_1:\ \alpha(H)=0$ holds. Therefore $\mathfrak c\subseteq\mathfrak c_1$.
		\end{proof}
		
	\section{Complex Structures on Quotients}
	
	\subsection{Complex Manifolds}
		
		Let $V$ be a real $2n$-dimensional vector space. A \textbf{complex structure}\index{complex structure! on a vector space} on $V$ is a linear transformation $J\colon V\to V$ such that $J^2=-\Id_V$. The idea is that $V$ is actually a complex vector space $V_J$, but one has forgotten how to multiply with complex scalars. More precisely, we define the complex vector space $V_J$ to be the vector space $V$ by extending real scalar multiplication to complex scalars via the formula
		$$(x+iy)v=xv+yJv\quad\forall(x+iy)\in\mathbb C\quad\forall v\in V_J.$$
		Then $V_J$ is a complex vector space of complex dimension $n$. In particular, complex structures only exist on even-dimensional vector spaces. If $V$ is a complex vector space, we denote by $J_V$ its natural complex structure $J_V:v\mapsto iv$.
		
		Let $V$ be a real vector space together with a complex structure $J$. We extend $J$ uniquely to a complex-linear transformation $J$ on $V_\mathbb C:=V\oplus i V$ and note that $J^2=-\Id_{V_\mathbb C}$. Then $V_\mathbb C=V^{1,0}\oplus V^{0,1}$ decomposes into the eigenspaces of $i$ and $-i$ respectively, which are given by
		$$V^{1,0}=\{X-iJX\ |\ X\in V\},\quad V^{0,1}=\{X+iJX\ |\ X\in V\}.$$
		
		\begin{prop} \label{prop:csccd}
			There is a one-to-one correspondence between complex structures on $V$ and decompositions $V_\mathbb C=P\oplus\overline P$ into complex subspaces.
		\end{prop}
		
		\begin{proof}
			We have already constructed the decomposition $V_\mathbb C=V^{1,0}\oplus V^{0,1}$. Furthermore we note $\overline{V^{1,0}}=V^{0,1}$. Conversely, if we have $V_\mathbb C=P\oplus\overline P$, then for any $v\in V$ there is a uniquely determined $z\in P$ such that $v=z+\overline{z}$. We can then define a linear transformation $J\colon V\to V$ by the formula
			$$Jv:=iz+\overline{iz}$$
			and observe $J^2=-\Id_V$. In addition, $J$ inverts the construction above.
		\end{proof}
		
		\begin{defi} \label{defi:cs}
			Let $M$ be a second countable Hausdorff topological space. We say that $M$ is a \textbf{complex manifold}\index{manifold!complex} of dimension $n$ if there is given a family of pairs $\{U_j,\varphi_j\}_{j\in J}$, where $\{U_j\}_{j\in J}$ is an open covering of $M$ and for each $j\in J$ the map $\varphi_j\colon U_j\to\varphi_j(U_j)$ is a homeomorphism of $U_j$ onto an open set in $\mathbb C^n$ such that
			$$\varphi_k\circ\varphi_j^{-1}\colon \varphi_j(U_j\cap U_k)\to\varphi_k(U_j\cap U_k)$$
			is a holomorphic map whenever $k,\ j\in J$ are such that $U_j\cap U_k\neq\emptyset$. If $U\subseteq\mathbb C^n$, then by a holomorphic map $f\colon U\to\mathbb C^m$ we mean a smooth map such that at each $p\in U$ the derivative $f^\prime(p)\colon\mathbb C^n\to\mathbb C^m$ is complex-linear.
		\end{defi}
		
		The family $\{U_j,\varphi_j\}_{j\in J}$ is called a \textbf{holomorphic atlas} for $M$. A holomorphic atlas on $M$ which is maximal with respect to inclusion is called a \textbf{complex structure}\index{complex structure!on a manifold} on $M$. Similar to the case of atlases of smooth manifolds, a holomorphic atlas uniquely determines a complex structure. We see that a complex manifold of dimension $n$ is a smooth manifold of dimension $2n$ if we interpret $\mathbb C^n$ as $\mathbb R^{2n}$.
		
		\begin{defi}
			Let $M$ be a complex manifold with a holomorphic atlas $\{U_j,\varphi_j\}_{j\in J}$. We say that $f\in C^\infty(M,\mathbb C)$ is \textbf{holomorphic}\index{function!holomorphic} if for all pairs $(U_j,\varphi_j)$ the maps $f\circ\varphi_j^{-1}$ are holomorphic. If $N$ together with $\{V_k,\psi_k\}$ is another complex manifold, then a smooth map $f\colon M\to N$ is called holomorphic\index{map!holomorphic} if for all $p\in M$ there are $j,\ k$ such that $p\in U_j$, $f(p)\in V_k$ and the map
			\begin{equation} \label{eq:holomchart}
				\psi_k\circ f\circ\varphi_j^{-1}\colon\varphi_j(U_j)\to\psi_k(V_k)
			\end{equation}
			is holomorphic. The complex vector space of holomorphic functions on $M$ is denoted by $\Hol(M,\mathbb C)$\index{$\Hol(M,\mathbb C)$, holomorphic functions}.
		\end{defi}
		
		Let $p\in M$ and interpret $M$ as a smooth manifold of dimension $2n$. Then the real tangent space $T_pM$ is the $2n$-dimensional vector space of derivations of $C^\infty(M,\mathbb R)$. If $(U_j,\varphi_j)$ is a chart with $p\in U_j$, then the derivative $\varphi^\prime(p)\colon T_pM\to\mathbb C^n$ is a real-linear isomorphism with which we can give $T_pM$ a complex structure $J_p$ obtained from $\mathbb C^n$ and make it a complex vector space. Using (\ref{eq:holomchart}) and the chain rule we find that $J_p$ is independent of the chart used. Moreover, if $N$ is another complex manifold, we see from the definition that a smooth map $f\in C^\infty(M,N)$ is holomorphic if and only if its derivative $f^\prime(p)$ is complex-linear with respect to $J_p$ and $J_{f(p)}$ for all $p\in M$ (cf. \cite[Proposition IX.2.3]{kn}).
		
		We denote by $(T_pM)_\mathbb C$ the complexified tangent space at $p$ and consider the decomposition $(T_pM)_\mathbb C=T_pM^{1,0}\oplus T_pM^{0,1}$ from Proposition \ref{prop:csccd}. We call $T_pM^{1,0}$ the \textbf{holomorphic tangent space}\index{tangent space!holomorphic} and $T_pM^{0,1}$ the \textbf{antiholomorphic tangent space}\index{tangent space!antiholomorphic}. If we view $(T_pM)_\mathbb C$ as the $2n$-dimensional complex vector space of derivations of $C^\infty(M,\mathbb C)$, then $T_pM^{1,0}$ corresponds to the subspace of derivations which vanish on antiholomorphic functions (cf. \cite[p.16]{gh}). Furthermore, we find that if $f\in C^\infty(M, N)$ and if we extend the derivative at $p$ to a complex-linear map $f^\prime(p)\colon(T_pM)_\mathbb C\to (T_{f(p)}M)_\mathbb C$ for all $p\in M$, then $f$ is holomorphic if and only if $f^\prime(p)(T_pM^{1,0})\subseteq T_{f(p)}N^{1,0}$ for all $p\in M$ (cf. \cite[Proposition IX.2.9]{kn}).
		
		\begin{defi}
			Let $M$ be a smooth manifold. An \textbf{almost complex structure}\index{complex structure! almost complex structure} on $M$ is a tensor field $\mathcal J$ of type $(1,1)$ which is, at every point $p\in M$, a complex structure on $T_pM$. Here one interprets the bilinear form $\mathcal J_p$ as an automorphism $v\mapsto\mathcal J_p(v,.)\in\Hom(T_pM^*,\mathbb R)\cong T_pM$. We call the pair $(M,\mathcal J)$ an \textbf{almost complex manifold}\index{manifold!almost complex} and note that $M$ necessarily has even dimension. If $\mathcal J$ is an almost complex structure on a homogeneous $G$-space $M$ with action $\sigma$, then $\mathcal J$ is called \textbf{invariant}\index{complex structure!invariant}, if for all $g\in G$
			$$\mathcal J_{\sigma(g)p}\left((\sigma(g))^\prime(p)v\right)=\sigma(g)^\prime(p)\left(\mathcal J_{p}(v)\right)\quad\forall p\in M\quad\forall v\in T_pM.$$
			An almost complex structure $\mathcal J$ is said to be \textbf{integrable}\index{complex structure!integrable}, if it has no torsion, that is, if (cf.\cite[p.123]{kn})
			$$N(\mathfrak X,\mathfrak Y):=[\mathcal J\mathfrak X,\mathcal J\mathfrak Y]-[\mathfrak X,\mathfrak Y]-\mathcal J[\mathfrak X,\mathcal J\mathfrak Y]-\mathcal J[\mathcal J\mathfrak X,\mathfrak Y]=0\quad\forall \mathfrak X,\ \mathfrak Y\in\mathcal X(M).$$
		\end{defi}
		
		\begin{theo}[Newlander-Nirenberg]\index{Newlander-Nirenberg Theorem}
			An almost complex manifold $(M,\mathcal J)$ is a complex manifold such that $\mathcal J$ coincides with the natural complex structure on the tangent spaces, if and only if $\mathcal J$ is integrable.
		\end{theo}
		
		\begin{proof}
			\cite[Theorem IX.2.5]{kn}.
		\end{proof}
		
		\begin{prop} \label{prop:iacs}
			Let $G$ be a Lie group and $H$ a closed Lie subgroup. Then the invariant almost complex structures $\mathcal J$ on $M:=G/H$ (cp. Theorem \ref{theo:quotishomo}) are in one-to-one correspondence with the set of linear transformations $J$ of $T_{eH}M=\mathfrak g/\mathfrak h$ such that
			\renewcommand{\theenumi}{(\roman{enumi})}
			\renewcommand{\labelenumi}{\theenumi}
			\begin{enumerate}
				\item $J^2=-\Id$
				\item $J\circ\tilde{\Ad}(h)=\tilde{\Ad}(h)\circ J\quad\forall h\in H$.
			\end{enumerate}
			Here $\tilde{\Ad}:\mathfrak g/\mathfrak h\to\mathfrak g/\mathfrak h$ is the map induced by $\Ad:\mathfrak g\to\mathfrak g$.
		\end{prop}
		
		\begin{proof}
			Let $\mathcal J$ be an invariant almost complex structure on $M$ and put $J:=\mathcal J_{eH}$. Then $J^2=-\Id$. The action $l_x$ of $x$ on $M$ is defined by $gH\mapsto xgH$. For $h\in H$ we find $hgH=hgh^{-1}H$ which means the action of $h$ is given by conjugation modulo $H$, i.e. $l_h=\tilde{I}_h$. Therefore $l_h^\prime(eH)=\tilde{\Ad}(h)$. Now using the invariance of $\mathcal J$ and noting that $H$ fixes $eH$ we find
			\begin{align*}
				J(\tilde{\Ad}(h)X) & = \mathcal J_{eH}(l_h^\prime(eH)X) \\
													 & = l_h^\prime(eH)\left(\mathcal J_{eH}(X)\right) \\
													 & = \tilde{\Ad}(h)\left(J(X)\right).
			\end{align*}
			
			Now let $J$ be given as in the proposition. We define an invariant almost complex structure $\mathcal J$ by moving $J$ around via the action. More precisely, we put			
			$$\mathcal J_{xH}( l_x^\prime(eH) X):=l_x^\prime(eH)I(X)  \quad\forall X\in T_{eH}M\quad \forall x\in G.$$
			There are several possible elements to reach $xH$ from $eH$, but if $x,\ y\in G$ are such that $xH=yH$, then $y^{-1}x\in H$. Furthermore $l_x=l_y\circ\tilde{I}_{y^{-1}x}$. Hence we find
			\begin{align*}
				l_x^\prime(eH)=l_y^\prime(eH)\circ\tilde{\Ad}(y^{-1}x),
			\end{align*}
			which shows that $\mathcal J$ is well-defined by the above formula.
		\end{proof}
		
		If an invariant $\mathcal J$ is given in either way, we define $\mathcal P$ to be the subbundle of $(TM)_\mathbb C$ which is at each point $gH$ the eigenspace to $i$ of $\mathcal J_{gH}$ (cf. Proposition \ref{prop:csccd}). If $\mathfrak p:=(T_{eH}M)^{1,0}$ is the $i$-eigenspace of $J$, then one could also define $\mathcal P$ to be the $G$-invariant subbundle of $(TM)_\mathbb C$ having $\mathcal P_{eH}=\mathfrak p$. This definition is possible, since $\mathfrak p$ is $\Ad(H)$-invariant by the invariance of $J$. Furthermore these two definitions are equivalent, since $\mathcal J$ is the invariant extension of $J$ to all tangent spaces.
		
		\begin{prop} \label{prop:iacsiiffclosed}
			Let $\mathcal J$ be an invariant almost complex structure on $M=G/H$. Then $\mathcal J$ is integrable if and only if $\mathcal P$ is closed under the Lie bracket.
		\end{prop}
		
		\begin{proof}
			By the eigenspace property $\mathcal P$ is given by $\mathcal P=\{\mathfrak X-i\mathcal J\mathfrak X\ |\ \mathfrak X\in \mathcal X(M)\}$. Suppose it is closed under the Lie bracket and let $\mathfrak X-i\mathcal J\mathfrak X$ and $\mathfrak Y-i\mathcal J\mathfrak Y\in\mathcal P$. Then we have
			\begin{align*}
				0 & = \mathcal J([\mathfrak X-i\mathcal J\mathfrak X,\mathfrak Y-i\mathcal J\mathfrak Y])-i([\mathfrak X-i\mathcal J\mathfrak X,\mathfrak Y-i\mathcal J\mathfrak Y]) \\
					& = \mathcal J[\mathfrak X,\mathfrak Y] -i\mathcal J[\mathfrak X,\mathcal J\mathfrak Y] -i\mathcal J[\mathcal J\mathfrak X,\mathfrak Y] -\mathcal J[\mathcal J\mathfrak X,\mathcal J\mathfrak Y]\\
					&\ -i[\mathfrak X,\mathfrak Y] -[\mathfrak X,\mathcal J\mathfrak Y] -[\mathcal J\mathfrak X,\mathfrak Y] +i[\mathcal J\mathfrak X,\mathcal J\mathfrak Y].
			\end{align*}
			The imaginary part of this last expression is exactly $N$ which then has consequently to be equal to zero.
			
			Suppose that $N$ is zero. Then $[\mathfrak X,\mathfrak Y]-[\mathcal J\mathfrak X,\mathcal J\mathfrak Y]=\mathcal J[\mathfrak X,\mathcal J\mathfrak Y]+\mathcal J[\mathcal J\mathfrak X,\mathfrak Y]$. Computing the Lie bracket for two arbitrary elements of $\mathcal P$ we find
			\begin{align*}
				[\mathfrak X-i\mathcal J\mathfrak X,\mathfrak Y-i\mathcal J\mathfrak Y] & = [\mathfrak X,\mathfrak Y] -i[\mathfrak X,\mathcal J\mathfrak Y] -i[\mathcal J\mathfrak X,\mathfrak Y] -[\mathcal J\mathfrak X,\mathcal J\mathfrak Y] \\
											& = -\mathcal J[\mathfrak X,\mathcal J\mathfrak Y]-\mathcal J[\mathcal J\mathfrak X,\mathfrak Y] -i[\mathfrak X,\mathcal J\mathfrak Y] -i[\mathcal J\mathfrak X,\mathfrak Y] \\
											& = \left(-\mathcal J[\mathfrak X,\mathcal J\mathfrak Y]-\mathcal J[\mathcal J\mathfrak X,\mathfrak Y]\right) - i\mathcal J\left(-\mathcal J[\mathfrak X,\mathcal J\mathfrak Y]-\mathcal J[\mathcal J\mathfrak X,\mathfrak Y]\right)\ \in\mathcal P.
			\end{align*}
		\end{proof}
		
		Up to now we found that invariant complex structures come from certain $\tilde{\Ad}(H)$-invariant decompositions of $(\mathfrak g/\mathfrak h)_\mathbb C$ which are in some sense closed under the Lie bracket. The next proposition gives the final characterization of invariant complex structures with which we will continue to work.
		
		\begin{prop} \label{prop:csalg} 
			Suppose $G$ is a Lie group and $H\subseteq G$ is closed Lie subgroup. The set of invariant complex structures on the homogeneous $G$-space $M=G/H$ is in natural bijection with the set of complex subalgebras $\mathfrak b$ of $\mathfrak g_\mathbb C$, having the following properties:
			\renewcommand{\theenumi}{(\roman{enumi})}
			\renewcommand{\labelenumi}{\theenumi}
			\begin{enumerate}
				\item $\mathfrak b$ contains $\mathfrak h_\mathbb C$, and $\Ad(H)$ preserves $\mathfrak b$
				\item the intersection of $\mathfrak b$ with its complex conjugate $\mathfrak b^-$ is precisely $\mathfrak h_\mathbb C$
				\item the dimension of $\mathfrak b/\mathfrak h_\mathbb C$ is half the dimension of $\mathfrak g_\mathbb C/\mathfrak h_\mathbb C$.
			\end{enumerate}
		\end{prop}
		
		\begin{proof}
			Let $\mathfrak b$ with the above properties be given. We denote by $\pi:\mathfrak g_\mathbb C\to(\mathfrak g/\mathfrak h)_\mathbb C$ the quotient map and put $\mathfrak p:=\pi(\mathfrak b)$. Then $\mathfrak p\cap\overline{\mathfrak p}=\{0\}$ and since $2\dim\mathfrak p=\dim(\mathfrak g/\mathfrak h)_\mathbb C$ we obtain $(\mathfrak g/\mathfrak h)_\mathbb C=\mathfrak p\oplus\overline{\mathfrak p}$. Furthermore, $\mathfrak p$ is $\tilde{\Ad}(H)$-invariant since
			\begin{align*}
				\tilde{\Ad}(h)\mathfrak p & = \tilde{\Ad}(h)\circ\pi(\mathfrak b) \\
																	& = \pi\circ\Ad(h)(\mathfrak b) \\
																	& = \pi(\mathfrak b)\\
																	& = \mathfrak p.
			\end{align*}
			Thus the complex structure $J$ obtained via Proposition \ref{prop:csccd} from $\mathfrak p$ is $\tilde{\Ad}(H)$-invariant and gives us an invariant almost complex structure $\mathcal J$ on $G/H$ according to Proposition \ref{prop:iacs}. Then Proposition \ref{prop:iacsiiffclosed} tells us that to show integrability it suffices to show that $\mathcal P$ is closed under the Lie bracket. This will be handled in the next proposition.
			
			Now for the other direction, let $\mathcal J$ be an invariant complex structure on $G/H$. We define $\mathfrak p$ and $\mathcal P$ as before. We put $\mathfrak b:=\pi^{-1}(\mathfrak p)$ which gives us a complex subspace of $\mathfrak g_\mathbb C$. It contains $\mathfrak h_\mathbb C=\ker\pi$. It is furthermore $\Ad(H)$-invariant, since $\Ad(h)(\pi^{-1}(\mathfrak p))=\pi^{-1}(\tilde{\Ad(h)}\mathfrak p)$. Since complex conjugation commutes with $\pi$, we also get $\mathfrak b\cap\overline{\mathfrak b}=\pi^{-1}(0)=\mathfrak h_\mathbb C$. The dimension property is also satisfied, since $\mathfrak b/\mathfrak h_\mathbb C\simeq\mathfrak p$ and the latter has half the dimension of $(\mathfrak g/\mathfrak h)_\mathbb C$. It remains to show that $\mathfrak b$ is a subalgebra of $\mathfrak g_\mathbb C$. For that consider the following proposition.
		\end{proof}
		
		For the following proposition, we consider any homogeneous $G$-space $M:=G/H$ with projection $\pi$. We denote by $\rho$ the right-translation on $G$.
	
	\begin{prop} \label{prop:intofbundleliealg}
		Let $\mathfrak b\subseteq\mathfrak g_\mathbb C$ an $\Ad(H)$-invariant subspace with $\mathfrak h_\mathbb C\subseteq\mathfrak b$ and denote $\mathfrak p:=\mathfrak b/\mathfrak h_\mathbb C\subset T_{eH}M_\mathbb C$. Let $\mathcal P\subseteq TM_\mathbb C$ be the $G$-invariant subbundle having $\mathcal P_{eH}=\mathfrak p$. Then
		$$\mathcal P\text{ is integrable}\quad\Leftrightarrow\quad\mathfrak b\text{ is a subalgebra}.$$
	\end{prop}
	
	We do some preparational work for the proof. We denote by $\mathcal X_\mathbb C(G)^{\rho(H)}$ the complex right-$H$-invariant vector fields on $G$. That is, $\mathfrak X\in\mathcal X_\mathbb C(G)^{\rho(H)}$ satisfies
	$$\mathfrak (Xf)(gh)=\mathfrak X(f\circ\rho_h)(g)\quad\forall g\in G,h\in H,f\in C^\infty(G).$$
	The derivative of the projection $d\pi$ pointwisely maps tangent space of $G$ into tangent spaces of $M$. If we want define the \textbf{push-forward} of a vector field $\mathfrak X$ on $G$ via $d\pi$, we must for the well-definedness demand that
	$$\mathfrak X(f\circ\pi)(g)=\mathfrak X(f\circ\pi)(gh)\quad\forall g\in G,h\in H,f\in C^\infty(M).$$
	Vector fields with that property are called \textbf{projectable}. We find that right-$H$-invariant vector fields are projectable, since $\pi=\pi\circ\rho_h$ for all $h\in H$. For projectable vector fields $\mathfrak X$, we define the push-forward map $\pi_*$ by
	$$\pi_*\mathfrak Xf(gH):=\mathfrak X(f\circ\pi)(g)\quad\forall gH\in M, f\in C^\infty(M).$$
	
	\begin{lemm} \label{lemm:existofrHipreimage}
		For each $\mathfrak V\in\mathcal X_\mathbb C(M)$ and each $x\in M$ there is $x\in U$ open and $\mathfrak X\in\mathcal X_\mathbb C(G)$ right $H$-invariant such that
		$$\pi_*\mathfrak X_{|U}=\mathfrak V_{|U}.$$
	\end{lemm}
	
	\begin{proof}
		From Proposition we know that on a neighborhood $U\subseteq M$ of $x$ there is a section $s:U\to G$. If we denote $V:=\pi^{-1}(U)$, then $V\simeq U\times H$ via $(x,h)\mapsto s(x)h$. We also observe that $V$ is closed under right-translation by $H$. We can then define a right-$H$-invariant vector field $\mathfrak X$ on $V$ from $\mathfrak V$ by the formula
		$$\mathfrak Xf(s(x)h):=\mathfrak V(f\circ s)(x)\quad\forall s(x)h\in V,f\in C^\infty(V).$$
		It is clear that $\mathfrak X$ is right-$H$-invariant (since the section corresponds to a fixed choice of an element in $s(x)H$). Furthermore on $U$ we have $\pi_*\mathfrak X=\mathfrak V$. Using cut-off functions we can define a vector field on $G$ which coincides with $\mathfrak X$ on a smaller set and hence has the desired properties.
	\end{proof}
	
	\begin{lemm}
		The set of right-$H$-invariant vector fields on $G$ is closed under the Lie bracket.
	\end{lemm}
	
	\begin{proof}
		Suppose $\mathfrak X$ and $\mathfrak Y$ right-$H$-invariant. Then
		\begin{align*}
			\mathfrak X(\mathfrak Yf)(gh) & = \mathfrak X((\mathfrak Yf)\circ \rho_h)(g) \\
																		& = \mathfrak X(\mathfrak Y(f\circ\rho_h))(g).
		\end{align*}
		Doing the same for $\mathfrak Y\mathfrak X$ we find that $[\mathfrak X,\mathfrak Y]$ is also right-$H$-invariant.
	\end{proof}
	
	\begin{lemm} \label{lemm:projishomo}
		Let $\mathfrak X$ and $\mathfrak Y$ on $G$ be right-$H$-invariant. Then
		$$\pi_*[\mathfrak X,\mathfrak Y]=[\pi_*\mathfrak X,\pi_*\mathfrak Y].$$
	\end{lemm}
	
	\begin{proof}
		\begin{align*}
			(\pi_*[\mathfrak X,\mathfrak Y]f)(gH) & = (\mathfrak X\mathfrak Y-\mathfrak Y\mathfrak X)(f\circ\pi)(g) \\
																						& = \mathfrak X(\mathfrak Y(f\circ\pi))(g)-\mathfrak Y(\mathfrak X(f\circ\pi))(g) \\
																						& = (\pi_*\mathfrak X)(\mathfrak Y(f\circ\pi)(g))-(\pi_*\mathfrak Y)(\mathfrak X(f\circ\pi)(g))\\
																						& = ((\pi_*\mathfrak X)(\pi_*\mathfrak Y)-(\pi_*\mathfrak Y)(\pi_*\mathfrak X))f(gH)\\
																						& = [\pi_*\mathfrak X,\pi_*\mathfrak Y]f(gH).
		\end{align*}
	\end{proof}
	
	We denote by $\mathcal B$ the preimage of $\mathcal P$ under the map $d\pi:TG_\mathbb C\to TM_\mathbb C$.
	
	\begin{lemm}
		The subbundle $\mathfrak B$ is left-$G$-invariant und right-$H$-invariant.
	\end{lemm}
	
	\begin{proof}
		The left-$G$-invariance follows from the equivariance of $\pi$ with respect to the left-action of $G$ on $G$ and $M$. For the right-$H$-invariance, we need to show that
		$$\rho_h^\prime(g)(\mathcal B_g)=\mathcal B_{gh}\quad\forall g\in G,\ \forall h\in H.$$
		By taking the definiton of $\mathcal B$ this reduces to showing for $v\in T_{gh}G_\mathbb C$
		$$d\pi^\prime(gh)v\in\mathcal P_{gH}\quad\Leftrightarrow\quad d\pi^\prime(g)(\rho_{h^{-1}}^\prime(gh))v\in\mathcal P_{gH}.$$
		But this follows from the chain rule applied to $\pi=\pi\circ\rho_{h}$ which gives
		$$d\pi^\prime(gh)=d\pi^\prime(g)\circ\rho_{h^{-1}}^\prime(gh).$$
	\end{proof}
	
	\begin{lemm} \label{lemm:subbspannedbyrHi}
		The subbundle is locally spanned by right-$H$-invariant vector fields.
	\end{lemm}
	
	\begin{proof}
		Let $g\in G$ choose $U$ neighborhood of $gH$ and section $s:U\to G$ with $s(gH)=g$ and $\pi\circ s=\Id_U$. Let $X_1,\dots,X_s$ be a basis of $\mathcal B_g$. Choose smooth sections $\mathfrak X_1,\dots,\mathfrak X_s$ of $\mathcal B$ along s(U) which have these values at $g$. By making $U$ smaller, we can assume that the vector fields $\mathfrak X_1,\dots,\mathfrak X_s$ form a basis of $\mathcal B$ in each point of $s(U)$. Now we can extend the vector fields uniquely to right-$H$-invariant vector fields on $s(U)H$ as it was done in Lemma \ref{lemm:existofrHipreimage}. Since $\mathcal B$ is right-H-invariant and the derivatives of $\rho(H)$ are isomorphisms, the extensions span $\mathcal B$ at every point of $s(U)H$.
	\end{proof}
	
	\begin{prop} \label{prop:intifint}
		The distribution $\mathcal B$ is integrable if and only if $\mathcal P$ is integrable.
	\end{prop}
	
	\begin{proof}
		Suppose $\mathcal B$ is integrable and let $\mathfrak V$ and $\mathfrak W$ be sections of $\mathcal P$. For any point in $M$ we can choose a neighborhood $U$ and right-$H$-invariant vector fields $\mathfrak X$ and $\mathfrak Y$ on $G$ such that $\pi_*\mathfrak X_{|U}=\mathfrak V_{|U}$ and $\pi_*\mathfrak Y_{|U}=\mathfrak W_{|U}$. Then Lemma \ref{lemm:projishomo} gives us
		$$\pi_*[\mathfrak X,\mathfrak Y]_{|U}=[\mathfrak V,\mathfrak W]_{|U}.$$
		Since $\mathfrak X$ and $\mathfrak Y$ are sections of $\mathcal B$, the right-hand-side is $\mathcal P$-valued. Therefore $[\mathfrak V,\mathfrak W]$ is again a section of $\mathcal P$. Hence $\mathcal P$ is integrable.
		
		Suppose $\mathcal P$ is integrable and let $\mathfrak X$ and $\mathfrak Y$ be sections of $\mathcal B$. Then according to Lemma \ref{lemm:subbspannedbyrHi} around each $g\in G$ these two vector fields are sums of right-$H$-invariant vector fields, say $\mathfrak X=\sum_{i=1}^s a_i\mathfrak X_i$ and $\mathfrak Y=\sum_{j=1}^s b_j\mathfrak Y_j$. Then their bracket is given by the term
		\begin{align*}
			[\mathfrak X,\mathfrak Y] & = \sum_{i,j=1}^s a_i((\mathfrak X_i)b_j\cdot\mathfrak Y_j+b_j(\mathfrak X_i\circ\mathfrak Y_j))-b_j((\mathfrak Y_j)a_i\cdot\mathfrak X_i-a_i(\mathfrak Y_j\circ\mathfrak X_i)) \\
																& = \sum_{i,j=1}^s a_ib_j\cdot[\mathfrak X_i,\mathfrak Y_j] + \sum_{i,j=1}^s a_i(\mathfrak X_i b_j)\cdot\mathfrak Y_j - \sum_{i,j=1}^s b_j(\mathfrak Y_j a_i)\cdot\mathfrak X_i.
		\end{align*}
		We apply Lemma \ref{lemm:projishomo} to $[\mathfrak X_i,\mathfrak Y_j]$ and use the integrability of $\mathfrak P$ to find that $\pi_*[\mathfrak X_i,\mathfrak Y_j]$ is a local section of $\mathcal P$. Hence $[\mathfrak X_i,\mathfrak Y_j]$ is $\mathcal B$-valued. Therefore $[\mathfrak X,\mathfrak Y]$ is a section of $\mathcal B$.
	\end{proof}
	
	\begin{proof}[Proof of Proposition \ref{prop:intofbundleliealg}]
		By the left-invariance of $\mathcal B$, it is (globally) spanned by left-invariant vector fields $\widetilde{X_1},\dots\widetilde{X_s}$, where $X_1,\dots,X_s$ is a basis of $\mathfrak b$. If $\mathfrak b$ is a Lie algebra, then $[\widetilde{\mathfrak X_i},\widetilde{\mathfrak X_j}]=\widetilde{[X_i,X_j]}$ is a again a section of $\mathcal B$ (c.f. Definition \ref{defi:liealg}). The integrabililty of $\mathfrak B$ and by Proposition \ref{prop:intifint} the integrabilty of $\mathcal P$ follows by writing sections of $\mathcal B$ as sums and doing the same calculation as in the previous proof.
		
		If on the other hand $\mathcal P$ is integrable, then $\mathcal B$ is integrable. Then if $X$ and $Y$ are in $\mathfrak b$, we have $[\tilde{X},\tilde{Y}](e)\in\mathcal B_e=\mathfrak b$. Hence $\mathfrak b$ is a subalgebra of $\mathfrak g_\mathbb C$.
	\end{proof}
		
		\begin{theo} \label{theo:eopocot}
			Let $M=G/C(T_1)$ be a generalized flag manifold and $T_1\subseteq T$ a maximal torus. Then to each $T_1$-admissible Weyl chamber of $T$, there exists a subalgebra $\mathfrak b$ of $\mathfrak g_\mathbb C$ satisfying the assumptions of Proposition \ref{prop:csalg}.
		\end{theo}
		
		\begin{proof}
			Let $\Delta_1^+$ be a positive root system of $T_1$ corresponding to a $T_1$-admissible Weyl chamber as in Proposition \ref{prop:existwc}. Define $\mathfrak c_{1\mathbb C}:=\mathfrak t_\mathbb C+\sum_{\alpha\in\Delta_1}\mathfrak g^\alpha$ as in Construction \ref{cons:altdefcent} and consider
			$$\mathfrak b:=\mathfrak c_{1\mathbb C}\oplus\sum_{\alpha\in(\Delta^+\backslash\Delta_1^+)}\mathfrak g^\alpha.$$
			From Proposition \ref{prop:rootprop} we see that $\mathfrak b$ is closed under the Lie bracket and hence it is a subalgebra of $\mathfrak g_\mathbb C$. Furthermore, its complex conjugate $\overline{\mathfrak b}$ is given by $\mathfrak c_{1\mathbb C}\oplus\sum_{\alpha\in(\Delta^+\backslash\Delta_1^+)}\mathfrak g^{-\alpha}$, because $\mathfrak c_{1\mathbb C}$ is invariant under complex conjugation and $\overline{\mathfrak g}_\alpha=\mathfrak g^{-\alpha}$. Therefore $\mathfrak b\cap\overline{\mathfrak b}=\mathfrak c_{1\mathbb C}$ and $2\dim\mathfrak b/\mathfrak c_{1\mathbb C}=\dim\mathfrak g_\mathbb C/\mathfrak c_{1\mathbb C}$. Furthermore, we have $\ad(\mathfrak c_{1\mathbb C})\mathfrak b=[\mathfrak c_{1\mathbb C},\mathfrak b]\subseteq\mathfrak b$. Hence for $H\in\mathfrak c_1=\mathfrak c_{1\mathbb C}\cap\mathfrak g$ holds $\Ad(\exp H)\mathfrak b=\exp(\ad H)\mathfrak b\subseteq\mathfrak b$ and we find that $\mathfrak b$ is $\Ad(C(T_1))$-invariant.
		\end{proof}
		
		\begin{coro} \label{coro:gfmcomplex}
			Any generalized flag manifold is a complex manifold with a natural $G$-invariant integrable almost complex structure.
		\end{coro}
		
	\subsection{Coadjoint Orbits} \label{sect:coadorb}
		
		We have given a definition of coadjoint orbits in Example \ref{exam:co}. Furthermore, if $\lambda\in\mathfrak g^*$ and
		$$\mathcal O_\lambda=\{\Ad^*(g)\lambda\ |\ g\in G\}$$
		is the associated coadjoint orbit\index{coadjoint!orbit}\index{$\mathcal O_\lambda$, coadjoint orbit of $\lambda$}, then the coadjoint action of $G$ on $\mathcal O_\lambda$ is transitive. Because of Corollary \ref{coro:fots} we can assume that $\lambda\in\mathfrak t^*$.
		
		We denote by $G_\lambda$\index{$G_\lambda$, stabilizer of $\lambda$} the stabilizer of $\lambda$ with respect to the coadjoint action. Proposition \ref{prop:stab} tells us that $G_\lambda$ is a Lie group with Lie algebra $\mathfrak g_\lambda=\{X\in\mathfrak g\ |\ \ad^*(X)\lambda=0\}$\index{$\mathfrak g_\lambda$, Lie algebra of $G_\lambda$}.
		
		Let $f\colon G\to\mathcal O_\lambda$ be the orbit map which maps $g\in G$ to the point $g\cdot\lambda\in\mathcal O_\lambda$. Then the induced $G$-equivariant bijection
		$$F\colon G/G_\lambda\to\mathcal O_\lambda, \quad gG_\lambda\mapsto \Ad^*(g)\lambda$$
		is exactly the map from Theorem \ref{theo:homoisquot}. Of course we have not yet shown that $\mathcal O_\lambda$ is a smooth manifold and hence $F$ is not a diffeomorphism. Nevertheless we know from Theorem \ref{theo:quotishomo} that $G/G_\lambda$ is a homogeneous $G$-space and can use $F$ to transport the differentiable structure from $G/G_\lambda$ onto $\mathcal O_\lambda$.
		
		We will sketch how the differentiable structure on $G/G_\lambda$ is obtained. For a detailed proof see \cite[Theorem 3.58]{war} or \cite[§8.4]{hga}.
		
		We have stated in Proposition \ref{prop:expprop} that $\exp\colon G\to G$ has derivative $\exp^\prime(0)=\Id_\mathfrak g$. With this property one can prove for the multiplication $m\colon G\times G\to G$ on $G$ that
		\begin{equation} \label{eq:derivmult}
			m^\prime(e,e)(X,Y)=X+Y\quad\forall X,\ Y\in G.
		\end{equation}
		We fix a vector space complement $\mathfrak a$ of $\mathfrak g_\lambda$ in $\mathfrak g$ such that $\mathfrak g=\mathfrak a\times\mathfrak g_\lambda$ and consider the map
		$$\theta\colon\mathfrak a\times\mathfrak g_\lambda\to G,\quad (X,Y)\mapsto \exp X \exp Y.$$
		Then (\ref{eq:derivmult}) shows that the derivative at $(0,0)$ is invertible and by the Inverse Function Theorem $\theta$ is a homeomorphism on a small neighborhood. Then we can choose an open subset $U$ of $\mathfrak a$ for which the map
		$$f\circ\theta\colon U\times\{0\}\to\mathcal O_\lambda$$
		is a homeomorphism onto its image. Moving the image around via translation yields an atlas for which $\mathcal O_\lambda$ is a smooth manifold with the following properties (cf. \cite[Theorem 8.4.6]{hga}).
		
		\begin{theo} \label{theo:hilgcoadorb}  \mbox{}
			\renewcommand{\theenumi}{(\roman{enumi})}
			\renewcommand{\labelenumi}{\theenumi}
			\begin{enumerate}
				\item $f\colon G\to\mathcal O_\lambda$ is a submersion and $\mathcal O_\lambda\to\mathfrak g^*$ is an immersion.
				\item $T_\lambda\mathcal O_\lambda=(\ad^*\mathfrak g)\lambda\stackrel{\text{via }F}{=}\mathfrak g/\mathfrak g_\lambda$.
				\item $T_{\Ad^*(g)\lambda}\mathcal O_\lambda=(\ad^*\mathfrak g)(\Ad^*(g)\lambda)\stackrel{\text{via }F}{=}\mathfrak g/\mathfrak g_{(\Ad^*(g)\lambda)}$.
			\end{enumerate}
		\end{theo}
		
	\subsection{Stabilizers of Functionals} \label{ss:stabofunc}
		
		Let $G$ be a compact connected Lie group and $T$ a maximal torus. We consider a functional $\lambda\in\mathfrak t^*$.
		
		We have stated in Subsection \ref{structtheo} that we view $\mathfrak t^*$ as a subset of $\mathfrak g^*$ via extending functionals by zero on the orthogonal complement of $\mathfrak t$ with respect to an $\Ad(G)$-invariant inner product. Let $t\in T$. The root space decomposition in Construction \ref{cons:rsd} tells us $\Ad(t)_{|\mathfrak t}=\Id_\mathfrak t$ and $\Ad(t)\mathfrak t^\bot\subseteq\mathfrak t^\bot$. Hence
		$$\Ad^*(t)\lambda(X+Y)=\lambda(X+\Ad(t^{-1})Y)=\lambda(X+Y)\quad\forall X\in\mathfrak t,\ \forall Y\in\mathfrak t^\bot.$$
		Therefore we obtain $T\subseteq G_\lambda$.
		
		\begin{defi}
			If $T$ is a maximal torus, we call a functional $\lambda\in\mathfrak t^*$ \textbf{regular}\index{functional!regular}, if $T=G_\lambda$. Else it is called \textbf{singular}\index{functional!singular}.
		\end{defi}
		
		\begin{prop} \label{prop:stabisconn}
			The stabilizer $G_\lambda$ is connected.
		\end{prop}
		
		\begin{proof} 
			This is clear for regular $\lambda$. For the general case, we show that the stabilizer $G_X$ of an element $X\in\mathfrak g$ with respect to the adjoint action is connected. This is sufficient since the adjoint and coadjoint action can be identified by using an $\Ad(G)$-invariant inner product.
			\begin{align*}
				G_X & = \{g\in G\ |\ \Ad(g)X=X\} \\
						& = \{g\in G\ |\ \Ad(g)tX=tX\quad\forall t\in\mathbb R\} \\
						& = \{g\in G\ |\ \exp(\Ad(g)tX)=\exp(tX)\quad\forall t\in\mathbb R\}
			\end{align*}
			To obtain equality in the last step, we need to see equivalence of the conditions. One direction is obvious. For the other, assume that
			\begin{align*}
				e & = \exp(\Ad(g)tX)(\exp(tX))^{-1} \\
					& = \exp(t(\Ad(g)X-X))
			\end{align*}
			Then $\{t(\Ad(g)X-X)\ |\ \exp(t(\Ad(g)X-X))=e,\  t\in\mathbb R\}$ is a linear subspace which has $0$ as an isolated point (because $\exp$ is locally invertible). Thus it must be zero and we retrieve the condition $\Ad(g)tX=tX$. From here we continue the calculation:
			\begin{align*}
				G_X & = \{g\in G\ |\ \exp(\Ad(g)tX)=\exp(tX)\quad\forall t\in\mathbb R\} \\
						& = \{g\in G\ |\ g \exp(tX)g^{-1}=\exp(tX)\forall t\in\mathbb R\} \\
						& = C(\exp(tX\ |\ t\in\mathbb R))
						& = C(\overline{\exp(tX\ |\ t\in\mathbb R)})
						& = C(T_1),
			\end{align*}
			where $T_1=\overline{\exp(tX\ |\ t\in\mathbb R)}$ is a torus. Thus $G_X$ is connected by Proposition \ref{prop:maxtcent}.
		\end{proof}
		
		We will find another criterion for $\lambda$ to be regular or singular. Recall that the Lie algebra of $G_\lambda$ is given by
		$$\mathfrak g_\lambda=\{X\in\mathfrak g\ |\ \lambda[X,Y]=0\quad\forall Y\in\mathfrak g\}.$$
		
		Let $(.,.)$ denote the inner product on $i\mathfrak t^*$ defined by the $\Ad(G)$-invariant inner product on $\mathfrak g$ and duality. Then Proposition \ref{prop:existrootvector} asserts the existence of non-zero elements $X_\alpha$, $\alpha\in\Delta$, such that
		$$\xi([X_\alpha,X_{-\alpha}])=(\xi,\alpha)\quad\forall\xi\in i\mathfrak t^*.$$
		Proposition \ref{prop:rootprop} tells us that the root spaces are one-dimensional. Hence an arbitrary element $X$ of $\mathfrak g_\mathbb C$ is of the form
		$$X=X_\mathfrak t+\sum_{\alpha\in\Delta}c_\alpha X_\alpha.$$
		Here $X_\mathfrak t=X_\mathfrak t(X)$ is a suitable element of $\mathfrak t_\mathbb C$ and the $c_\alpha$ are scalars depending on $X$.
		
		We are ready to determine $\mathfrak g_\lambda$. For the following calculation we use Proposition \ref{prop:rootprop} and that $\lambda_{|g_\alpha}=0$ for all $\alpha\in\Delta$. Pick $\beta\in\Delta$, then
		\begin{align*}
			\lambda([X,X_\beta]) & = \lambda([X_\mathfrak t+\sum_{\alpha\in\Delta}c_\alpha X_\alpha,X_\beta]) \\
													 & = c_{-\beta}\lambda([X_{-\beta},X_\beta]) \\
													 & = -i c_{-\beta}(i\lambda,\beta).
		\end{align*}
		This calculation shows that $X_\beta\in\mathfrak g_{\lambda\mathbb C}$ if and only if $(i\lambda,\beta)=0$. Since the root space $\mathfrak g^\beta$ is one-dimensional and hence spanned by $X_\beta$, in this case we have $\mathfrak g^\beta\subset\mathfrak g_{\lambda\mathbb C}$. Since the criterion is independent of the sign of $\beta$, we then also have $\mathfrak g^{-\beta}\subset\mathfrak g_{\lambda\mathbb C}$. We define the \textbf{singular roots}\index{root!singular}\index{$\Delta_{sing}(\lambda)$, singular roots of $\lambda$} of $\lambda$ as
		$$\Delta_{sing}(\lambda):=\{\beta\in\Delta\ |\ (i\lambda,\beta)=0\}$$
		and $\Delta_{sing}^+(\lambda):=\Delta^+\cap\Delta_{sing}$. Then we can express our consideration in the following way.
		
		\begin{prop} \label{prop:regsing}
			Let $T$ be a maximal torus in $G$ and $\lambda\in\mathfrak t^*$. Then
			$$\mathfrak g_{\lambda\mathbb C}=\mathfrak t_\mathbb C\oplus\sum_{\beta\in\Delta_{sing}(\lambda)}\mathfrak g^\beta$$
			and
			$$\mathfrak g_\lambda=\mathfrak t\oplus\sum_{\beta\in\Delta_{sing}^+(\lambda)} \mathfrak g\cap(\mathfrak g^\beta+\mathfrak g^{-\beta}).$$
			Therefore $\lambda$ is regular if and only if $\Delta_{sing}(\lambda)=\emptyset$.
		\end{prop}
		
		Let $\lambda\in\mathfrak t^*$ and $\alpha,\ \beta\in\Delta_{sing}(\lambda)$. Then
		$$(i\lambda,\alpha+\beta)=(i\lambda,\alpha)+(i\lambda,\beta)=0.$$
		If $\alpha+\beta\in\Delta$, then $\alpha+\beta\in\Delta_{sing}(\lambda)$ and therefore we can apply Construction \ref{cons:altdefcent} with $\Delta_1=\Delta_{sing}(\lambda)$. This gives us a torus $T_1\subseteq T$ such that $C(T_1)$ has Lie algebra $\mathfrak g_\lambda$. Since $G_\lambda$ is connected after Proposition \ref{prop:stabisconn}, we have $G_\lambda=C(T_1)$ because of Proposition \ref{theo:existsubgroups}. We have thus proved the following important proposition.
		
		\begin{prop} \label{prop:stabiscent}
			Let $G$ be a compact connected Lie group, $T$ a maximal torus and $\lambda\in\mathfrak t^*$. Then there exists a torus $T_1$ such that $G_\lambda=C(T_1)$. If $\lambda$ is regular, then $T_1=T$ and hence $G_\lambda=T$.
		\end{prop}
		
		\begin{coro} \label{coro:coadgfm}
			Let $G$ be a compact connected Lie group, $T$ a maximal torus and $\mathcal O_\lambda$ a coadjoint orbit. We can assume $\lambda\in\mathfrak t^*$ because of Corollary \ref{coro:fots}. From Subsection \ref{sect:coadorb} we know that $\mathcal O_\lambda$ has the form $G/G_\lambda$. If we apply Proposition \ref{prop:stabiscent} to find a torus $T_1$ such that $G_\lambda=C(T_1)$, then we have shown that any coadjoint orbit $\mathcal O_\lambda$ has the structure of a generalized flag manifold $G/C(T_1)$.
		\end{coro}
		
	\chapter{Representation Theory}
	
	\section{Complex Line Bundles}
	
	\subsection{Line Bundles with Hermitian Connection}
		
		\begin{defi} 
			Let $M$ be a smooth manifold. A \textbf{complex line bundle}\index{line bundle!complex} $L$ over $M$ is a smooth manifold $L$ together with a smooth surjective map $\pi\colon L\to M$ such that
			\renewcommand{\theenumi}{(\roman{enumi})}
			\renewcommand{\labelenumi}{\theenumi}
			\begin{enumerate}
				\item for each $p\in M$ the \textbf{fiber}\index{fiber} $L_p:=\pi^{-1}(p)$\index{$L_p$, fiber over $p$} is a one-dimensional complex vector space,
				\item for each $p\in M$ there is a neighborhood $U$ of $p$ and a diffeomorphism $\varphi$ of the form
					\begin{equation} \label{eq:vb}
						\varphi\colon\pi^{-1}(U)\to U\times\mathbb C,\quad\varphi(q)=(\pi(q),\phi(q))
					\end{equation}
					such that $\phi_{|L_p}\colon L_p\to\mathbb C$ is a complex-linear isomorphism for each fiber.
			\end{enumerate}
			 A pair $(U,\varphi)$ is called a \textbf{local trivialization}. If $V\subseteq M$ is open then a map $s\colon V\to L$ satisfying $s(p)\in L_p$ for all $p\in V$ is called a \textbf{section}\index{section} on $V$. It is called \textbf{non-vanishing} if $s(p)\neq0\in L_p$ for all $p\in V$. The $C^\infty(V,\mathbb C)$-module of smooth sections on $V$ is denoted by $C^\infty(V,L)$\index{$C^\infty(V,L)$, smooth sections on $V$ of $L$}.	Let $\{U_j,\varphi_j\}_{j\in J}$ be a covering by local trivializations. For each $j\in J$ we obtain a non-vanishing smooth section $s_j$ on $U_j$ defined by the formula
			$$s_j(p)=\phi_j^{-1}(1).$$
			Then the family $\{U_j,s_j\}_{j\in J}$ is called a \textbf{local system} for $L$. By the corresponding \textbf{transition functions}\index{transition functions} we mean the family of functions $c_{jk}\in C^\infty(U_j\cap U_k,\mathbb C)$\index{$c_{jk}$, transition function}, for $j,\ k$ such that the intersection is not empty, defined by
			$$c_{jk}(q)z=\phi_k^{-1}\circ\phi_j(q)z.$$
			From that we find for the local system
			$$c_{jk}s_j=s_k.$$
			The transition functions satisfy the following \textbf{cocycle conditions}\index{cocycle conditions}
			\begin{equation} \label{eq:cc}
				c_{jk}=c_{kj}^{-1},\quad c_{jk}c_{kl}=c_{jl},
			\end{equation}
			where in the latter case $j,k,l$ are such that $U_j\cap U_k\cap U_l\neq\emptyset$.
		\end{defi}
		
		\begin{defi} 
			Let $L_1$ and $L_2$ be complex line bundles over a smooth manifold $M$. We say that $L_1$ and $L_2$ are \textbf{equivalent}\index{line bundle!equivalence} if there exists a diffeomorphism $\tau:L_1\rightarrow L_2$ such that for all $p\in M$ the map
			\begin{equation} \label{eq:eolb}
				\tau_{|(L_1)_p}:(L_1)_p\rightarrow(L_2)_p
			\end{equation}
			is a complex-linear isomorphism. This defines an equivalence relation on the set of complex line bundles over $M$. If $L$ is a complex line bundle over $M$, then we denote by $\left[L\right]$ the equivalence class of $L$. The set of all equivalence classes is denoted by $\mathcal L(M)$\index{$\mathcal L(M)$, equivalence classes of line bundles}.
		\end{defi}
		
		\begin{prop} \label{prop:tfgiveclb}
			Let $\{U_j\}_{j\in J}$ be an open covering of $M$ and assume that for each pair $j,\ k$ such that $U_j\cap U_k$ is not empty, we have a function $c_{jk}\in C^\infty(U_j\cap U_k)$ such that if $U_j\cap U_k\cap U_l$ is not empty the cocycle conditions (\ref{eq:cc}) hold. Then up to equivalence there is a unique complex line bundle with a local system $\{U_j,s_j\}_{j\in J}$ which has transition functions $c_{jk}$.
		\end{prop}
		
		\begin{proof}[Sketch of proof]
			The complex line bundle can be realized by taking the disjoint union $\bigcup_{j\in J}U_j\times\mathbb C$ and factor out the equivalence relation
			$$(p_j,z_j)\simeq(p_k,z_k)\quad\Leftrightarrow\quad p_j=p_k\text{ and }z_j=c_{jk}(p_j)z_k,$$
			where $(p_j,z_j)\in U_j\times\mathbb C$ and $(p_k,z_k)\in U_k\times\mathbb C$. Confer \cite[Theorem 1.2.6]{walhomo}.
		\end{proof}
		
		\begin{prop} \label{prop:eqiftfrel}
			Let $L_1$ and $L_2$ be complex line bundles with local systems $\{U_j,s_j^1\}_{j\in J}$ and $\{U_j,s_j^2\}_{j\in J}$ and transition functions $c_{jk}^1$, $c_{jk}^2$, respectively. Then $L_1$ and $L_2$ are equivalent if and only if for each $j$ there are smooth functions  $g_j\colon U_j\to\mathbb C^*$ such that $g_k c_{jk}^1 g^{-1}_j=c_{jk}^2$ on each non-empty intersection $U_j\cap U_k$.
		\end{prop}
		
		\begin{proof}
			\cite[Proposition 1.2.7]{walhomo}.
		\end{proof}
		
		\begin{defi}
			Let $M$ be a smooth manifold. Then a \textbf{line bundle with connection}\index{line bundle!with connection} over $M$ is a pair $(L,\nabla)$ such that $L$ is a complex line bundle over $M$ and $\nabla$, which is called \textbf{connection}\index{connection} or \textbf{covariant derivative}\index{covariant derivative}, is an assignment\index{$\nabla$, affine connection}
			$$\nabla\colon\mathcal X_{\mathbb C}(M)\to\End_\mathbb C(C^\infty(M,L))$$
			which satisfies
			\renewcommand{\theenumi}{(\roman{enumi})}
			\renewcommand{\labelenumi}{\theenumi}
			\begin{enumerate}
				\item $\nabla_{f\mathfrak X+g\mathfrak Y}=f\nabla_\mathfrak X+g\nabla_\mathfrak Y\quad (C^\infty(M,\mathbb C)\text{-linearity})$,
				\item $\nabla_\mathfrak X fs=\mathfrak X(f)s+f\nabla_\mathfrak X s\quad(\text{Leibniz rule})$
			\end{enumerate}
			for all functions $f,\ g\in C^\infty(M,\mathbb C)$, vector fields $\mathfrak X,\ \mathfrak Y\in\mathcal X_\mathbb C(M)$, and sections $s\in C^\infty(M,L)$.
		\end{defi}
		
		\begin{defi} 
			Let $(L_1,\nabla^1)$ and $(L_2,\nabla^2)$ be complex line bundles with connection over a smooth manifold $M$. We say that $(L_1,\nabla^1)$ and $(L_2,\nabla^2)$ are \textbf{equivalent as line bundles with connection}\index{line bundle!equivalence with connection}, if there exists a diffeomorphism $\tau:L_1\rightarrow L_2$ satisfying the property (\ref{eq:eolb}) of an equivalence of line bundles and in addition
			$$\tau\circ(\nabla_{\mathfrak X}^1s)=\nabla_{\mathfrak X}^2(\tau\circ s)\quad\forall s\in C^\infty(M,L_1)\quad\forall\mathfrak X\in\mathcal X_{\mathbb C}(M).$$
		\end{defi}
		
		\begin{prop} \label{prop:covderisloc} 
			Let $p\in M$, then the value $(\nabla_\mathfrak Xs)(p)$ depends only on the values of $s$ and $\mathfrak X$ in an arbitrarily small neighborhood of $p$.
		\end{prop}
		
		\begin{proof}
			Let $U$ be open and $p\in U$. Suppose $s=0$ on $U$. Using Theorem \ref{theo:countpara} we can find an open neighborhood $p\in V\subseteq U$ and a smooth function $f$ such that $f=0$ on $V$ and $f=1$ on $M \backslash U$. Then $s=fs$ and we find
			\begin{align*}
				\nabla_\mathfrak X s(p) & = \nabla_\mathfrak X f s(p) \\
														 		& = \mathfrak X(f)(p) s(p) +f(p) \nabla_\mathfrak X s(p) = 0.
			\end{align*}
			Hence any two section $s,\ t$ which coincide on a small neighborhood of $p$ have $\nabla_\mathfrak X s(p)=\nabla_\mathfrak X t(p)$ for all $\mathfrak X\in\mathcal X_\mathbb C(M)$. The proof for $\mathfrak X$ is similar.
		\end{proof}
		
		\begin{prop} \label{prop:smoothlocsect}
			Let $L$ be a complex line bundle over $M$. Let $U\subseteq M$ be open and let $s\in C^\infty(U,L)$ a non-vanishing smooth section on $U$. If $t\in C^\infty(U,L)$ is another smooth section on $U$, then the function $f\colon U\rightarrow\mathbb C$ defined by
			$$f(p) s(p)= t(p)\quad\forall p\in M$$
			is smooth. In particular, any smooth section on $U$ is given by the product of $s$ with a smooth function. We will use $\frac{t}{s}$ as a symbol for $f$.
		\end{prop}
		
		\begin{proof}
			$f$ is smooth since on a local trivialization it is the quotient of two smooth functions. Since $s$ is non-vanishing, the quotient is well-defined for any section $t$.
		\end{proof}
		
		\begin{cons} \label{cons:fam1forms}
			Let $\{U_j,s_j\}_{j\in J}$ be a local system of $L$ and consider a fixed $j$. Because of Proposition \ref{prop:covderisloc} the restriction of $\nabla$ to $U_j$ is well-defined. We can define a 1-form $\alpha_j\colon\mathcal X_\mathbb C(U_j)\to C^\infty(U_j,\mathbb C)$\index{$\alpha_j$, local form} by the formula
			$$\left\langle \alpha_j,\mathfrak X\right\rangle=\frac{1}{2\pi i}\frac{\nabla_\mathfrak X s_j}{s_j},$$
			where we have used the result and the notation from Proposition \ref{prop:smoothlocsect}.
		\end{cons}
		
		\begin{prop} \label{prop:proplocform}
			Fix $j,\ k$ in a local system such that $U_j\cap U_k\neq\emptyset$ and apply Construction \ref{cons:fam1forms}. Then
			\begin{equation} \label{eq:proplocform}
				\alpha_k=\alpha_j+\frac{1}{2\pi i}\frac{dc_{jk}}{c_{jk}}.
			\end{equation}
		\end{prop}
		
		\begin{proof}
			Let $\mathfrak X\in\mathcal X_\mathbb C(M)$. Then
			\begin{align*}
				2\pi i\left\langle \alpha_j,\mathfrak X\right\rangle & = \frac{\nabla_\mathfrak X s_j}{s_j} \\
																														 & = \frac{\nabla_\mathfrak X (c_{jk}s_k)}{c_{jk}s_k} \\
																														 & = \frac{\mathfrak X(c_{jk})}{c_{jk}}+\frac{\nabla_\mathfrak X s_k}{s_k}\\
																														 & = \left\langle \frac{dc_{jk}}{c_{jk}}+2\pi i\alpha_k,\mathfrak X\right\rangle.
			\end{align*}
		\end{proof}
		
		We call the family $\{U_j,s_j,\alpha_j\}_{j\in J}$ a \textbf{local form}\index{connection!local form} of $\nabla$. Suppose that $s\in C^\infty(M,L)$. Then from Proposition \ref{prop:smoothlocsect} we know that there is an $f_j$ with $s_{|U_j}=f_js_j$ for each $j$. Proposition \ref{prop:covderisloc} says that $(\nabla_\mathfrak Xs)(p)=(\nabla_\mathfrak Xs_{|U_j})(p)$ for each $p\in U_j$ and hence we have on a fixed $U_j$:
		\begin{equation} \label{eq:locform}
			\nabla_\mathfrak X s=\nabla_\mathfrak X f_js_j=\mathfrak X(f_j)s_j+2\pi if_j\alpha_j(\mathfrak X)s_j\quad\forall\mathfrak X\in\mathcal X_\mathbb C(M).
		\end{equation}
		Therefore $\nabla$ is uniquely determined by $\{U_j,s_j,\alpha_j\}_{j\in J}$.
		
		\begin{prop} \label{prop:locformglobform}
			Let $L$ be a complex line bundle over $M$, let $\{U_j,s_j\}_{j\in J}$ be a local system of $L$, and let $\{\alpha_j\}_{j\in J}$ be a family of local 1-forms satisfying \ref{eq:proplocform}. Then there is a uniquely determined connection $\nabla$ on $L$ having $\{U_j,s_j,\alpha_j\}_{j\in J}$ as a local form.
		\end{prop}
		
		\begin{proof}
			\cite[Corollary 1.4.1]{limaaa}.
		\end{proof}
		
		We consider a third way of looking at connections on $L$.
		
		We denote $\mathbb C^*:=\mathbb C\backslash\{0\}$. The 1-form $\frac{1}{2\pi i}\frac{dz}{z}$ on $\mathbb C^*$ is invariant under multiplication with elements of $\mathbb C^*$. For $p\in M$, any two identifications of $(L_p)\backslash\{0\}$ with $\mathbb C^*$ are $\mathbb C^*$-multiples of each other. Therefore, there is a unique 1-form $\beta_p$ on $(L_p)\backslash\{0\}$ such that for any $\mathbb C^*$-equivariant map
		$$\tau\colon\mathbb C^*\to L_p\backslash\{0\}$$
		one has $\tau^*(\beta_p)=\frac{1}{2\pi i}\frac{dz}{z}$.
		
		\begin{defi}
			Let $L$ be a complex line bundle over $M$. We define the open set\index{$L^*$, bundle without zero section}
			$$L^*:=\bigcup_{p\in M}(L_p\backslash\{0\}).$$
			Let $\alpha\in\Omega^1(L^*)$ be a 1-form. Then $\alpha$ is called a \textbf{connection form}\index{connection form} if it satisfies
			\renewcommand{\theenumi}{(\roman{enumi})}
			\renewcommand{\labelenumi}{\theenumi}
			\begin{enumerate}
				\item $\alpha$ is invariant under the action of $\mathbb C^*$ on $L$ by multiplication in the fibers; and
				\item for all $p\in M$ one has $\alpha_{|L_p\backslash\{0\}}=\beta_p$.
			\end{enumerate}
		\end{defi}
		
		\begin{prop} \label{prop:connform}
			Let $L$	be a complex line bundle over $M$. Then there is a one-to-one correspondence between connections in $L$ and connection forms on $L^*$. If $\nabla$ is a connection in $L$, $\alpha$ the corresponding connection form, and $\{U_j,s_j,\alpha_j\}$ a local form of $\nabla$ then one has for any $j$
			$$\alpha_j=s_j^*(\alpha).$$
		\end{prop}
		
		\begin{proof}
			\cite[Proposition 1.5.1]{limaaa}.
		\end{proof}
		
		\begin{defi}
			Let $(L,\nabla)$ be a line bundle with connection over $M$ and let $\{U_k,s_j,\alpha_j\}_{j\in J}$ be a local form. Then Proposition \ref{prop:proplocform} tells us that on any non-empty intersection $U_j\cap U_k$ we have $\alpha_k-\alpha_j=\frac{1}{2\pi i}\frac{dc_{jk}}{c_{jk}}$. Since
			\begin{align*}
				d\frac{dc_{jk}}{c_{jk}} & = d(\frac{1}{c_{jk}})\wedge dc_{jk} - \frac{1}{{c_{jk}}}d^2c_{jk} \\
																& = d^2(\log\circ c_{jk}) \\
																& = 0,
			\end{align*}
			we find $d\alpha_j=d\alpha_k$. Hence there exists a unique 2-form $\omega$ such that
			$$\omega_{|U_j}=d\alpha_j\quad\forall j\in J.$$
			It is called the \textbf{curvature}\index{curvature} of $\nabla$.
		\end{defi}
		
		\begin{prop}
			Up to a constant factor, our definition of the curvature coincides with the usual definition of a curvature of a connection. More precisely, we have
				$$2\pi i\omega(\mathfrak X,\mathfrak Y)s=(\nabla_\mathfrak X\nabla_\mathfrak Y-\nabla_\mathfrak Y\nabla_\mathfrak X-\nabla_{[\mathfrak X,\mathfrak Y]})s\quad\forall s\in C^\infty(M,L).$$
		\end{prop}
		
		\begin{proof}
			In a trivialization $(U_j,s_j)$ where $s=f_js_j$ we find using (\ref{eq:locform}):
			\begin{align*}
				&\quad\ (\nabla_\mathfrak X\nabla_\mathfrak Y-\nabla_\mathfrak Y\nabla_\mathfrak X-\nabla_{[\mathfrak X,\mathfrak Y]})s \\
				& = \left(\mathfrak X(\mathfrak Y(f_j)) + 2\pi i\mathfrak X(f_j)\alpha_j(\mathfrak Y)+2\pi if_j\mathfrak X(\alpha_j(\mathfrak Y)) + 2\pi i\mathfrak Y(f_j)\alpha_j(\mathfrak X)-4\pi^2 f_j\alpha_j(\mathfrak Y)\alpha_j(\mathfrak X) \right)s_j \\
				& - \left(\mathfrak Y(\mathfrak X(f_j)) + 2\pi i\mathfrak Y(f_j)\alpha_j(\mathfrak X)+2\pi if_j\mathfrak Y(\alpha_j(\mathfrak X)) + 2\pi i\mathfrak X(f_j)\alpha_j(\mathfrak Y)-4\pi^2 f_j\alpha_j(\mathfrak X)\alpha_j(\mathfrak Y) \right)s_j \\
				& - \left( \mathfrak X(\mathfrak Y(f_j))-\mathfrak Y(\mathfrak X(f_j))+2\pi i f_j\alpha_j([\mathfrak X,\mathfrak Y]) \right)s_j\\
				& = 2\pi i\left( (\mathfrak X(\alpha_j(\mathfrak Y))- \mathfrak Y(\alpha_j(\mathfrak X)) - \alpha_j([\mathfrak X,\mathfrak Y]))f_j s_j\right)\\
				& = 2\pi i (d\alpha_j(\mathfrak X,\mathfrak Y))s \\
				& = 2\pi i \omega(\mathfrak X,\mathfrak Y)s.
			\end{align*}
		\end{proof}
		
		\begin{prop} \label{prop:curvature}
			Let $\alpha$ be the connection form of $\nabla$. The curvature $\omega$ is a closed 2-form and we have the relation
			$$d\alpha=(\pi_{|L^*})^*\omega$$
		\end{prop}
		
		\begin{proof}
			\cite[Proposition 1.6.1]{limaaa}.
		\end{proof}
		
		\begin{prop} \label{prop:elbwchtsc}
			Let $(L_1,\nabla^1)$ and $(L_2,\nabla^2)$ be complex line bundles with connection over a smooth manifold $M$. If $(L_1,\nabla^1)$ and $(L_2,\nabla^2)$ are equivalent as line bundles with connection, then $\nabla^1$ and $\nabla^2$ have the same curvature.
		\end{prop}
		
		\begin{proof}
			Let $\tau$ be the equivalence diffeomorphism and let $\{U_j,s_j,\alpha_j\}_{j\in J}$ be a local form of $\nabla^1$. Then $\{U_j,\tau\circ s_j\}_{j\in J}$ is a local system of $L_2$ and we have
			\begin{align*}
				\nabla_{\mathfrak X}^2(\tau\circ s_j) & = \tau\circ(\nabla_\mathfrak X^1 s_j) \\
																							& = \tau\circ(\alpha_j(\mathfrak X)s_j) \\
																							& = \alpha_j(\mathfrak X)(\tau\circ s_j)\quad\forall \mathfrak X\in\mathcal X_\mathbb C(M).
			\end{align*}
			Therefore $\{U_j,\tau\circ s_j,\alpha_j\}_{j\in J}$ is a local form of $\nabla^2$ and the curvatures coincide.
		\end{proof}
		
		\begin{defi}
			A \textbf{Hermitian structure} $H$ in $L$\index{Hermitian structure} is a family $\{H_p\}_{p\in M}$, where $H_p$ is a Hermitian inner product on the fiber $\pi^{-1}(p)$. It is furthermore smooth, in the sense that for any smooth sections $s$ and $t\in C^\infty(M,L)$ the function $p\mapsto H_p(s(p),t(p))$ is smooth. A connection $\nabla$ on $L$ and a Hermitian structure $H$ are said to be \textbf{compatible}\index{connection!compatible}, if
			$$\mathfrak X( H(s,t) )= H(\nabla_\mathfrak X s,t)+H(s,\nabla_\mathfrak X t)\quad\forall\mathfrak X\in\mathcal X(M).$$
			In this case we call $\nabla$ a \textbf{Hermitian connection}\index{connection! Hermitian}. We denote by $\{|H|^2_p\}_{p\in M}$\index{$|H|^2_p$} the family of maps $|H|^2_p\colon L_p\to\mathbb R,\quad z\mapsto |H(z,z)|$.
		\end{defi}
		
		\begin{defi} \label{defi:elnabla}
			Let $(L_1,\nabla^1,H^1)$ and $(L_2,\nabla^2,H^2)$ be line bundles with Hermitian connection over smooth manifolds $M$ and $N$ respectively. A \textbf{isomorphism of line bundles}\index{line bundle! isomorphism of line bundles} is a diffeomorphism $\tau\colon L_1\to L_2$ such that there exists a (uniquely determined) diffeomorphism $\check{\tau}\colon M\to N$ making the following diagram commutative
			$$\begin{xy}
     				\xymatrix{
         			L_1 \ar[r]^-\tau \ar[d]^-{\pi_1}    &   L_2 \ar[d]^-{\pi_2} \\
         			M \ar[r]_-{\check{\tau}} & N
         			}
   		\end{xy}$$
			and such that for all $p\in M$ the map $\tau_{|(L_1)_p}\colon (L_1)_p\to (L_2)_{\check{\tau}(p)}$ is a linear isomorphism. We call $\tau$ a \textbf{isomorphism of line bundles with connection}, if $\tau$ satisfies in addition
			\begin{equation} \label{eq:elnablaprop}
				\tau\circ(\nabla^1_\mathfrak X s)\circ \check{\tau}^{-1}=\nabla^2_{\check{\tau}_*\mathfrak X}(\tau\circ s\circ \check{\tau}^{-1})\quad\forall s\in C^\infty(N,L_2)\quad\forall \mathfrak X\in\mathcal X(M).
			\end{equation}
			$\tau$ is called \textbf{isometric}\index{line bundle!isometric map}, if $|H^2|^2\circ\tau=|H^1|^2$. We denote by $E(L,\nabla)$\index{$E(L,\nabla)$, isometric bundle maps} the subgroup of $D(L)$ of diffeomorphisms which are isometric isomorphisms of line bundles with Hermitian connection. Since Kostant shows in \cite[Proposition 1.9.1]{limaaa} that compatible Hermitian structures are unique up to factor, we see that the definition of $E(L,\nabla)$ is independent of the chosen $H$.
		\end{defi}
		
	\subsection{Homogeneous Line Bundles}
		
		\begin{defi}
			Let $G$ be a Lie group and $M$ a homogeneous $G$-space with smooth action $\sigma\colon G\times M\to M$. A complex line bundle $L$ over $M$ is called a \textbf{homogeneous line bundle}\index{line bundle! homogeneous} if there is a smooth transitive action $\sigma_L\colon G\times L\to L$\index{$\sigma_L$, action on line bundle} such that
			\renewcommand{\theenumi}{(\roman{enumi})}
			\renewcommand{\labelenumi}{\theenumi}
			\begin{enumerate}
				\item $\sigma_L(g)L_p=L_{\sigma(g)p}\quad\forall p\in M\ \forall g\in G$
				\item The mapping $\sigma_L(g)_{|L_p}\colon L_p\to L_{\sigma(g)p}$ is a complex-linear isomorphism for all $p\in M$ and $g\in G$.
			\end{enumerate}
			If $(L,\nabla,H)$ is a line bundle with Hermitian connection, then we say that it is a \textbf{homogeneous line bundle with Hermitian connection}, if $\sigma_L(G)\subseteq E(L,\nabla)$.
		\end{defi}
		
		\begin{defi} \label{defi:actonsect}
			Let $L$ be a homogeneous line bundle over $M$. Then we define an action $\sigma_S\colon G\times C^\infty(M,L)\to C^\infty(M,L)$\index{$\sigma_S$, action on sections} on the space of smooth sections by the formula
			$$(\sigma_S(g)s)(p)=\sigma_L(g)(s(\sigma(g^{-1})p))\quad\forall p\in M.$$
			This means we first manipulate the argument of $s$ and then move the value back into the right fiber. We observe that $G$ acts by linear transformations and, in fact, we obtain an infinite-dimensional representation of $G$ in $C^\infty(M,L)$ in this way.
		\end{defi}
		
		We give a basic construction which describes all homogeneous line bundles over a homogeneous $G$-space $M$. First of all, we assume $M$ to be of the form $M=G/H$ which is no restriction as we know from Theorem \ref{theo:homoisquot}.
		
		\begin{cons} \label{cons:homolb}
			Let $\chi\colon H\to\mathbb C^*$ be a character of $H$. Then $H$ acts on $G\times\mathbb C$ from the right via
			$$(g,z)\cdot h=(gh,\chi(h^{-1})z)\quad\forall g\in G\quad \forall z\in\mathbb C$$
			The space of orbits under this action is denoted by $G\times_\chi\mathbb C$\index{$G\times_\chi\mathbb C$, glued line bundle}. We define the map $\pi\colon G\times_\chi\mathbb C\to G/H$ by $\pi([g,z]):=gH$. This is well-defined since any representative of $[g,z]$ is mapped to the coset $gH$. If we give $G\times_\chi\mathbb C$ the quotient topology, then $\pi$ is continuous. 
		\end{cons}
		
		\begin{prop}
			The pair $(G\times_\chi\mathbb C,\pi)$ together with the action $g_1.[g_2,z]:=[g_1g_2,z]$ is a homogeneous line bundle over $G/H$.
		\end{prop}
		
		\begin{proof}
			\cite[5.2.2]{walhomo}.
		\end{proof}
		
		\begin{prop} \label{prop:sectfunc}
			Let $G$ be a Lie group, $H$ a closed subgroup with a character $\chi\colon H\to\mathbb C^*$ and $G\times_\chi\mathbb C$ the homogeneous line bundle from Construction \ref{cons:homolb}. Consider the following subspace of $C^\infty(G,\mathbb C)$.
			$$C^\infty_\chi(G,\mathbb C):=\{f\in C^\infty(G,\mathbb C)\ |\ f(gh)=\chi(h^{-1})f(g)\quad\forall h\in H\}$$\index{$C^\infty_\chi(G,\mathbb C)$, sections as functions}
			$G$ acts on $C^\infty_\chi(G,\mathbb C)$ by left-translation in the argument. With respect to this action and the action $\sigma_S$ on $C^\infty(G/H,G\times_\chi\mathbb C)$ from Definition \ref{defi:actonsect} the map
			$$\Phi_\chi\colon C^\infty_\chi(G,\mathbb C)\to C^\infty(G/H,G\times_\chi\mathbb C),\quad f\mapsto\left(gH\mapsto[g,f(g)]\right)$$
			is $G$-equivariant bijection.
		\end{prop}
		
		\begin{proof}
			\cite[Proposition 2.12]{hang}.
		\end{proof}
		
		We have found in Theorem \ref{theo:homoisquot} a uniform description for homogeneous $G$-spaces. We now do the same for homogeneous line bundles.
		
		\begin{theo} \label{theo:hlbisglued}
			Let $L$ be a homogeneous line bundle over $G/G_p$. We define the map $\chi\colon G_p\to\mathbb C$ by the formula
			$$\chi(g)x=\sigma_L(g)x\quad\forall g\in G.$$
			Then $L$ is equivalent to $G\times_\chi\mathbb C$.
		\end{theo}
		
		\begin{proof}
			\cite[Lemma 5.2.3]{walhomo}.
		\end{proof}
		
	\subsection{Holomorphic Line Bundles}
		
		\begin{defi} 
			Let $M$ be a complex manifold. A \textbf{holomorphic line bundle}\index{line bundle!holomorphic} over $M$ is a complex line bundle over $M$ for which there exists a local system $\{U_j,s_j\}_{j\in J}$ such that the transition functions are holomorphic. We call this a \textbf{holomorphic local system}.
		\end{defi}
		
		\begin{prop}
			If $L$ is a holomorphic line bundle over $M$ and $\{U_j,s_j\}_{j\in J}$ is a holomorphic local system, then there is a unique complex structure (cp. Definition \ref{defi:cs}) on $L$ such that for the associated family of local trivializations the maps
			$$\varphi\colon\pi^{-1}(U)\to U\times\mathbb C,\quad\varphi(q)=(\pi(p),\phi(p))$$
			and their inverses are holomorphic. Furthermore, the sections $s_j$ are holomorphic mappings between complex manifolds and any section $C^\infty(M,L)$ is holomorphic if and only if it has the local form $s=f s_j$ where $f$ is a holomorphic function. The vector space of holomorphic sections is denoted by $\Gamma_{hol}(M,L)$\index{$\Gamma_{hol}(M,L)$, holomorphic sections}.
		\end{prop}
		
		\begin{proof}[Sketch of proof]
			Let $\{U_j,\varphi_j\}_{j\in J}$ be a family of local trivializations belonging to the holomorphic local system $\{U_j,s_j\}_{j\in J}$. Then the inverse $\varphi^{-1}:U_j\times\mathbb C\to\pi^{-1}(U_j)$ is given by $(p,z)\mapsto zs_j(p)$. If $(U_j,\psi_j)$ is a holomorphic chart for $M$, then we obtain a chart $(\tilde{U_j},\tilde{\varphi_j})$ for $L$ by setting $\tilde{U_j}:=\pi^{-1}(U_j)$ and $\tilde{\varphi_j}(v):=(\psi_j(p),z)\in\mathbb C^{n+1}$, where $v=zs_j(p)$. If $U_j,\ U_k$ have non-empty intersection, then the change of charts is given by $(z_1,\dots,z_n,z)\mapsto (\psi_k\circ\psi_j^{-1}(z_1,\dots,z_n),c_{jk}\circ\psi_j^{-1}z)$ and hence holomorphic. One can show that the resulting holomorphic atlas has the desired properties. Confer \cite[p.132]{gh}.
		\end{proof}
		
		
		\begin{prop} \label{prop:voganhlb}
			Let $G$ be a Lie group with closed Lie subgroup $H$ and let $\chi\colon H\to\mathbb C^*$ be a character of $H$. Suppose $G/H$ carries an invariant complex structure given by $\mathfrak b$ (cp. Proposition \ref{prop:csalg}) and let $G\times_\chi\mathbb C$ be the homogeneous line bundle over $G/H$ from Construction \ref{cons:homolb}. Then to make $G\times_\chi\mathbb C$ a holomorphic vector bundle amounts to giving a Lie algebra representation $\mu$ of $\overline{\mathfrak b}$ on $\mathbb C$, satisfying
			\begin{enumerate}
				\item the complex linear extension of the differential of the group representation $\chi:H\to\mathbb C^*$ agrees with the Lie algebra representation $\mu$ restricted to $\mathfrak h_\mathbb C$; and
				\item for $h$ in $H$, $X\in\overline{\mathfrak b}$, and $w\in\mathbb C$
					$$\chi(h)(\mu(X)w)=\mu(\Ad(h)X)(\chi(h)w).$$
			\end{enumerate}
			The space $\Gamma(G/H,G\times_\chi\mathbb C)$ of holomorphic sections may be identified with the subspace $C^\infty_{\chi,\mu}(G,\mathbb C)$ of $C^\infty_{\chi}(G,\mathbb C)$ (cp. Proposition \ref{prop:sectfunc}) given by
			$$C^\infty_{\chi,\mu}(G,\mathbb C):=\{f\in C^\infty_\chi(G,\mathbb C)\ |\ \tilde{X}(f)(g)=\mu(X)(f(g))\quad\forall X\in\overline{\mathfrak b}\}.$$\index{$C^\infty_{\chi,\mu}(G,\mathbb C)$, holomorphic sections}
		\end{prop}
		
		\begin{proof}
			A proof of this theorem is given in \cite[Theorem 3.6]{tirwolf}.
		\end{proof}
		
		\begin{prop} \label{prop:hlb}
			Let $G$ be a compact connected Lie group, $T_1$ a torus and $\chi\colon C(T_1)\to \mathbb C$ a character of $C(T_1)$. We give $G/C(T_1)$ a complex structure by applying Theorem \ref{theo:eopocot}. Then the homogeneous line bundle $G\times_\chi\mathbb C$ over $G/C(T_1)$ is a holomorphic line bundle.
		\end{prop}
		
		\begin{proof}[Sketch of proof]
			The complex structure on $G/C(T_1)$ comes from a subalgebra $\mathfrak b$ after we have chosen a positive Weyl chamber. Then $\overline{\mathfrak b}$ is of the following form.
			$$\overline{\mathfrak b}=\mathfrak c_{1\mathbb C}+\sum_{\alpha\in(\Delta^+\backslash\Delta^+_1)}\mathfrak g^{-\alpha}$$
			We define $\mu$ to be the induced Lie algebra representation $L\chi$ of $\mathfrak c_{1\mathbb C}$ extended to $\overline{\mathfrak b}$ by letting the right summand act trivially. It remains to show that $\ref{prop:voganhlb}(ii)$ holds. Then the claim follows with this proposition. Confer \cite[Proposition 6.3.3]{walhomo}.
		\end{proof}
		
		\begin{prop} \label{prop:repinholo}
			Let $\sigma_S$ be the action of $G$ on the space of smooth sections from Definition \ref{defi:actonsect}. Then the restriction of $\sigma_S$ to the space of holomorphic sections $\Gamma(G/C(T_1),G\times_{\chi}\mathbb C)$ is a representation of $G$.
		\end{prop}
		
		\begin{proof}
			We already know that $G$ acts linearly on the smooth sections. It remains to show that the space of holomorphic sections is invariant under the action of $G$. We do that by showing that $C^\infty_{\chi,\mu}(G,\mathbb C)$ is invariant.
			
			Let $f\in C^\infty_{\chi,\mu}(G,\mathbb C)$ and $h\in G$. Then $h\cdot f\in C_\chi^\infty(G,\mathbb C)$ and
			\begin{align*}
				(\tilde{X} (h\cdot f))(g) & = \tilde{X}(g)(f\circ\lambda_{h^{-1}}) \\
																	& = (\lambda_{h^{-1}}^\prime(g)\tilde{X}(g))(f) \\
																	& = \tilde{X}(h^{-1}g)(f) \\
																	& = \mu(X)(f(h^{-1}g)) \\
																	& = \mu(X)((h\cdot f)(g))\quad\forall\mathfrak X\in\mathfrak b^-.
			\end{align*}
			Therefore $h\cdot f\in C^\infty_{\chi,\mu}(G,\mathbb C)$.
		\end{proof}
		
	\section{Borel-Weil Theorem}
		
		\begin{theo}[Cartan-Weyl]\index{Cartan-Weyl Theorem} \label{theo:cwt}
			Let $G$ be a compact connected Lie group, $T$ a maximal torus and let $\Delta^+(\mathfrak g_\mathbb C,\mathfrak t_\mathbb C)$ be a system of positive roots in $\Delta(\mathfrak g_\mathbb C,\mathfrak t_\mathbb C)$. Then apart from equivalence the irreducible finite-dimensional representations $(\pi,V)$ of $G$ stand in one-one correspondence with the dominant analytically integral linear functionals $\lambda\in\mathfrak t$, the correspondence being that $\lambda$ is the highest weight of $\pi$. More precisely, if we put
			$$V^{\mathfrak n}:=\{v\in V\ |\ L\pi(X)v=0\quad\forall X\in\mathfrak n\},$$
			then $V^{\mathfrak n}$ is invariant under $T$, and the corresponding representation $\pi^+$ of T on $V^{\mathfrak n}$ is irreducible and has dominant weight $\lambda$.
		\end{theo}
		
		\begin{proof}
			\cite[Theorem 5.110]{knapp}.
		\end{proof}
		
		Example \ref{exam:char} and Corollary \ref{coro:irreunitabel} tell us that an irreducible unitary representation of $T$ can be viewed as a character $\tau\colon T\to S^1$. Therefore the theorem suggests that the equivalence classes of unitary irreducible representations are parametrized by characters of a maximal torus.
		
		\begin{cons} \label{cons:borel-weil}
			Let $G$ be a compact connected Lie group, $T$ a maximal torus in $G$ and $\tau\colon T\to S^1$ a character. Let $\lambda:=L\tau\in\mathfrak t^*$. Note that according to Proposition \ref{prop:uniqueforconnected} the character $\tau$ is uniquely determined by $\lambda$. We apply Construction \ref{cons:homolb} to obtain a homogeneous line bundle $G\times_\tau\mathbb C$ over $G/T$.	We choose a positive Weyl chamber to obtain a suitable subalgebra $\mathfrak b$ as in Theorem \ref{theo:eopocot}. We apply Proposition \ref{prop:csalg} to give $G/T$ a complex structure. Then we can use Proposition \ref{prop:hlb} which gives $G\times_\tau\mathbb C$ the structure of a holomorphic line bundle. Using Proposition \ref{prop:repinholo} we have a representation $(\pi_\lambda,\Gamma(G/T,G\times_\tau\mathbb C))$ of $G$ in the space of holomorphic sections.
		\end{cons}
		
		The following is a classical result in representation theory of compact connected Lie groups.
		
		\begin{theo}[Borel-Weil]\index{Borel-Weil Theorem}
			Let $G$ be a compact connected Lie group, $T$ a maximal torus and $\tau\colon T\to S^1$ a character. Apply Construction \ref{cons:borel-weil} to obtain a representation $(\pi_\lambda,\Gamma(G/T,G\times_\tau\mathbb C))$ of $G$. Then $\Gamma(G/T,G\times_\tau\mathbb C)$ is non-zero if and only if $L\tau$ is dominant. In this case, $\pi_\lambda$ is the irreducible representation attached to $L\tau$ by the Cartan-Weyl Theorem.
		\end{theo}
		
		\begin{proof}
			\cite[Theorem 4.12.5]{dk}.
		\end{proof}
		
		\begin{cons} \label{cons:gborel-weil}
			Let $G$ be a compact connected Lie group, $T$ a maximal torus in $G$, $T_1\subseteq T$ a torus and $\chi\colon C(T_1)\to S^1$ a character. We apply Construction \ref{cons:homolb} to obtain a homogeneous line bundle $G\times_\chi\mathbb C$ over $G/C(T_1)$.  We choose a $T_1$-admissible Weyl chamber to obtain a suitable Lie algebra $\mathfrak b$ as in Theorem \ref{theo:eopocot}. Then Proposition \ref{prop:csalg} tells us that there exists a complex structure on $G/C(T_1)$. We can use Proposition \ref{prop:hlb} to give $G\times_\chi\mathbb C$ the structure of a holomorphic line bundle. Using Proposition \ref{prop:repinholo} we obtain a representation $(\pi_\chi,\Gamma(G/C(T_1),G\times_\chi\mathbb C))$ of $G$ in the space of holomorphic sections.
		\end{cons}
		
		The following theorem is taken from \cite[Theorem 6.3.7]{walhomo}. We remark that it is actually true for finite-dimensional representations $(\sigma,W)$ of a centralizer of a torus $T_1$ instead of just one-dimensional representations $(\sigma,\mathbb C)$.
		
		\begin{theo}\label{theo:gbwt}
			Let $G$ be a compact connected Lie group and let $T_1$ be a torus in $G$. Let $(\sigma,\mathbb C)$ be an irreducible unitary representation of $C(T_1)$. If there exists no irreducible representation $(\pi,V)$ of $G$ so that the representation $(\pi_{|C(T_1)},V^{\mathfrak n})$ is equivalent with $(\sigma,\mathbb C)$, then $\Gamma(G/C(T_1),G\times_\sigma\mathbb C)=\{0\}$. If there exists such a $(\pi,V)$, then $(\pi_\sigma,\Gamma(G/C(T_1),G\times_\sigma\mathbb C))$ and $(\pi,V)$ are equivalent.
		\end{theo}
		
		\begin{prop} \label{prop:11char}
			Let $G$ be a compact connected Lie group, $T$ a maximal torus and $\lambda\in\mathfrak t^*$. Then there is a one-to-one correspondence between characters $\tau\colon T\to S^1$ and $\chi\colon G_\lambda\to S^1$ which satisfy $L\tau=L\chi_\lambda=:2\pi i\lambda$. As before (cp. Subsection \ref{structtheo}) we interpret $\lambda$ as a functional on $\mathfrak g_\lambda$ which is zero on the orthogonal complement of $\mathfrak t$.
		\end{prop}
		
		\begin{proof}
			The following proof involves some concepts which were not introduced in this work and which are of no further relevance later on. For more detailed definitions confer \cite{knapp}.
			
			We have investigated the relation between $\lambda$, $T$ and $G_\lambda$ in Subsection \ref{ss:stabofunc}. There we have shown that $T\subseteq G_\lambda$ and that the Lie algebra of $G_\lambda$ is given by
			$$\mathfrak g_\lambda=\mathfrak t + \sum_{\alpha\in\Delta_{sing}^+(\lambda)}\mathfrak g\cap(\mathfrak g^\alpha+\mathfrak g^{-\alpha}).$$
			Therefore if $\chi:G_\lambda\to S^1$ is a character of $G_\lambda$, then $\chi_{|T}$ is a character of $T$ satisfying $L(\chi_{|T})=L\chi_{|\mathfrak t}$. Hence the first direction is realized by restricting $\chi$ to $T$.
			
			Now let $\tau:T\to S^1$ with $L\tau=2\pi i\lambda$ be given. Then we extend $2\pi i \lambda$ to a functional $\tilde{\lambda}$ on $\mathfrak g_\lambda$ which is zero on the root spaces. Let $\widetilde{G_\lambda}$ be a universal covering group of $G_\lambda$ (see \cite[§I.11]{knapp} for definition and properties). Then $\widetilde{G_\lambda}$ is a simply connected Lie group with Lie algebra $\mathfrak g_\lambda$ and there is a Lie group homomorphism $\pi:\widetilde{G_\lambda}\to G_\lambda$ which is also a local homeomorphism. According to Proposition \ref{prop:liftexistforsc} there exists a unique character $\tilde{\chi}$ of $\widetilde{G_\lambda}$ such that the following left diagram is commutative. Our goal is to find a character $\chi$ of $G_\lambda$ such that the right diagram is commutative. Because then $L\chi=2\pi i\tilde{\lambda}$ and the claim follows with the argumentation from the first part of the proof.
			
			$$\begin{xy}
				\xymatrix{
					\widetilde{G_\lambda} \ar[r]^-{\tilde{\chi}} & S^1 & \widetilde{G_\lambda} \ar[r]^-{\tilde{\chi}} \ar[d]_-{\pi}& S^1 \\
					\mathfrak g_\lambda \ar[u]^-{\exp_{\widetilde{G_\lambda}}} \ar[r]^-{2\pi i \tilde{\lambda}} & \mathbb R \ar[u]_-{e^{i.}} & G_\lambda \ar@{.>}[ru]_-{\chi}
				}			
			\end{xy}$$
			
			To prove the existence of $\chi$, we will show $\tilde{\chi}_{|\ker\pi}=1$. In this case we can define $\chi(a)$ via the preimages of $a$ under the covering map. More precisely, if we interpret $G_\lambda\simeq \widetilde{G_\lambda}/\ker\pi$ (cf. \cite[Proposition 1.101]{knapp}), then $\chi$ is well-defined by
			$$\chi(g\ker\pi):=\tilde{\chi}(g),\quad\forall g\in\widetilde{G_\lambda}.$$
			First of all we show $\ker\pi\subseteq Z(\widetilde{G_\lambda})$. As a kernel of a Lie group homomorphism $\ker\pi$ is normal and since $\pi$ is a covering, it is furthermore discrete. Let $n\in\ker\pi$. Then the function $f_n:g\mapsto gng^{-1}$ on $\widetilde{G_\lambda}$ is continuous and has values in $\ker\pi$ since the latter is normal. But since $\ker\pi$ is a discrete set and $\widetilde{G_\lambda}$ is connected the function $f_n$ is constant with value $f_n(e)=n$. But this means in fact that $\ker\pi$ is contained in the center of $\widetilde{G_\lambda}$.
			
			Since $G_\lambda$ is compact, its Lie algebra decomposes into $\mathfrak g_\lambda=Z(\mathfrak g_\lambda)\oplus\mathfrak k$, where $\mathfrak k$ denotes the semi-simple Lie algebra $[\mathfrak g_\lambda,\mathfrak g_\lambda]$ (cf. \cite[Corollary 4.25]{knapp}). On group level, $G_\lambda$ decomposes into the commuting product $G_\lambda=Z(G_\lambda)K$, where $K$ is a semi-simple Lie group with Lie algebra $\mathfrak k$ (cf. \cite[Theorem 4.29]{knapp}). Then we can assume that $\widetilde{G_\lambda}$ is of the form $\widetilde{G_\lambda}=\widetilde{K}\times Z(\mathfrak g_\lambda)$. Here $\widetilde{K}$ is a universal covering group of $K$. Then Weyl's Theorem (\cite[Theorem 4.69]{knapp}) tells us that $\widetilde{K}$ is compact.
			
			We note that $T$ is a maximal torus in $G_\lambda$. Hence we obtain a decomposition $\mathfrak t=Z(\mathfrak g_\lambda)\oplus\mathfrak t_K$ with the same arguments as above. Hence we can choose a maximal torus $T_{\tilde{K}}$ in $\widetilde{K}$ which has Lie algebra $\mathfrak t_K$. We then consider the subgroup $A:=T_{\tilde{K}}\times Z(\mathfrak g_\lambda)$ of $\widetilde{G_\lambda}$. Via the decomposition of the Lie algebra $\mathfrak t$ we can assume that $T=Z(G_\lambda)T_K$, where $T_K$ is a subgroup of $K$ with Lie algebra $\mathfrak t_K$. Hence $\pi_A:A\to T$ is a covering. In particular, $A$ has Lie algebra $\mathfrak t$.
			
			The center $Z(\widetilde{G_\lambda})$ is of the form $Z(\widetilde{G_\lambda})=Z(\widetilde{K})\times Z(\mathfrak g_\lambda)$. Since the center $Z(\widetilde{K})$ of a compact Lie group lies in every maximal torus (cf. \cite[Corollary 4.47]{knapp}), we find $Z(\widetilde{G_\lambda})\subseteq A$. With the inclusion shown above, we finally arrive at $\ker\pi\subseteq A$.
			
			We consider the following two characters of $A$. Firstly, we have the restriction of $\tilde{\chi}$ to $A$. Its derivative is given by $L(\tilde{\chi}_{|A})={(L\tilde{\chi})}_{|\mathfrak t}=\lambda$. Secondly, we consider $\tau\circ\pi_{|A}$. Since $\pi$ is a covering map, this character has also differential equal to $\lambda$. Therefore $\tilde{\chi}_{|A}=\tau\circ\pi_{|A}$ according to Proposition \ref{prop:uniqueforconnected}.
			
			Since $\tau$ maps $e$ to 1 and because $\ker\pi\subseteq A$, we find
			\begin{align*}
				1 = \tau\circ\pi_{|\ker\pi} & = (\tau\circ\pi_{|A})_{|\ker\pi} \\
																		& = (\tilde{\chi}_{|A})_{|\ker\pi} \\
																		& = \tilde{\chi}_{|\ker\pi}.
			\end{align*}
			With the argumentation above, this proves the existence of $\chi$.
		\end{proof}
		
		\begin{rema}
			Let us make two remarks at this point.
			
			The above proposition is only non-trivial for the case that the functional $\lambda\in\mathfrak t^*$ is singular. If $\lambda$ is regular, then $T$ and $G_\lambda$ coincide and the statement is a tautology (cp. Subsection \ref{ss:stabofunc}).
			
			If $\tau:T\to S^1$ is a character and $\lambda:=\frac{1}{2\pi i}L\tau$, then $\tau$ is uniquely determined by $\lambda$ according to Proposition \ref{prop:uniqueforconnected}. One often interprets the existence of $\tau$ as a special property of $\lambda$. If a functional $\lambda$ lifts to a character $\tau$ of $T$, then $\lambda$ is called \textbf{analytically integral}\index{functional!analytically integral} (cp. \cite[§IV.7]{knapp}).
		\end{rema}
		
		\begin{prop} \label{prop:borelweildisguised}
			Let $G$ be a compact connected Lie group, $T$ a maximal torus in $G$ and $\tau\colon T\to S^1$ a character. Let $\lambda:=\frac{1}{2\pi i}L\tau\in\mathfrak t^*$. Apply Proposition \ref{prop:11char} to obtain a character $\chi$ of $G_\lambda=C(T_1)$. We apply Construction \ref{cons:borel-weil} to $\tau$ and choose a $T_1$-admissible Weyl chamber $P_1$. We apply Construction \ref{cons:gborel-weil} to $\chi$ and choose a corresponding Weyl chamber. Then $(\pi_\lambda,\Gamma(G/T,G\times_\tau\mathbb C))$ and $(\pi_\chi,\Gamma(G/C(T_1),G\times_\chi\mathbb C))$ are equivalent.
		\end{prop}
		
		\begin{proof}
			We apply Theorems \ref{theo:cwt} and \ref{theo:gbwt}. The unitary representation $(\tau,\mathbb C)$ of $C(T)=T$ is one-dimensional and hence irreducible. If $\lambda$ is a dominant weight, then the Cartan-Weyl Theorem tells us that there is a representation $(\pi,V)$ with highest weight $\lambda$ such that $(\pi_{|T},V^{\mathfrak n})$ is equivalent to $(\tau,\mathbb C)$. Then the same is true for $(\chi,\mathbb C)$ since it has also highest weight $L\chi_{|T}=\lambda$ (cp. \cite[Theorem 5.110]{knapp}). Then Theorem \ref{theo:gbwt} tells us that $\tau$ and $\chi$ induce equivalent representations of $G$ in the respective spaces of holomorphic sections with highest weight $\lambda$. On the other hand, if $\lambda$ is not dominant, then there exists no such representation and the spaces of holomorphic sections are trivial.
		\end{proof}
		
	\chapter{The Orbit Method}
	
	\section{Geometry of Coadjoint Orbits}
	
	\subsection{Symplectic Structure}
		
		\begin{defi}
			Let $V$ be a finite-dimensional vector space and let $\omega$ be a skew-symmetric bilinear non-degenerate form on $V$. Then $V$ is called a \textbf{symplectic vector space}\index{symplectic! vector space} with \textbf{symplectic form}\index{symplectic!form on a vector space} $\omega$. For a symplectic form $\omega$ on $V$ and a set $W\subseteq V$ one defines the orthogonal complement in the usual way as
			$$W^{\bot_\omega}:=\{v\in V\ |\ \omega(v,W)=0\}.$$
			A subspace $P\subseteq V$ is called \textbf{isotropic}\index{isotropic subspace} if $\omega(P,P)=0$. It is called \textbf{Lagrangian}\index{Lagrangian subspace}, if $P^{\bot_\omega}=P$.
		\end{defi}
		
		\begin{defi}
			Let $M$ be a smooth manifold, and let $\omega$ be a differential 2-form on $M$. Then $(M,\omega)$ is called a \textbf{symplectic manifold}\index{symplectic!manifold}\index{manifold!symplectic} if
			\renewcommand{\theenumi}{(\roman{enumi})}
			\renewcommand{\labelenumi}{\theenumi}
			\begin{enumerate}
				\item $\omega_p$ is a symplectic form on $T_pM$ for all $p\in M$; and
				\item $\omega$ is closed.
			\end{enumerate}
		\end{defi}
		
		\begin{prop} \label{prop:maxlagrangedim}
			Let $V$ be a finite-dimensional symplectic vector space, then $V$ has even dimension $2m$. Lagrangian subspaces exist and are isotropic subspaces of dimension $m$.
		\end{prop}
		
		\begin{proof}
			\cite[Lemma 1.2.4]{woo}.
		\end{proof}
		
		\begin{coro}
			Symplectic manifolds have even dimension.
		\end{coro}
		
		Let $(M,\omega)$ be a symplectic manifold. Because $\omega_p$ is a non-degenerate form on each tangent space $T_p$, we obtain a canonical identification of $T_pM$ with $T_pM^*$ (cf. \cite[§1.1]{woo}). Therefore, the contraction of $\omega$ with a non-zero vector field uniquely determines a non-zero 1-form. In this way we obtain an isomorphism
		$$\mathcal X(M)\to\Omega^1(M),\quad\mathfrak X\mapsto i_\mathfrak X\omega.$$
		
		\begin{defi}
			Let $f\in C^\infty(M,\mathbb R)$. The \textbf{Hamiltonian vector field}\index{vector field!Hamiltonian} $\mathfrak X_f$\index{$\mathfrak X_f$, Hamiltonian vector field} associated to $f$ is the vector field on $M$ uniquely determined by the relation
			$$i_{\mathfrak X_f}\omega=df.$$
			The set of all Hamiltonian vector fields is denoted by $\Ham(M,\omega)$\index{$\Ham(M,\omega)$, Hamiltonian vector fields}. The bilinear form $\{.,.\}$\index{$\{.,.\}$, Poisson bracket} on $C^\infty(M,\mathbb R)$ defined by the formula
			$$\{f,g\}:=\omega(\mathfrak X_f,\mathfrak X_g)\quad\forall f,\ g\in C^\infty(M,\mathbb R)$$
			is called the \textbf{Poisson bracket}\index{Poisson bracket}.
		\end{defi}
		
		\begin{prop}
			$(\Ham(M,\omega),[.,.])$ and $(C^\infty(M,\mathbb R),\{.,.\})$ are Lie algebras and $f\mapsto\mathfrak X_f$ is a Lie algebra homomorphism.
		\end{prop}
		
		\begin{proof}
			\cite[Proposition 7.3.4]{hga}.
		\end{proof}
		
		\begin{defi}
			Let $(M,\omega)$ be a symplectic manifold. Then a smooth map $\tau\colon M\to M$ is called \textbf{symplectomorphism}\index{symplectomorphism} or \textbf{canonical transformation}\index{canonical transformation}, if for all $p\in M$
			$$\omega(\tau^\prime(p)v,\tau^\prime(p)w)_{\tau(p)}=\omega(v,w)_p\quad\forall v,\ w\in T_pM.$$
			A vector field on $M$ is called \textbf{symplectic}\index{symplectic!vector field}\index{vector field!symplectic} if the local flows are symplectomorphisms. The set of symplectic vector fields is denoted by $\Symp(M,\omega)$\index{$\Symp(M,\omega)$, symplectic vector fields}.
		\end{defi}
		
		\begin{prop} \label{prop:sympvf}
			\mbox{}
			\renewcommand{\theenumi}{(\roman{enumi})}
			\renewcommand{\labelenumi}{\theenumi}
			\begin{enumerate}
				\item $\mathfrak X\in\Symp(M,\omega)\quad\Leftrightarrow\quad\mathcal L_\mathfrak X\omega=0\quad\Leftrightarrow\quad d(i_\mathfrak X\omega)=0$
				\item $\Ham(M,\omega)\subset\Symp(M,\omega)$
				\item $[\Symp(M,\omega),\Symp(M,\omega)]\subseteq\Ham(M,\omega)$
			\end{enumerate}
		\end{prop}
		
		\begin{proof}
			\cite[Proposition 7.3.2, Proposition 7.3.3]{hga}.
		\end{proof}
		
		\begin{defi}
			Let $G$ be a Lie group, $(M,\omega)$ a symplectic manifold, and $\sigma\colon G\times M\to M$ a smooth action. Then $\sigma$ is called a \textbf{symplectic action}\index{symplectic!action}\index{action!symplectic} if for each $g\in G$ the diffeomorphism $\sigma(g)$ is a symplectomorphism (cp. Definition \ref{defi:action}). It is called \textbf{strongly symplectic}\index{action!strongly symplectic} or \textbf{almost Hamiltonian}\index{action!almost Hamiltonian}, if for all $X\in\mathfrak g$ the form $i_{L{\sigma}(X)}\omega$ is exact. This means for all $X\in\mathfrak g$ we have $L{\sigma}(X)\in\Ham(M,\omega)$. It is called \textbf{Hamiltonian}\index{action!Hamiltonian}, if there is a homomorphism $\mu$ of Lie algebras $\mathfrak g$ and $C^\infty(M,\mathbb R)$ such that $\mathfrak X_{\mu(X)}=L{\sigma}(X)$ for all $X\in\mathfrak g$. A smooth manifold together with a smooth transitive Hamiltonian action is called a \textbf{Hamiltonian $G$-space}\index{Hamiltonian $G$-space}.
		\end{defi}
		
		\begin{prop}
			If $G$ is a connected Lie group and $(M,\omega)$ a symplectic manifold, then a smooth action of $G$ on $M$ is symplectic if and only if $L{\sigma}(X)\in\Symp(M,\omega)$ for all $X\in\mathfrak g$.
		\end{prop}
		
		\begin{proof}
			\cite[Proposition 8.2.1]{hga}.
		\end{proof}
		
		From Proposition \ref{prop:sympvf} we know that the property $L{\sigma}(X)\in\Symp(M,\omega)$ is equivalent to $d(i_{L{\sigma}(X)}\omega)=0$. Hence the definition of a strongly symplectic action is really a tightening of the definition of symplectic action.
		
	\subsection{Kähler Manifolds}
		
		
		\begin{defi}
			Let $M$ be a smooth manifold. A \textbf{Riemannian metric}\index{Riemannian metric} $\mathcal G$ on $M$ is a smooth differentiable real-valued pairing $\mathcal G$ such that $\mathcal G_p$ is an inner product on each tangent space $T_pM$. Note that a Riemannian metric allows for concepts such as measuring of arcs lengths by integration over the lengths of tangent vectors. If $(M,\mathcal J)$ is an almost complex manifold and $\mathcal G$ satisfies $\mathcal G(\mathcal J\mathfrak X,\mathcal J\mathfrak Y)=\mathcal G(\mathfrak X,\mathfrak Y)$ for all vector fields $\mathfrak X,\ \mathfrak Y\in\mathcal X(M)$, then $\mathcal G$ is called a \textbf{Hermitian metric}\index{Hermitian metric}. This means a Hermitian metric defines a Hermitian inner product on the tangent spaces $T_pM$ with respect to the almost complex structure. We call $\mathcal G$ a \textbf{Kähler metric}\index{Kähler!metric}, if the \textbf{fundamental 2-form} $\Phi: (\mathfrak X,\mathfrak Y)\mapsto \mathcal G(\mathfrak X,\mathcal J\mathfrak Y)$ is closed. If $M$ is a complex manifold and $\mathcal G$ is a Kähler metric, then we call $(M,\mathcal G)$ a \textbf{Kähler manifold}\index{Kähler!manifold}\index{manifold! Kähler}.
		\end{defi}
		
		Let $(V,\omega)$ be real symplectic $2n$-dimensional vector space with a complex structure $J$. We say that $J$ and $\omega$ are \textbf{compatible} if $J$ is a canonical transformation for $\omega$, that means
		$$\omega(Jv,Jw)=\omega(v,w)\quad\forall v,\ w\in V.$$
		If $J$ and $\omega$ are compatible, we can define a non-degenerate symmetric bilinear form $G$ on $V$ by putting
		\begin{equation} \label{eq:kaehlermetric}
			G(X,Y):=2\omega(X,JY)\quad\forall X,\ Y\in V.
		\end{equation}
		Then we obtain a Hermitian inner product on the complex vector space $V_J$ by the formula
			$$(X,Y)_J:=2\omega(X,JY)+i2\omega(X,Y).$$
		
		We will now construct the decomposition of $V_\mathbb C$ from Proposition \ref{prop:csccd} in a different manner. We define an injective map $p\colon V\to V_\mathbb C, \quad X\mapsto\frac{1}{2}(X-iJX)$. This allows us to identify $V$ with the the subspace $P_J:=\{X-iJX\}$ and the identification satisfies that $(.,.)_J$ coincides with the inner product $(.,.)_{P_J}:(Z,W)\mapsto 4i\omega(Z,\overline W)$ on $P_J$. Moreover, we find that $P_J$ is a Lagrangian subspace of $(V_\mathbb C,\omega)$, since it has dimension $2n$ and
		\begin{align*}
			\omega(X-iJX,Y-iJY) & = \omega(X,Y)-\omega(X,iJY)-\omega(iJX,Y)+\omega(iJX,iJY)\\
												  & = \omega(X,Y)+\omega(iJX,Y)-\omega(iJX,Y)-\omega(X,Y) = 0.
			\end{align*}
		
		Conversely, let us assume we have a Lagrangian subspace $P$ of $(V_\mathbb C,\omega)$ which satisfies $P\cap\overline{P}=\{0\}$. Then $V_\mathbb C=P\oplus\overline P$ and Proposition \ref{prop:csccd} yields a complex structure $J$ on $V$. We observe:
		
		\begin{prop}
			$\omega$ and the complex structure $JX:=iZ-i\overline Z$ for $X=Z+\overline{Z}$ from Proposition \ref{prop:csccd} are compatible.
		\end{prop}
		
		\begin{proof} \label{prop:sympcscompat}
			Let $Y=W+\overline{W}$. The vectors $Z,\ W$ are elements of $P$. Thus $\omega(Z,W)=\omega(\overline Z,\overline W)=0$. Then we can do the following calculation.
			\begin{align*}
				\omega(JX,JY) & = \omega(iZ+\overline{iZ},iW+\overline{iW}) \\
											& = -\omega(Z,W) + \omega(\overline Z,W) + \omega(Z, \overline W) - \omega(\overline Z,\overline W)\\
											& = \omega(Z,W) + \omega(\overline Z,W) + \omega(Z, \overline W) + \omega(\overline Z,\overline W)\\
											& = \omega(Z+\overline{Z},W+\overline{W}) \\
											& = \omega(X,Y) \qedhere
			\end{align*}
		\end{proof}
		
		This idea gives rise to the following definition.
		
		\begin{defi}
			Let $(M,\omega)$ be a symplectic manifold. A \textbf{Kähler polarization}\index{Kähler!polarization} $\mathcal P$ of $M$ is a complex subbundle of $TM_\mathbb C$ such that
			\renewcommand{\theenumi}{(\roman{enumi})}
			\renewcommand{\labelenumi}{\theenumi}
			\begin{enumerate}
				\item For all $p\in M$ we have that $\mathcal P_p$ is a Lagrangian subspace of $(T_pM_\mathbb C,\omega_p)$;
				\item For all $p\in M$ we have $\mathcal P_p\cap\overline{P_p}=\{0\}$; and
				\item the induced almost complex structure $\mathcal J$ is integrable.
			\end{enumerate}
			We say a complex vector field $\mathfrak X$ is tangent to $\mathcal P$, if $\mathfrak X(p)\in\mathcal P_p$ for all $p\in M$. In this case we will write $\mathfrak X\in\mathcal P$. A smooth complex function $f\in C^\infty(M,\mathbb C)$ is said to be \textbf{polarized}\index{function!polarized}, if
			$$\overline{\mathfrak X}(f)=0\quad\forall \mathfrak X\in\mathcal P.$$
			The space of polarized functions is denoted by $C^\infty_{\mathcal P}(M,\mathbb C)$\index{$C^\infty_{\mathcal P}(M,\mathbb C)$, polarized functions}.
		\end{defi}
		
		\begin{prop} \label{prop:kaehlercomplex}
			A Kähler polarization gives a symplectic manifold $(M,\omega)$ the structure of a Kähler manifold.
		\end{prop}
		
		\begin{proof}
			We have verified in Proposition \ref{prop:sympcscompat} that on each tangent space the almost complex structure $\mathcal J$ and the symplectic form $\omega$ are compatible. Hence we can pointwise define a Hermitian inner product by the formula (\ref{eq:kaehlermetric}). The smoothness and the invariance of the corresponding Riemannian metric $\mathcal G(.,.):=2\omega(.,\mathcal J.)$ are inherited from $\mathcal J$ and $\omega$. Since the fundamental 2-form is $\omega$ and hence closed, we have found that $\mathcal G$ is a Kähler metric. Furthermore, since $\mathcal J$ is integrable, we can apply the Newlander-Nirenberg Theorem to obtain a complex structure.
		\end{proof}
		
		\begin{prop} \label{prop:holoispola}
			If $\mathcal P$ is a Kähler polarization of $(M,\omega)$, then the holomorphic functions of the Kähler manifold $M$ are exactly the polarized functions with respect to $\mathcal P$.
		\end{prop}
		
		\begin{proof}
			We observe that for each $p\in M$ the subspaces $\mathcal P_p$ and $\overline{\mathcal P_p}$ of $T_pM_\mathbb C$ are the holomorphic and the antiholomorphic tangent space respectively. We have stated in the definition (cf. \cite[p.16]{gh}) of these spaces that a function $f$ on $M$ is holomorphic if and only if $v(f_p)=0$ for all $p\in M$ and $v\in\overline{\mathcal P_p}$.
		\end{proof}
		
	\subsection{Structure of Coadjoint Orbits} \label{sect:strcoadorb}
		
		Let $G$ be a Lie group and let $\lambda\in\mathfrak g^*$. We consider the coadjoint orbit $\mathcal O_\lambda=G/G_\lambda$. We have shown in Subsection \ref{sect:coadorb} that coadjoint orbits are smooth manifolds. Furthermore, if $G$ is a compact connected Lie group, we have seen in Corollary \ref{coro:coadgfm} that we can view $\mathcal O_\lambda$ as a generalized flag manifold. We remark that in case $G$ is compact and connected then so is $\mathcal O_\lambda$ since it is the continuous image of $G$ under the orbit map.
		
		The coadjoint action of $G$ on $\mathcal O_\lambda$ gives us for each $X\in\mathfrak g$ a vector field $\ad^*(X)$ on $\mathcal O_\lambda$ (cf. Definition \ref{defi:action} and Subsection \ref{ss:adcoad}). Theorem \ref{theo:hilgcoadorb} tells us that for each $\eta\in\mathcal O_\lambda$ the linear map $L_\eta\colon\mathfrak g\to T_\eta\mathcal O_\lambda,\quad X\mapsto\ad^*(X)\eta$ is surjective. Moreover, we can identify $T_\eta\mathcal O_\lambda\cong\mathfrak g/\mathfrak g_\eta$, where $\mathfrak g_\eta=\ker L_\eta$ is the Lie algebra of the stabilizer of $\eta$.
		
		The bilinear form $B_\eta\colon\mathfrak g\times\mathfrak g\to\mathbb R,\quad (X,Y)\mapsto\eta([X,Y])=(\ad^*(X)\eta)(Y)$ has null space $\mathfrak g_\eta$ and hence descends to a non-degenerate skew-symmetric bilinear form
		\begin{equation} \label{eq:sympform}
			\omega(\eta)\colon\mathfrak g/\mathfrak g_\eta\times\mathfrak g/\mathfrak g_\eta,\quad \omega(\eta)(\ad^*(X)\eta,\ad^*(Y)\eta)=\eta([X,Y]).
		\end{equation}		
		We define pointwisely a differential 2-form $\omega_\eta:=\omega(\eta)$ not yet known to be smooth.
		
		Let $\mu\colon\mathfrak g\to C^\infty(\mathcal O_\lambda,\mathbb R)$ be the assignment
		$$\mu(X)(\eta):=\left\langle \eta, X\right\rangle\quad\forall X\in\mathfrak g\quad\forall\eta\in\mathcal O_\lambda.$$
		
		\begin{prop} \label{prop:muprop1}
			$\mu$ satisfies $\ad^*(X)\mu(Y)=\mu([X,Y])$.
		\end{prop}
		
		\begin{proof}
			\begin{align*}
				(\ad^*(X)\mu(Y))(\eta) & = \frac{d}{dt} \mu(Y)(\Ad^*(\exp-tX)\eta) \\
															 & = \frac{d}{dt} \left\langle \Ad^*(\exp-tX)\eta, Y\right\rangle \\
															 & = \frac{d}{dt} \left\langle \eta, \Ad(\exp tX) Y\right\rangle \\
															 & = \left\langle \eta,[X,Y]\right\rangle \\
															 & = \mu([X,Y])
			\end{align*}
		\end{proof}
		
		\begin{prop}
			The 2-form $\omega$ defined by $\omega_\eta=\omega(\eta)$ is smooth.
		\end{prop}
		
		\begin{proof}
			For any $\eta\in\mathcal O_\lambda$ and $X,\ Y\in\mathfrak g$, we have
			\begin{align*}
				\omega_\eta(\ad^*(X)\eta,\ad^*(Y)\eta) & = \left\langle \eta,[X,Y]\right\rangle \\
																							 & = \mu([X,Y])(\eta)
			\end{align*}
			If we pick a $X_1,\dots,X_k\in\mathfrak g$ such that $\{\ad^*(X_1)(\eta),\dots,\ad^*(X_k)(\eta)\}$ form a basis of $T_\eta\mathcal O_\lambda$, then there is a neighborhood $U$ of $\eta$ in $\mathcal O_\lambda$ such that the vector fields $\ad^*(X_1),\dots\ad^*(X_k)$ span the tangent space at each point in $U$. Then we can use
			\begin{equation} \label{eq:omegaforvf}
				\omega(\ad^*(X),\ad^*(Y))=\mu([X,Y])
			\end{equation}
			 to express $\omega$ in terms of smooth functions.
		\end{proof}
		
		\begin{prop}
			The coadjoint action on $(\mathcal O_\lambda,\omega)$ is strongly symplectic.
		\end{prop}
		
		\begin{proof}
			We combine Proposition \ref{prop:muprop1} and (\ref{eq:omegaforvf}).
			\begin{align*}
				\left\langle i_{\ad^*(X)}\omega, \ad^*(Y)\right\rangle & = - \omega(\ad^*(Y),\ad^*(X)) \\
																															 & = \mu([Y,X]) \\
																															 & = \ad^*(Y)(\mu(X)) \\	
																															 & = \left\langle d(\mu(X)), \ad^*(Y)\right\rangle
			\end{align*}
		\end{proof}
		
		\begin{theo} \label{theo:coadhamilton}
			$\omega$ is a symplectic 2-form on $\mathcal O_\lambda$ and the coadjoint action is Hamiltonian with Lie algebra homomorphism $\mu$.
		\end{theo}
		
		\begin{proof}
			\cite[Theorem 8.4.9]{hga}.
		\end{proof}
		
		\begin{rema}
			For a Hamiltonian action $\sigma:G\times M\to M$ with Lie algebra homomorphism $\mu:\mathfrak g\to C^\infty(M,\mathbb R)$ the map $P:\mathcal O_\lambda\to\mathfrak g^*$ defined by
			$$<P(m),X>:=\mu(X)(m)$$
			is called the \textbf{moment map}\index{moment map} of $(M,\sigma)$. In the case of coadjoint orbits this map is the inclusion $\mathcal O_\lambda\hookrightarrow\mathfrak g^*$.
		\end{rema}
		
		\begin{prop} \label{prop:coadkaehler}
			Let $G$ be a compact connected Lie group and $\lambda\in\mathfrak g^*$. Then $\mathcal O_\lambda$ is a Kähler manifold.
		\end{prop}
		
		\begin{proof}
			According to Corollary \ref{coro:fots} we can assume that $\lambda\in\mathfrak t^*$ where $T$ is a maximal torus in $G$. We use Proposition \ref{prop:stabiscent} to find a torus $T_1$ such that $G_\lambda=C(T_1)$. Hence $\mathcal O_\lambda\cong G/C(T_1)$ is a generalized flag manifold. Let $\mathfrak b$ be a suitable subalgebra constructed as in Theorem \ref{theo:eopocot}, that is
			$$\mathfrak b:=\mathfrak g_{\lambda\mathbb C}\oplus\sum_{\alpha\in(\Delta^+\backslash\Delta_1^+)}\mathfrak g^\alpha.$$
			We put $\mathfrak p:=\sum_{\alpha\in(\Delta^+\backslash\Delta_1^+)}\mathfrak g^\alpha$ and use $(T_\lambda\mathcal O_\lambda)_\mathbb C\cong\mathfrak g_\mathbb C/\mathfrak g_{\lambda\mathbb C}$ to obtain a decomposition $T_\lambda\mathcal O_{\lambda\mathbb C}=\mathfrak p\oplus\overline{\mathfrak p}$.
			
			From Proposition \ref{prop:rootprop} and the definition of a $T_1$-admissible Weyl chamber we see that $\mathfrak p$ is closed under the Lie bracket. In particular $[\mathfrak p,\mathfrak p]\subseteq\mathfrak t^\bot$. Therefore we find
			$$\omega_\lambda(\mathfrak p,\mathfrak p)=\lambda([\mathfrak p,\mathfrak p])=0.$$
			This means that $\mathfrak p$ is in addition a Lagrangian subspace. We define $\mathcal P$ to be the $\Ad^*(G)$-invariant complex subbundle which has $\mathcal P_\lambda=\mathfrak p$. Then the associated almost complex structure $\mathcal J$ is the integrable almost complex structure from the proof of Proposition \ref{prop:csalg}. Since $\omega$ is $\Ad^*(G)$-invariant, we know that $\mathcal P$ is Lagrangian in every complexified tangent space and hence it is a Kähler polarization. Then $\mathcal P$ gives $\mathcal O_\lambda$ the structure of a Kähler manifold according to Proposition \ref{prop:kaehlercomplex}.
		\end{proof}
		
	\section{Prequantization}
	
	\subsection{\v{C}ech Cohomology}
		
		Let $M$ be a smooth manifold. We briefly recall the definition of the \v{C}ech cohomology $\check{H}^2(M,A)$, where $A=\mathbb Z$ or $A=\mathbb R$.
		
		Let $\mathcal U=\{U_j\}_{j\in J}$ be an open covering of $M$ and $k\in\mathbb N_0$. A \textbf{$\mathcal U$-$k$-cochain}\index{cochain} with values in $A$ is an assignment $c$ that attaches to each collection of $k+1$ sets $U_{j_0},\dots,U_{j_k}$ in the covering with non-empty intersection a number $c_{j_0,\dots,j_k}\in A$. We denote the set of all $\mathcal U$-$k$-cochains by $C^k(\mathcal U,A)$. It carries a natural additive group structure, where the sum of two chains $c,\ d$ is defined pointwisely as $(c+d)_{j_0,\dots,j_k}:=c_{j_0,\dots,j_k}+d_{j_0,\dots,j_k}$. We define the coboundary operator $\delta_k\colon C^k(\mathcal U,A)\to C^{k+1}(\mathcal U,A)$ by the formula
		$$(\delta_k c)_{j_0,\dots,j_{k+1}}:=\sum_{l=0}^{l+1}(-1)^l c_{j_0,\dots \hat{j_l},\dots,j_k}\quad \forall c\in C^k(\mathcal U,A).$$
		where $\hat{j}$ indicates that this element is omitted. It is a straightforward calculation to see that applying the coboundary operator twice always yields zero. That means the coboundary operator satisfies $\delta_{k+1}\circ\delta_k=0$. In particular, we have $\im\delta_k\subseteq\ker\delta_{k+1}$ and for $k\in\mathbb N$ we can define the following.
		$$\check{H}^k(\mathcal U,A):=\ker\delta_{k}/\im\delta_{k-1}.$$
		
		We introduce an ordering on the set of open coverings. Let $\mathcal V=\{V_l\}_{l\in L}$ and define $\mathcal V>\mathcal U$ if and only if $\mathcal V$ is a refinement of $\mathcal U$. In this case there is a map $\mu\colon L\to J$ such that $V_l\subseteq U_{\mu(l)}$ for all $V_l\in\mathcal V$. This yields a map $\mu_k\colon C^k(\mathcal U,A)\to C^k(\mathcal V,A)$ by putting
		$$(\mu_k c)_{l_0,\dots,l_k}:={c_{\mu(l_0),\dots,\mu(l_k)}}.$$
		Since $\mu_k$ commutes with the coboundary operator, it induces a map $\mu^k\colon \check{H}^k(\mathcal U,A)\to \check{H}^k(\mathcal V,A)$. One can prove that it is independent of the choice of $\mu$ and that if $\mathcal W>\mathcal V$, then the induced map for the refinement $\mathcal W>\mathcal U$ is just the concatenation of the induced maps of the intermediate steps (cf. \cite[5.33]{war}).
		
		We can then define the \textbf{\v{C}ech cohomology}\index{Cech cohomology} $\check{H}^k(M,A)$\index{$\check{H}^2(M,A)$, \v{C}ech cohomology} as the direct limit over the open coverings. That is, we form the disjoint union of all $H^k(\mathcal U,A)$ and factor out the following equivalence relation. If $c\in H^k(\mathcal U,A)$ and $d\in H^k(\mathcal V,A)$, then
		$$a\simeq b\quad\Leftrightarrow\quad \exists \mathcal W:\ \mathcal W>\mathcal U,\ \mathcal W>\mathcal V\text{ and }\mu^k_{\mathcal U,\mathcal W}(c)=\mu^k_{\mathcal U,\mathcal V}(d).$$
		This definition is cumbersome, but for well-behaved spaces we can circumvent it with the following proposition. Note that Theorem \ref{theo:countpara} asserts the existence of contractible coverings for smooth manifolds.
		
		\begin{prop}
			If $\mathcal U$ is a contractible covering of $M$, then
			$$\check{H}^k(M,A)\cong\check{H}^k(\mathcal U,A).$$
		\end{prop}
		
		\begin{proof}
			\cite[p.40]{gh}.
		\end{proof}
		
		Therefore we will consider contractible coverings $\mathcal U$ when dealing with $\check{H}^2(M,\mathbb Z)$ and $\check{H}^2(M,\mathbb R)$. We note that the inclusion $i\colon\mathbb Z\hookrightarrow\mathbb R$ gives a canonical homomorphism $i\colon C^k(\mathcal U,\mathbb Z)\to C^k(\mathcal U,\mathbb R)$. Moreover, the map between the cochains yields a homomorphism $\check{H}^2(M,\mathbb Z)\to\check{H}^2(M,\mathbb R)$ (cf. \cite[5.33]{war}).
		
	\subsection{Integrality Condition}
		
		\begin{theo} \label{theo:lbisochern} 
			Let $M$ be a smooth manifold. There is a one-to-one correspondence between $\mathcal L(M)$ and $\check{H}^2(M,\mathbb Z)$.
		\end{theo}
		
		\begin{proof}
			We will give a map $\kappa\colon\mathcal L(M)\to\check{H}^2(M,\mathbb Z)$. The presentation follows \cite[§1.2]{limaaa}. Let $L$ be a line bundle over $M$ and $\{U_j,s_j\}$ a local system with transition functions $c_{jk}$. We use Theorem \ref{theo:countpara} to assume that $\mathcal U:=\{U_j\}_{j\in J}$ is a contractible covering. 
			
			Let $U_j\cap U_k$ be non-empty. Then it is contractible (cf. Definition \ref{defi:coverings})and we have a smooth function $c_{jk}\colon U_j\cap U_k\to\mathbb C$. Because $c_{jk}$ is continuous and non-vanishing, the image $c_{jk}(U_j\cap U_k)$ is contractible hence simply connected and it does not contain zero. Therefore we can pick logarithms to define the following functions.
			$$f_{jk}\colon U_j\cap U_k\to\mathbb C,\quad f_{jk}:=\frac{1}{2\pi i}\log c_{jk}$$
			If the intersection $U_j\cap U_k\cap U_l$ is non-empty we can define a continuous function
			\begin{equation} \label{eq:cocycfunc}
				a_{jkl}:=f_{jk}+f_{kl}-f_{jl}.
			\end{equation}
			on $U_j\cap U_k\cap U_l$. Using the cocycle conditions (\ref{eq:cc}) we find $\exp(2\pi ia_{jkl})=c_{jk}c_{kl}c_{jl}^{-1}=1$ and therefore $a_{jkl}$ has values in $\mathbb Z$. Since it is continuous and $\mathbb Z$ is discrete, it is actually a locally constant function. We can therefore define a $\mathcal U$-$2$-cochain by
			$$a:(U_j,U_k,U_l)\mapsto a_{jkl}.$$
			We calculate
			\begin{align*}
				(\delta_2 a)_{jklm} & = a_{klm}						  -a_{jlm}			 		   +a_{jkm}			  		  -a_{jkl} \\
														& = f_{kl}+f_{lm}-f_{km}-f_{jl}-f_{lm}+f_{jm}+f_{jk}+f_{km}-f_{jm}-f_{jk}-f_{kl}+f_{jl} \\
														& = 0.
			\end{align*}
			This shows that $a\in\ker\delta_2$ and hence it defines an element $[a]$ in $\check{H}^2(\mathcal U,\mathbb Z)\cong\check{H}^2(M,\mathbb Z)$. Wallach states in \cite[Lemma 2.3]{wal} that $[a]$ is independent of the choices made in its definition. Hence we have defined a map $\kappa\colon [L]\mapsto [a]$.
			
			Let $\{h_j\}_{j\in J}$ be a partition of unity subordinate to $\mathcal U$ with same index. If $a\in C^2(\mathcal U,\mathbb Z)$ is a cochain, and $j,\ k$ are such that $U_j\cap U_k\neq\emptyset$ we can define smooth functions $f_{jk}\in C^\infty(U_j\cap U_k)$ by the formula
			$$f_{jk}:=\sum_{l\in J} a(U_j,U_k,U_l) h_l.$$
			If $\delta^2a=0$, then the functions $f_{jk}$ satisfy (\ref{eq:cocycfunc}). Since $a$ is integer-valued, the functions $c_{jk}:=\exp(2\pi if_{jk})$ then satisfy the cocycle conditions and hence define a complex line bundle $L$ over $M$. From the construction we see $\kappa([L])=[a]$ which proves surjectivity.
			
			Let $L_1$ and $L_2$ be complex line bundles with local systems $\{U_j,s_j^1\}_{j\in J},\ \{U_j,s_j^2\}_{j\in J}$ respectively. We denote by $c_{jk}^1,\ c_{jk}^2$ and $f_{jk}^1,\ f_{jk}^2$ the transition functions and the respective logarithms of $L_1$ and $L_2$ from the construction. Let us assume $\kappa([L_1])=\kappa([L_2])$, then the functions $f_{jk}:=f_{jk}^1-f_{jk}^2$ satisfy $f_{jk}+f_{kl}-f_{jl}=0$. Using the partition of unity, we define smooth functions on $U_j$
			$$t_j:=\sum_{k\in J}f_{kj} h_k$$
			which satisfy $t_k-t_l=f_{jl}$. Finally we observe that the family of functions $g_j\colon U_j\to\mathbb C^*$, $g_j:=\exp(2\pi i)t_j$ fulfills the requirement of Proposition \ref{prop:eqiftfrel}. This shows $[L_1]=[L_2]$ and thus we have proven injectivity.
		\end{proof}
		
		\begin{defi}
			Let $L$ be a line bundle over a smooth manifold $M$ and $\kappa$ the map from the proof of Theorem \ref{theo:lbisochern}. Then the element $\kappa(\left[L\right])$ is called the \textbf{Chern class}\index{Chern class of a line bundle} of $L$ and $\left[L\right]$ accordingly. Theorem \ref{theo:lbisochern} states that conversely, the Chern class determines its complex line bundle uniquely up to equivalence.
		\end{defi}
		
		\begin{rema}
			If $L_1$ and $L_2$ are complex line bundles over $M$ then we define $L_1\otimes L_2$ as the fiberwise tensor product which again is a complex line bundle over $M$. Then $\mathcal L(M)$ together with the multiplication $\cdot:(\left[L_1\right],\left[L_2\right])\mapsto[L_1\otimes L_2]$ is a group known as the \textbf{Picard group}\index{Picard group} of $M$. With respect to this group structure the map $\kappa$ from the proof of Theorem \ref{theo:lbisochern} is a group isomorphism (cf. \cite[p.133]{gh} and \cite[p.91]{limaaa}).
		\end{rema}
		
		\begin{theo}[de Rham Isomorphism]\index{de Rham Isomorphism} \label{theo:derhamiso}  
			Let $M$ be a smooth manifold. The \v{C}ech cohomology $\check{H}^2(M,\mathbb R)$ is isomorphic with the de Rham cohomology $H^2_{dR}(M,\mathbb R)$.
		\end{theo}
		
		\begin{proof}[Sketch of proof]
			We only construct a map from $H^2_{dR}(M,\mathbb R)$ to $\check{H}^2(M,\mathbb R)$. Many arguments are similar to the proof of Theorem \ref{theo:lbisochern}. For a general proof of the de Rham Theorem confer \cite[§5]{war} or \cite[p.43]{gh}. A proof for our specific situation can be found in \cite[Theorem 1.3]{wal}.
			
			Let $\omega\in[\omega]$ be a real closed 2-form and $\mathcal U:=\{U_j\}_{j\in J}$ a contractible open covering of $M$. We apply the Poincaré Lemma to obtain for each $j\in J$ a real 1-form $\alpha_j$ satisfying
			$$\omega_{|U_j}=d\alpha_j\in\Omega^1(U_j).$$
			If $j,\ k$ are such that $U_j\cap U_k\neq\emptyset$ then the forms $d{\alpha_j},\ d{\alpha_k}$ coincide on the intersection and therefore $(\alpha_k-\alpha_j)_{|U_j\cap U_k}$ is closed. The intersection $U_j\cap U_k$ is contractible (cf. Definition \ref{defi:coverings}) and hence we can apply the Poincaré Lemma once more to obtain functions $f_{jk}$ such that
			$$(\alpha_k-\alpha_j)_{|U_j\cap U_k}=df_{jk}.$$
			Then on a non-empty intersection $U_j\cap U_k\cap U_l\neq\emptyset$ we have
			$$df_{jk}+df_{kl}-df_{jl}=\alpha_j-\alpha_k+\alpha_k-\alpha_l-\alpha_j+\alpha_l=0.$$
			Therefore the function $a_{jkl}:=f_{jk}+f_{kl}-f_{jl}$ defined on the intersection is a real constant. We define a $\mathcal U$-$2$-cochain $a$ by putting $a:(U_j,U_k,U_l)\mapsto a_{jkl}$. Then we have $a\in\ker \delta^2$ which follows from the same calculation as in the proof of Theorem \ref{theo:lbisochern}. Hence $a$ determines a class $[a]\in\check{H}^2(M,\mathbb R)$.
		\end{proof}
		
		\begin{defi}
			Let $M$ be a smooth manifold and $\omega$ a closed 2-form on $M$. Then we denote by $\mathcal L_C(M,\omega)$\index{$\mathcal L_C(M,\omega)$} the set of equivalence classes of line bundles with connection which have curvature $\omega$ and which admit a compatible Hermitian structure. This is well-defined because of Proposition \ref{prop:elbwchtsc}.
		\end{defi}
		
		\begin{prop} \label{prop:cherncurv}
			Let $L_1$ and $L_2$ be complex line bundles over $M$ with Hermitian connections $\nabla^1$, $\nabla^2$ and curvatures $\omega^1$, $\omega^2$ respectively. Then the image of $\left[\omega^1\right]$ under the inverse of the de Rham isomorphism is the Chern class of $L_1$. In particular, if $\left[\omega^1\right]=\left[\omega^2\right]$, then $L_1$ and $L_2$ are equivalent as complex line bundles.
		\end{prop}
		
		\begin{proof}
			Let $(L,\nabla)$ be a line bundle with Hermitian connection over $M$ and let $\omega$ be its curvature. Proposition \ref{prop:curvature} states that $\omega$ is closed and hence it determines a class $[\omega]\in H^2_{dR}(M,\mathbb R)$. Let $\mathcal U:=\{U_j\}_{j\in J}$ be a contractible covering of $L$ and let $\{U_j,s_j,\alpha_j\}_{j\in J}$ be a local form of $\nabla$. By rescaling the sections $s_j$ we can assume that $|H|^2\circ s_j=1$. Then the following calculations shows that the forms $\alpha_j$ are real. On $U_j$ we have:
			\begin{align*}
				0 = \mathfrak X(|H|^2\circ s_j)		 & = H(\nabla_\mathfrak X s_j,s_j) + H(s_j,\nabla_\mathfrak X s_j) \\
																					 & = H(\nabla_\mathfrak X s_j,s_j) + \overline{H(\nabla_\mathfrak X s_j,s_j)} \\
																					 & = 2\pi i (H(\alpha_j(X)s_j,s_j) + \overline{H(\alpha_j(X)s_j,s_j)}) \\
																					 & = 2\pi i (\alpha_j(\mathfrak X) - \overline{\alpha_j(\mathfrak X)}) H(s_j,s_j)\quad\forall\mathfrak X\in\mathcal X(M).
			\end{align*}
			Let $c_{jk}$ denote the transition functions of the local system $\{U_j,s_j\}_{j\in J}$. As in the proof of Theorem \ref{theo:lbisochern} we define $f_{jk}:=\frac{1}{2\pi i}\log c_{jk}$ to determine an integer $\mathcal U$-$2$-cochain $a$ by putting $a_{jkl}:=f_{jk}+f_{kl}-f_{jl}$. Then $[a]$ is the Chern class of $L$. On the other hand, we have from Proposition \ref{prop:proplocform} that $\alpha_j-\alpha_k=\frac{1}{2\pi i}\frac{dc_{jk}}{c_{jk}}=df_{jk}$. Moreover by definition of the curvature $\omega_{|U_j}=d\alpha_j$. This construction inverts the assignment in the de Rham isomorphism and hence $[\omega]$ is the preimage of the Chern class.
		\end{proof}
		
		\begin{defi}
			We call a closed 2-form $\omega$ \textbf{integral}\index{form!integral 2-form}, if its class $[\omega]\in H^2_{dR}(M,\mathbb R)$ lies in the image of $\check{H}^2(M,\mathbb Z)$ under the de Rham isomorphism. A coadjoint orbit with an integral symplectic form is called an \textbf{integral orbit}\index{orbit!integral}.
		\end{defi}
		
		\begin{theo}[Integrality Condition]\index{Integrality Condition} \label{theo:intcond} 
			Let $M$ be a smooth manifold and $\omega$ a closed real 2-form on $M$. Then $\mathcal L_C(M,\omega)$ is not empty if and only if $\left[\omega\right]\in H^2_{dR}(M,\mathbb R)$ is integral.
		\end{theo}
		
		\begin{proof}
			We have shown in Proposition \ref{prop:cherncurv} that if $(L,\nabla)\in\mathcal L_C(M,\omega)$ then the class $[\omega]$ is integral. For the other direction let $\omega$ be an integral closed 2-form on $M$ and $\{U_j\}_{j\in J}$ a contractible covering. We proceed as in the proof of Theorem \ref{theo:derhamiso} to obtain functions $f_{jk}$ on non-empty intersections $U_j\cap U_k$ such that $f_{jk}+f_{kl}-f_{jl}$ is constant where defined. The integrality of $\omega$ then means that this constant is an integer. We define functions $c_{jk}:=\exp(2\pi if_{jk})$. Then the integrality yields that the functions $c_{jk}$ satisfy the cocycle conditions (\ref{eq:cc}). Then Proposition \ref{prop:tfgiveclb} tells us that there is a complex line bundle $L$ with a local system $\{U_j,s_j\}_{j\in J}$ having the $c_{jk}$ as transition functions. This also means that the Chern class of $L$ is the preimage of $[\omega]$ under the de Rham isomorphism. By construction, the 1-forms and the transition functions satisfy
			$$\alpha_j-\alpha_k=df_{jk}=\frac{1}{2\pi i}\frac{dc_{jk}}{c_{jk}}$$
			and by Proposition \ref{prop:locformglobform} there is a uniquely determined connection $\nabla$ on $L$ with curvature $\omega$. We define a Hermitian structure $H$ on $L$ by the formula
			$$H_p(x,y):= \frac{x}{s_j(p)}\overline{\left(\frac{y}{s_j(p)}\right)}\quad\forall p\in U_j\quad\forall x,y\in L_p.$$
			Here $j$ is such that $p\in U_j$. The definition is independent of the choice of $j$, since the transition functions have absolute value equal to 1. Therefore if we replace $s_j$ by $s_k=c_{jk}s_j$ we find that the value of $H_p$ is unchanged. The following calculations show that $H$ is compatible with $\nabla$. We consider the following expression restricted to $U_j$ where each sections $s,\ t\in C^\infty(M,L)$ are given by $s=fs_j$ and $t=gs_j$.
			\begin{align*}
				\mathfrak X(H(s,t)) & =	\mathfrak X(H(fs_j,gs_j)) \\
														& = \mathfrak X(f\overline{g})\\
														& = (\mathfrak X f)\overline{g}+f(\mathfrak X\overline g)
			\end{align*}
			On the other hand,
				$$H(\nabla_\mathfrak X s,t) = H\left( (\mathfrak X f)s_j + 2\pi i f \alpha_j(\mathfrak X)s_j, g s_j \right) = (\mathfrak X f)\overline g + 2\pi i f \alpha_j(\mathfrak X)\overline{g}$$
			and
				$$H(s, \nabla_\mathfrak X t) = H\left( fs_j, (\mathfrak X g)s_j + 2\pi i g \alpha_j(\mathfrak X)s_j \right) = f (\mathfrak X \overline g) - f 2\pi i \overline g \alpha_j(\mathfrak X)$$
			since the forms $\alpha_j$ are real. Therefore
				$$\mathfrak X(H(s,t))=H(\nabla_\mathfrak X s,t)+H(s,\nabla_\mathfrak X t).$$
		\end{proof}
		
		\begin{rema}
			If $L_1$ is a complex line bundle over $M$ and $\left[\omega^1\right]$ is the image of the Chern class of $L_1$ under the de Rham isomorphism, then any $\omega\in\left[\omega^1\right]$ is integral by definition. According to Theorem \ref{theo:intcond} there exists a complex line bundle $L_2$ with connection which has curvature $\omega$. Using Proposition \ref{prop:cherncurv} we find that $L_1$ and $L_2$ have the same Chern class and are therefore equivalent as complex line bundles.
		\end{rema}
		
	\subsection{Prequantum Bundle} \label{ss:prequant}
		
		Let $G$ be a compact connected Lie group and let $(\mathcal O_\lambda,\omega)$ be an integral coadjoint orbit. Again we assume 
		$\lambda\in\mathfrak t^*$ for a maximal torus $T$. From Theorem \ref{theo:intcond} we know that there is a line bundle $(L,\nabla,H)$ with Hermitian connection and curvature $\omega$ over $\mathcal O_\lambda$. Let $\alpha$ be the associated connection form (cf. Proposition \ref{prop:connform}). The following construction is from \cite[§2]{limaaa}.
		
		Let $x\in L^*$. We denote by $\Ver_x(L)$\index{$\Ver_x(L)$} the 2-dimensional tangent space of $L^*_{\pi(x)}$ at $x$ and we set $\Hor_x(L):=\ker\alpha_x$\index{$\Hor_X(L)$}. Then Kostant shows in \cite[Proposition 2.6.1]{limaaa}
		$$T_x(L)=\Ver_x(L)\oplus \Hor_x(L)\quad\forall x\in L^*.$$
		Let $e(L)$\index{$e(L)$} be the Lie algebra of vector fields on $L^*$ which commute with the action of $\mathbb C^*$ on $L^*$. Then Kostant uses the above decomposition to show that $e(L)$ is parametrized by $\delta\colon C^\infty(M,\mathbb C)\times\mathcal X(M)\to e(L)$. This is realized by showing that the horizontal part corresponds to real vector fields on $M$ and that the vertical part may be described by picking a non-zero direction in each fiber (cf. \cite[Proposition 2.9.1]{limaaa}). We also obtain that $\mathfrak Z:=\delta(\phi,\mathfrak X)$ is globally integrable if and only if $\mathfrak X$ is (cf. \cite[Theorem 2.10.1]{limaaa}). We denote by $e(L,\nabla)$\index{$e(L,\nabla)$} the set of $\mathfrak Z\in e(L)$ such that $\mathcal L_{\mathfrak Z}\alpha=0$ and $\mathfrak Z(|H|^2)=0$.
		
		\begin{prop} \label{prop:derivsigmal}
			If $\mathfrak Z\in e(L)$ is globally integrable, then $\mathfrak Z\in e(L,\nabla)$ if and only if the corresponding 1-parameter group of diffeomorphisms lies in $E(L,\nabla)$. An element $\mathfrak Z=\delta(\phi,\mathfrak X)$ lies in $e(L,\nabla)$ if and only if $\phi$ is real and $i_{\mathfrak X}\omega=d\phi$.
		\end{prop}
		
		\begin{proof}
			\cite[Theorem 3.3.1]{limaaa}.
		\end{proof}
		
		We have shown in Theorem \ref{theo:coadhamilton} that the coadjoint action of $G$ on $\mathcal O_\lambda$ is Hamiltonian. Then for any $X\in\mathfrak g$ we obtain a smooth real function $\mu(X)\in C^\infty(\mathcal O_\lambda,\mathbb R)$ and a Hamiltonian vector field $X_{\mu(X)}$ on $\mathcal O_\lambda$ associated to the 1-form $i_{\mathfrak X_{\mu(X)}}\omega$. Then using the parametrization $\delta$ we can associate to $X\in\mathfrak g$ an element $\delta(\mu(X),\mathfrak X_{\mu(X)})\in e(L,\nabla)$. The existence of a mapping between $\mathfrak g$ and $e(L,\nabla)$ is a necessary condition for a smooth action of $G$ on $L$ by isometric isomorphisms of line bundles with Hermitian connection. We will proof the existence of such an action which has this mapping as differential for the case that $G$ is simply connected (cp. \cite[Theorem 4.5.1]{limaaa}).
		
		\begin{theo}
			Let $(M,\omega)$ be a symplectic manifold with integral symplectic form and let $(L,\nabla,H)\in\mathcal L_C(M,\omega)$. Let $G$ be a simply connected Lie group and assume that $M$ is a Hamiltonian $G$-space with strongly symplectic action $\sigma$ and Lie algebra homomorphism $\mu$. Then $\sigma$ may be lifted to a smooth action $\sigma_L\colon G\to E(L,\nabla)$ such that $L$ is a homogeneous line bundle. Moreover, $\sigma_L$ may be chosen such that
			$$L{\sigma}_L(X)=\delta(\mu(X),\mathfrak X_{\mu(X)})\quad\forall X\in\mathfrak X.$$
		\end{theo}
		
		\begin{proof}[Sketch of proof.]
			The assignment $X\mapsto\delta(\mu(X),\mathfrak X_{\mu(X)})$ is a homomorphism of Lie algebras (cf. \cite[Theorem 4.2.1]{limaaa}). Since the action is Hamiltonian, we have $\mathfrak X_{\mu(X)}=L{\sigma}(X)$. This vector field is globally integrable and therefore $\delta(\mu(X),\mathfrak X_{\mu(X)})$ is globally integrable (cf. \cite[Theorem 2.10.1]{limaaa}). Then \cite[§IV, Theorem III]{pal}  tells us that there is a smooth action $\sigma_L$ of $G$ on $L^*$ whose derivative is $\delta(\mu(X),\mathfrak X_{\mu(X)})$. From Proposition \ref{prop:derivsigmal} we see that $\sigma_L(G)\subseteq E(L,\nabla)$.
		\end{proof}
		
		In the following, we will work with integral coadjoint orbits $(\mathcal O_\lambda,\omega)$ which admit a line bundle with Hermitian connection and whose coadjoint action lifts to $L$ in the way described above. We will give a criterion when integral coadjoint orbits allow for lifts.
		
		If one is more interested in the coadjoint orbits than in the group $G$, then the situation is simpler. The following consideration from \cite[6.3.2]{walhomo} shows that any coadjoint orbit of a compact connected Lie group is also a coadjoint orbit of a simply connected compact Lie group. 
		
		\begin{prop}
			Let $G$ be a compact connected Lie group and let $Z:=\{z\in G\ |\ gz=zg\quad\forall g\in G\}$ be center of $G$. Then $Z$ is contained in any maximal torus and in particular, we have $Z\subseteq G_\lambda$. Then the Isomorphism Theorem tells us $\mathcal O_\lambda=G/G_\lambda=(G/Z)/(G_\lambda/Z)$. Let $\tilde{G}$ be a simply connected universal covering group of $G$ (cp. \cite[§1.11]{knapp}) and let $\tilde{G_\lambda}$ be the connected subgroup corresponding to $G_\lambda$. Then $\mathcal O_\lambda=\tilde{G}/\tilde{G_\lambda}$ and $\tilde{G}$ is compact.
		\end{prop}
				
	\section{Quantization}
	
	\subsection{Integrality Condition Revisited}
		
		\begin{defi}
			Let $(\mathcal O_\lambda,\omega)$ be a coadjoint orbit and denote by $\sigma$ the coadjoint action of $G$ on $\mathcal O_\lambda$. We say that $(\mathcal O_\lambda,\omega)$ is \textbf{strongly integral}\index{orbit!strongly integral} if there exists a $(L,\nabla,H)\in\mathcal L_C(\mathcal O_\lambda,\omega)$ and a smooth action $\sigma_L\colon G\times L\to L$ such that $\sigma_L(G)\subseteq E(L,\nabla)$ making $L$ a homogeneous  line bundle. Furthermore we demand $L{\sigma}_L(X)=\delta(\mu(X),\mathfrak X_{\mu(X)})$ for all $X\in\mathfrak X$ (cp. Subsection \ref{ss:prequant} or \cite[§3]{limaaa}). Here $\mu:\mathfrak g\to C^\infty(\mathcal O_\lambda)$ is the Lie algebra homomorphism of the Hamiltonian action $\sigma$ and $L{\sigma}(X)=\mathfrak X_{\mu(X)}$ holds for all $X\in\mathfrak g$.
		\end{defi}
		
		The following theorem is adapted from \cite[Theorem 5.7.1]{limaaa}.
		
		\begin{theo} \label{prop:intcondrev}
			Let $G$ be a compact connected Lie group and $T$ a maximal torus. Let $(\mathcal O_\lambda,\omega)$ be a coadjoint orbit and assume $\lambda\in\mathfrak t^*$ which is no restriction according to Corollary \ref{coro:fots}. Then $(\mathcal O_\lambda,\omega)$ is strongly integral if and only there exists a character $\tau\colon T\to S^1$ such that $L\tau=2\pi i\lambda\colon\mathfrak t\to\mathbb R$.
		\end{theo}
		
		\begin{proof}[Sketch of proof]
			Let $(L,\nabla,H)$ be a line bundle with Hermitian connection and curvature $\omega$. Let $\sigma$ denote the coadjoint action and $\sigma_L\colon G\to E(L,\nabla)$ its lifted version. The subgroup $\sigma_L(G_\lambda)$ induces isometric linear automorphisms on the fiber $L_\lambda$. Hence we obtain a character $\chi\colon G_\lambda\to S^1$ by defining $\chi(a)$ via the formula
			$$\chi(a)v=\sigma_L(g)v\quad\forall v\in L_\lambda.$$
			This is well-defined since the fiber is one-dimensional. Since the action is linear and isometric, it suffices to calculate $\chi(a)$ for a single $v$ and $\chi(a)$ has absolute value 1. Kostant shows that for another complex line bundle with Hermitian connection which is equivalent to $L$ the equivalence diffeomorphism intertwines the actions of $G$. Therefore $\chi$ actually only depends on $[L,\nabla]$.
			
			We now show that $L\chi=2\pi i\lambda$. Fix $v\in L_\lambda$ and let $X\in\mathfrak g_\lambda$. Then the vector field $L{\sigma}(X)=\ad^*(X)$ vanishes at $\lambda$. Let $\gamma(t):=\sigma_L(\exp -tX)v$ be the local flow of $L{\sigma}_L(X)=\delta(\mu(X),L{\sigma}(X))$ through $v$. Then from \cite[Theorem 2.10.1]{limaaa}, which describes the local flows of vector fields parametrized by $\delta$, we find that $\gamma$ is of the form
			$$\gamma(t)=e^{2\pi i t \lambda(X)}v.$$
			Therefore we have established $\chi(\exp tX)v=e^{2\pi i t \lambda(X)}v$ which proves the claim. Since $G_\lambda$ is connected (cf. Proposition \ref{prop:stabisconn}), the character $\chi$ is uniquely determined by $\lambda$ (cf. Proposition \ref{prop:uniqueforconnected}) and independent of the choice of $(L,\nabla,H)$. Using Proposition \ref{prop:11char} we find that $2\pi i\lambda$ lifts to a character $\tau$ of $T$.
			
			Now let us assume we are given a character $\tau\colon T\to S^1$ of $T$ such that $L\tau=2\pi i\lambda$. We extend $\tau$ to $\chi\colon G_\lambda\to S^1$ by using Proposition \ref{prop:11char}. Then we can apply Construction \ref{cons:homolb} to obtain the homogeneous line bundle $L:=G\times_\chi\mathbb C$ over $G/G_\lambda$. It has a canonical $G$-action defined by $\sigma_L(a)(g,z)G_\lambda:=(ag,z)G_\lambda$. On the other hand, $\mathbb C^*$ operates on $L$ by $c\cdot(g,z)G_\lambda:=(g,cz)G_\lambda$ and this assignment commutes with the action of $G$. Therefore, if we define $H:=G\times\mathbb C^*$ we obtain a surjection
			$$\nu: H\to L^*,\quad\nu(g,z):=z\cdot\sigma_L(g)(e,1)G_\lambda=(g,z)G_\lambda.$$
			$H$ acts on itself by left translation and multiplication, respectively. We define an $H$-invariant form $\delta$ by
			$$\delta:=(\tilde{\lambda},\frac{1}{2\pi i}\frac{dz}{z}).$$
			Here $\tilde{\lambda}$ denotes the left-invariant 1-form on $G$ whose value at $e$ is $\lambda$. Kostant shows in \cite[p.200]{limaaa} that there is a unique connection form $\alpha$ on $L^*$ which is $H$-invariant and satisfies $\nu^*\alpha=\delta$. As in Proposition \ref{prop:connform} there is a unique connection $\nabla$ on $L$ associated to $\alpha$. Then \cite[Proposition 5.7.2]{limaaa} asserts that the curvature of $\nabla$ is $\omega$.
		\end{proof}
		
	\subsection{Orbit Representations}
		
		Let $G$ be a compact connected Lie group and let $(\mathcal O_\lambda,\omega)$ be a strongly integral coadjoint orbit. Choose a homogeneous line bundle $(L,\nabla,H)$ over $\mathcal O_\lambda$ with curvature $\omega$ such that $G$ acts on $L$ by connection preserving isometries. Moreover, let $\mathcal P$ be an invariant Kähler polarization whose existence we proved in Proposition \ref{prop:coadkaehler}.
		
		\begin{prop} \label{prop:Gaotsops}
			Let $\sigma_S$ be the associated action of $G$ on the space of smooth sections from Definition \ref{defi:actonsect}. Then the space of polarized sections $C^\infty_\mathcal P(M,L):=\{s\in C^\infty(M,L)\ |\ \nabla_{\mathfrak X}s=0\quad\forall \mathfrak X\in\mathcal P \}$\index{$C^\infty_\mathcal P$, polarized sections} is $\sigma_S(G)$-invariant.
		\end{prop}
		
		\begin{proof}
			Let $s$ be a polarized section, that is, $\nabla_\mathfrak X s=0$ for all $\mathfrak X\in\mathcal P$. Note that the invariance of $\mathcal P$ means that $\mathcal P$ is defined by $\sigma(G)$-invariant vector fields $\mathfrak X\in\mathcal X_\mathbb C(\mathcal O_\lambda)$ such that $\mathfrak X(\lambda)\in\mathcal P_\lambda\subset (T_\lambda\mathcal O_\lambda)_{\mathbb C}$. Since $\sigma_L(G)\subseteq E(L,\nabla)$, we can use  (\ref{eq:elnablaprop}) from Definition \ref{defi:elnabla}.
			\begin{align*}
				\nabla_\mathfrak X(\sigma_S(g)s) & = \nabla_\mathfrak X (\sigma_L(g)\circ s\circ\sigma(g^{-1})) \\
																				 & = \sigma_L(g)\circ(\nabla_\mathfrak X s)\circ\sigma(g^{-1}) \\
																				 & = 0
			\end{align*}
			for all $\mathfrak X\in\mathcal P$ and all $g\in G$.
		\end{proof}
		
		\begin{prop}
			There exists a $G$-invariant measure on $\mathcal O_\lambda$.
		\end{prop}
		
		\begin{proof}
			Helgason proves in \cite[§I, Theorem 1.9]{helgason} that there exists an invariant measure on the coset space $G/H$ if and only if
			$$|\det\Ad_G(h)|=|\det\Ad_H(h)|,\quad\forall h\in H.$$
			We note that $h\mapsto |\det\Ad(h)|$ is a continuous group homomorphism and that in our case $H=G_\lambda$ is compact. Therefore the image is a compact subgroup of $(\mathbb R,\cdot)$ and the above expressions are always equal to 1.
		\end{proof}
		
		\begin{prop} \label{prop:coadunit}
			The space of polarized sections carries a natural Hilbert space structure such that the action of $G$ is a unitary representation.
		\end{prop}
		
		\begin{proof}
			We define a Hermitian inner product on $C^\infty_\mathcal P(M,L)$ by putting
			$$\left\langle s,t\right\rangle:=\int_{\mathcal O_\lambda} H(s,t)(\eta)d\mu(\eta).$$
			This is well-defined, since $H(s,t)$ is smooth and hence measurable. Since $\mathcal O_\lambda$ is compact, the integral is finite. To see that this bilinear form is positive definite, we note that a non-negative function which is not the zero function has positive integral since $\mu$ is a Borel measure.  Using $\sigma_L(G)\subset E(L,\nabla)$ we find:
			\begin{align*}
				\left\langle \sigma_S(g)s,\sigma_S(g)t\right\rangle & = \int_{\mathcal O_\lambda} H(\sigma_L(g)s(\sigma(g^{-1})\eta),\sigma_L(g)t(\sigma(g^{-1})\eta))d\mu(\eta) \\
				& = \int_{\mathcal O_\lambda} H(s(\sigma(g^{-1})\eta),t(\sigma(g^{-1})\eta))d\mu(\eta) \\
				& = \int_{\mathcal O_\lambda} H(s,t)(\eta)d\mu(\eta) \\
				& = \left\langle s,t\right\rangle
			\end{align*}
			Thus $G$ acts unitarily with respect to $\left\langle .,.\right\rangle$. Woodhouse proves in \cite[p.286]{woo} that the space of polarized sections is topologically closed. We will show in the following theorem that in our case this space is finite-dimensional.
		\end{proof}
		
		\begin{cons} \label{cons:orbitmethod}
			Let $G$ be a compact connected Lie group, $T$ a maximal torus in $G$ and $\tau\colon T\to\mathbb C$. Let $2\pi i\lambda:=L\tau\in\mathfrak t^*$. Then the coadjoint orbit $(\mathcal O_\lambda,\omega)$ is strongly integral according to Theorem \ref{prop:intcondrev}. Choose a homogeneous line bundle $(L,\nabla,H)$ over $\mathcal O_\lambda$ with curvature $\omega$ such that $G$ acts on $L$ by connection preserving isometries. Let $T_1$ be torus such that $G_\lambda=C(T_1)$ and $\mathcal O_\lambda=G/G(T_1)$ which is possible according to Corollary \ref{coro:coadgfm}. We choose a $T_1$-admissible Weyl chamber to obtain a suitable subalgebra $\mathfrak b$ of $\mathfrak g_\mathbb C$ as described in Theorem \ref{theo:eopocot}. We use $\mathfrak b$ and Proposition \ref{prop:coadkaehler} which gives $\mathcal O_\lambda$ the structure of a Kähler manifold. We define the space of polarized sections $C^\infty_\mathcal P(M,L)$ and use Propositions \ref{prop:Gaotsops} and \ref{prop:coadunit} to obtain a unitary representation $(\pi_{\mathcal O_\lambda},C^\infty_\mathcal P(M,L))$.
		\end{cons}
		
		\begin{prop} \label{prop:holopola}
			Let $G$ be a compact connected Lie group, $T$ a maximal torus in $G$, $\tau\colon T\to\mathbb C$ and $2\pi i\lambda:=L\tau$. Let $\chi\colon G_\lambda\to S^1$ be the character of $G_\lambda=C(T_1)$ associated by Proposition \ref{prop:11char}. Then Construction \ref{cons:gborel-weil} for $C(T_1)$ and Construction \ref{cons:orbitmethod} for $\mathcal O_\lambda=G/C(T_1)$ are both applicable. If we choose in both constructions the same Weyl chamber, then $(\pi_\chi,\Gamma(G/T,G\times_\chi\mathbb C))$ and $(\pi_{\mathcal O_\lambda},C^\infty_\mathcal P(M,L))$ are equivalent.
		\end{prop}
		
		\begin{proof}
			We will identify both representations step by step. We already know that $\mathcal O_\lambda=G/C(T_1)$ and hence the base spaces are equal. Furthermore, both spaces are equal as complex manifolds, since the almost complex structure $\mathcal J$ on $G/C(T_1)$ and the polarization $\mathcal P$ on $\mathcal O_\lambda$ originate from the same subalgebra $\mathfrak b$ (cp. Proposition \ref{prop:coadkaehler}). In the proof of Theorem \ref{prop:intcondrev} we have seen that $L$ and $G\times_{\chi}\mathbb C$ can be identified as homogeneous line bundles (cp. Theorem \ref{theo:hlbisglued}). Therefore, similar to Proposition \ref{prop:holoispola}, it remains to show that the holomorphic sections are exactly the polarized sections of $L$. We can do that by relating both concepts to the Cauchy-Riemann equations.
			
			Knapp shows in \cite[p. 94]{knapp} the following statement. Let $f:M\to N$ be a smooth map between complex manifolds with almost complex structures $\mathcal J^M$ and $\mathcal J^N$. By definition $f$ is holomorphic in $p\in M$ if and only if
			$$\mathcal J^N_{f(p)}\circ f^\prime(p)=f^\prime(p)\circ J^M_p.$$
			Then the following is an equivalent condition. For every pair of holomorphic charts $(U,\varphi)$ at $p$ and $(V,\psi)$ at $f(p)$ with local coordinates $\varphi:=(x_1+iy_1,\dots,x_m+iy_m)$ and $\psi:=(u_1+iv_1,\dots,u_n+iv_n)$, where we assume $f$ to be given by $(u_j+iv_j)(x_1+iy_1,\dots,x_m+iy_m)$, $1\le j\le n$, the equations
			$$\frac{\partial u_k}{\partial x_j}(p)-\frac{\partial v_k}{\partial y_j}(p)=0,\quad \frac{\partial u_k}{\partial y_j}(p)+\frac{\partial v_k}{\partial x_j}(p)=0$$
			hold for $1\le j\le m$, $1\le k\le n$. Then Woodhouse asserts in \cite[§9.1]{woo} that the polarized sections are exactly the holomorphic sections, since $\mathcal P$ is locally spanned by $\frac{\partial}{\partial z_j}=\frac{1}{2}(\frac{\partial}{\partial x_j}-i \frac{\partial}{\partial y_j})$.
		\end{proof}
		
		\begin{theo}[Orbit Method]
			Let $G$ be a compact connected Lie group and $(\mathcal O_\lambda,\omega)$ a strongly integral coadjoint orbit. Then Construction \ref{cons:orbitmethod} yields an irreducible unitary representation $(\pi_{\mathcal O_\lambda},C^\infty_\mathcal P(M,L))$. Furthermore, all elements of $\hat{G}_u$ arise in that way.
		\end{theo}
		
		\begin{proof}
			By Corollary \ref{coro:fots} we can assume $\lambda\in\mathfrak t^*$ for a maximal torus $T$. Furthermore we use the result from Corollary \ref{coro:coadgfm} to assume $\mathcal O_\lambda\simeq C(T_1)$ for a torus $T_1\subseteq T$. It was shown in Theorem \ref{prop:intcondrev} that $\mathcal O_\lambda$ is integral if and only if $\lambda$ is the differential of a unitary character. We can then simultaneously apply Propositions \ref{prop:borelweildisguised} and \ref{prop:holopola} where we choose a $T_1$-admissible Weyl chamber and obtain two compatible root systems as in Theorem \ref{theo:eopocot}. Then we have that $(\pi_{\mathcal O_\lambda},C^\infty_\mathcal P(M,L))$ is equivalent to $(\pi_\lambda,\Gamma(G/T,G\times_\tau\mathbb C))$, where the latter is the representation from Construction \ref{cons:borel-weil}. Then the Borel-Weil Theorem yields that $(\pi_{\mathcal O_\lambda},C^\infty_\mathcal P(M,L))$ is irreducible. The unitary structure was shown in Proposition \ref{prop:coadunit}. In addition to the irreducibility we obtain from the Borel-Weil Theorem that every irreducible unitary representation is equivalent to the representation associated to some strongly integral coadjoint orbit.
		\end{proof}

	\label{derindex}
	\printindex


\newpage
	\begin{thebibliography}{\hspace{1.5 cm}} \label{biblio}
		
		\bibitem[AK71]{ausk} Auslander, L., Kostant, B. - \textit{Polarization and unitary representations of solvable Lie groups}, Inventiones Mathematicae 14, Springer, 1971
		\bibitem[Bla83]{bla} Blattner, R.J. - \textit{On geometric quantization},	Lecture Notes in Mathematics, Springer, Berlin, 1983
		\bibitem[CG90]{corgreen} Corwin, L.J., Greenleaf, F.P. - \textit{Representations of nilpotent Lie groups and their applications}, Cambridge University Press, Cambridge, 1990
		\bibitem[DK00]{dk} Duistermaat, J.J., Kolk, J.A.C. - \textit{Lie Groups}, Springer, Berlin, 2000
		\bibitem[Ech98]{spa} Echeverria, A. et al. - \textit{Mathematical foundations of geometric quantization}, EXTRACTA MATH. 13, 1998
		\bibitem[GH78]{gh} Griffiths, P., Harris, J. - \textit{Principles of algebraic geometry}, Wiley, New York, 1978
		\bibitem[Hal74]{halmos} Halmos, P.R. - \textit{Measure Theory}, Springer, New York, 1974
		\bibitem[He84]{helgason} Helgason, S. - \textit{Groups And Geometric Analysis}, Academic Press, Orlando, 1984
		\bibitem[Hil06]{hga} Hilgert, J. - \textit{Global Analysis}, lecture notes, Universität Paderborn, 2006
		\bibitem[Hil07]{hang} Hilgert, J. - \textit{Reproducing Kernels in Representation Theory}, in \textit{Symmetries in Complex Analysis}, B. Gilligan and G. Roos eds. Submitted.
		\bibitem[HN91]{hine} Hilgert, J., Neeb, K.H. - \textit{Lie-Gruppen und Lie-Algebren}, Vieweg, Braunschweig, 1991
		\bibitem[HNP94]{hhp} Hilgert, J., Neeb, K.H., Plank, W. - \textit{Symplectic convexity theorems and coadjoint orbits}, Compositio Mathematica, 94 no. 2, 1994, p. 129-180
		\bibitem[Kir04]{kir} Kirillov, A.A. - \textit{Lectures on the orbit method}, American Mathematical Society, Providence, R.I., 2004
		\bibitem[Kir76]{kireg} Kirillov, A.A. - \textit{Elements of the theory of representations}, Springer, Berlin, 1976
		\bibitem[Kir99]{kirmd} Kirillov, A.A. - \textit{Merits and Demerits of the Orbit Method}, Bulletin of the AMS, vol. 36, no. 3, 1999, p. 433-488
		\bibitem[Kna02]{knapp} Knapp, A.W. - \textit{Lie Groups Beyond an Introduction}, Birkhäuser, Boston, Mass., 2002
		\bibitem[KN69]{kn} Kobayashi, S., Nomizu, K. - \textit{Foundations of differential geometry Vol.2}, Interscience Publ., New York, 1969
		\bibitem[Kos70]{limaaa} Kostant, B. - \textit{Quantization and unitary representations}, Lecture Notes in Mathematics 170, Springer, Berlin, 1970
		\bibitem[Li80]{lips} Lipsman, R. - \textit{Orbit theory and harmonic analysis on Lie groups with co-compact nilradical}, J. Math. Pure and Appl. 59, 1980
		\bibitem[Nar85]{nara} Narasimhan, R. - \textit{Analysis on real and complex manifolds}, North-Holland, Amsterdam, 1985
		\bibitem[Pal56]{pal} Palais, R.S. - \textit{A global formulation of the Lie theory of transformation groups}, Memoirs of the American Mathematical Society 22, American Math. Soc., Providence, RI, 1968
		\bibitem[Pon05]{poncin} Poncin, N. - \textit{Méthodes géométriques et mathématique physique}, lecture notes, Université de Metz, 2005
		\bibitem[TW70]{tirwolf} Tirao, J.A., Wolf, J.A. - \textit{Homogeneous holomorphic vector bundles},  Indiana Univ. Math. J.  20  1970/1971 15-31
		\bibitem[Ve83]{vergne} Vergne, M. - \textit{Representations of Lie groups and the orbit method}, Emmy Noether in Bryn Mawr (Bryn Mawr, Pa., 1982), p. 59-101, Springer, New York, 1983
		\bibitem[Vog87]{vogur} Vogan, D.A. - \textit{Unitary representations of reductive Lie groups}, Princeton Univ. Pr., Princeton, NJ, 1987		
		\bibitem[Wal73]{walhomo} Wallach, N.R. - \textit{Harmonic analysis on homogeneous spaces}, Dekker, New York, 1973
		\bibitem[Wal77]{wal} Wallach, N.R. - \textit{Symplectic geometry and Fourier analysis}, Math Sci Press, Brookline, Mass., 1977
		\bibitem[War83]{war} Warner, F.W. - \textit{Foundations of differentiable manifolds and Lie groups}, Springer, New York, 1983
		\bibitem[Wel80]{wells} Wells, R.O. - \textit{Differential analysis on complex manifolds}, Springer, New York, 1980
		\bibitem[Woo92]{woo} Woodhouse, N.M.J. - \textit{Geometric quantization}, Clarendon Press, Oxford, 1992
	\end{thebibliography}
\end{document}